%% file: compsem-big.tex
\input{jbw-switch}
\DeclareSwitch{ictac}

\ictactrue

\documentclass[runningheads]{llncs}

\usepackage{amssymb}
\usepackage{graphicx}
\usepackage{url}
\usepackage{proof}

\newcommand{\MEM}[1]{\textup{$\mathrm{#1}$}}
\newcommand{\META}[1]{\mathit{#1}}

\newcommand{\METArel}[0]{\META{rel}}
\newcommand{\METAfun}[0]{\META{fun}}
\newcommand{\MEMdom}[1]{\MEM{dom}(#1)}
\newcommand{\MEMran}[1]{\MEM{ran}(#1)}

\newcommand{\mypair}[2]{\langle #1, #2 \rangle}
\newcommand{\func}[2]{#1 \f #2}
\newcommand{\fv}[1]{\MEM{fv}(#1)}
\newcommand{\dom}[1]{\MEM{dom}(#1)}

\def\be{\beta}
\def\wh{h}
\def\r{\restriction}
\def\u{\uparrow}

\def\l{\lambda}
\def\o{\omega}

\def\G{\Gamma}
\def\D{\Delta}
\def\B{\nabla}
\def\f{\rightarrow}
\def\fx{\leadsto}
\def\v{\vdash}
\def\<{\langle}
\def\>{\rangle}
\def\su{\blacktriangleright}
\def\F{\displaystyle\frac}

\def\myrhobe{\stackrel{\rho_\be}{\rightarrow}}

\def\leftrhobe{\stackrel{\rho_\be}{\leftarrow}}

\def\myrhobet{\stackrel{\rho_{\be\eta}}{\rightarrow}}

\def\leftrhobet{\stackrel{\rho_{\be\eta}}{\leftarrow}}

\def\myrhor{\stackrel{\rho_r}{\rightarrow}}
\def\myrrhor{\stackrel{\rho_r}{\rrightarrow}}
\def\myrrhor{\stackrel{\rho_r}{\rrightarrow}}
\def\lleftrhor{\stackrel{\rho_r}{\lleftarrow}}
\def\leftrhor{\stackrel{\rho_r}{\leftarrow}}
\def\deg{\mbox{d}}
\newcommand{\mbto}{\mathbin{\to}}
\newcommand{\legalenv}[1]{\mathrm{OK}(#1)}

\def\rrightarrow{\rightarrow \hspace{-.8em} \rightarrow}
\def\lleftarrow{\leftarrow \hspace{-.8em} \leftarrow}

\begin{document}

\title{A complete realisability semantics for intersection types and arbitrary expansion variables}

\titlerunning{A complete realisability semantics}

\author{Fairouz Kamareddine\inst{1}
\and Karim Nour\inst{2}
\and Vincent Rahli\inst{1}
\and J. B. Wells\inst{1}}

\institute{ULTRA Group (Useful Logics, Types, Rewriting, and their Automation), \email{http://www.macs.hw.ac.uk/ultra/}\and{Universit\'e de Savoie, Campus Scientifique, 73378 Le Bourget du Lac, France, \email{nour@univ-savoie.fr}}}

\maketitle

\begin{abstract}
  \emph{Expansion} was introduced at the end of the 1970s for
  calculating \emph{principal typings} for $\lambda$-terms in
  intersection type systems.  \emph{Expansion variables} (E-variables)
  were introduced at the end of the 1990s to simplify and help
  mechanise expansion.  Recently, E-variables have been further
  simplified and generalised to also allow calculating other type
  operators than just intersection.  There has been much work on
  semantics for intersection type systems, but only one such work on
  intersection type systems with E-variables.  That work established
  that building a semantics for E-variables is very challenging.
  Because it is unclear how to devise a space of meanings for
  E-variables, that work developed instead a space of meanings for
  types that is hierarchical in the sense of having many degrees
  (denoted by indexes). However, although the indexed calculus helped
  identify the serious problems of giving a semantics for expansion
  variables, the sound realisability semantics was only complete when
  one single E-variable is used and furthermore, the universal type
  $\omega$ was not allowed.  In this paper, we are able to overcome
  these challenges.  We develop a realisability semantics where we
  allow an arbitrary (possibly infinite) number of expansion variables
  and where $\omega$ is present.  We show the soundness and
  completeness of our proposed semantics.
\end{abstract}

\section{Introduction}

Expansion is a crucial part of a procedure for calculating
\emph{principal typings} and thus helps support compositional type
inference.
For example, the $\lambda$-term
$M=(\lambda{x}.x(\lambda{y}.yz))$ can be assigned the typing
\(
  \Phi_1 = \langle(z : a) \vdash (((a\mbto b)\mbto b)\mbto c)\mbto c\rangle
\),
which happens to be its principal typing.
The term $M$ can also be assigned the typing
\(
  \Phi_2 = \langle(z : a_1 \sqcap a_2) \vdash (((a_1 \mbto b_1)\mbto b_1) \sqcap ((a_2 \mbto b_2)\mbto b_2)\mbto c )\mbto c\rangle
\),
and an expansion operation can obtain $\Phi_2$ from $\Phi_1$.
Because the early definitions of expansion were complicated~\cite{Cop+Dez-Cia+Ven:HBC-1980},
E-variables were introduced in order to make the calculations
easier to mechanise and reason about.
For example, in System E~\cite{Car+Pol+Wel+Kfo:ESOP-2004}, the above typing
$\Phi_1$ is replaced by
\(
  \Phi_3 = \langle(z : e a) \vdash e((((a\mbto b)\mbto b)\mbto c)\mbto c)\rangle
\),
which differs from $\Phi_1$ by the insertion of the E-variable $e$ at
two places, and $\Phi_2$ can be obtained from $\Phi_3$ by substituting
for $e$ the \emph{expansion term}:\\
\(
  E = (a := a_1, b := b_1) \sqcap (a := a_2, b := b_2)
\).

Carlier and Wells~\cite{Car+Wel:ITRS-2004} have surveyed the history
of expansion and also E-variables.  Kamareddine, Nour, Rahli and
Wells~\cite{report} showed that E-variables pose serious challenges
for semantics.  In the list of open problems published in 1975 in
\cite{LNCS/1875}, it is suggested that an arrow type expresses
functionality.  Following this idea, a type's semantics is given as a
set of closed $\lambda$-terms with behaviour related to the
specification given by the type.  In many kinds of semantics, the
meaning of a type $T$ is calculated by an expression $[T]_\nu$ that
takes two parameters, the type $T$ and a valuation $\nu$ that assigns
to type variables the same kind of meanings that are assigned to
types. In that way, models based on term-models have been built for
intersection type systems
\cite{Hindley:ISOP-1982,krivine-lctm-1990,kamnour} where intersection
types (introduced to type more terms than in the Simply Typed Lambda
Calculus) are interpreted by set-theoretical intersection of meanings.
To extend this idea to types with E-variables, we need to devise some
space of possible meanings for E-variables.  Given that a type $e\,T$
can be turned by expansion into a new type $S_1(T) \sqcap S_2(T)$,
where $S_1$ and $S_2$ are arbitrary substitutions (or even arbitrary
further expansions), and that this can introduce an unbounded number
of new variables (both E-variables and regular type variables), the
situation is complicated.

This was the main motivation for \cite{report} to develop a space of
meanings for types that is hierarchical in the sense of having many
degrees. When assigning meanings to types, \cite{report} captured
accurately the intuition behind E-variables by ensuring that each use
of E-variables simply changes degrees and that each E-variable acts as
a kind of capsule that isolates parts of the $\lambda$-term being
analysed by the typing.

The semantic approach used in \cite{report} is realisability semantics
along the lines in Coquand~\cite{Coquand:TLCA-2005} and Kamareddine
and Nour~\cite{kamnour}.  Realisability allows showing
\emph{soundness} in the sense that the meaning of a type $T$ contains
all closed $\lambda$-terms that can be assigned $T$ as their result
type.  This has been shown useful in previous work for characterising
the behaviour of typed $\lambda$-terms \cite{krivine-lctm-1990}.  One
also wants to show the converse of soundness which is called
\emph{completeness} (see Hindley~\cite{Hin2,Hin3,Hindley:BSTT-1997}),
i.e., that every closed $\lambda$-term in the meaning of $T$ can be
assigned $T$ as its result type. Moreover, \cite{report} showed that
if more than one E-variable is used, the semantics is not
complete. Furthermore, the degrees used in \cite{report} made it
difficult to allow the universal type $\omega$ and this limited the
study to the $\lambda I$-calculus.  In this paper, we are able to
overcome these challenges.  We develop a realisability semantics where
we allow the full $\lambda$-calculus, an arbitrary (possibly infinite)
number of expansion variables and where $\omega$ is present, and we
show its soundness and completeness.  We do so by introducing an
indexed calculus as in \cite{report}. However here, our indices are
finite sequences of natural numbers rather than single natural
numbers.

In Section~\ref{secpure} we give the full $\lambda$-calculus indexed
with finite sequences of natural numbers and show the confluence of
$\beta$, $\beta \eta$ and weak head reduction on the indexed
$\lambda$-calculus.  In Section~\ref{sectypes} we introduce the type
system for the indexed $\lambda$-calculus (with the universal type
$\o$).  In this system, intersections and expansions cannot occur
directly to the right of an arrow.  In Section~\ref{srsec} we
establish that subject reduction holds for $\v$.  In
Section~\ref{sexpsec} we show that subject $\beta$-expansion holds for
$\v$ but that subject $\eta$-expansion fails.  In Section~\ref{redsec}
we introduce the realisability semantics and show its soundness for
$\v$.  In Section~\ref{complesec} we establish the completeness of $\v$
by introducing a special interpretation.  We conclude in
Section~\ref{concsec}.
\ifictac
Due to space limitations, we omit the details of the proofs. Full
proofs however can be found in the expanded version of this article
(currently at~\cite{kamnourrahliwells}) which will always be available
at the authors' web pages.
\else
Omitted proofs can be found in the appendix.
\fi

\section{The pure $\l^{{\cal L}_{\mathbb N}}$-calculus}
\label{secpure}

In this section we give the $\lambda$-calculus indexed with finite
sequences of natural numbers and show the confluence of $\beta$,
$\beta \eta$ and weak head reduction.

Let $n, m, i, j, k, l$ be metavariables which range over the set of
natural numbers ${\mathbb N} = \{0, 1, 2, \dots\}$.  We assume that if
a metavariable $v$ ranges over a set $s$ then $v_{i}$ and $v', v'',$
etc.\ also range over $s$.  A binary relation is a set of pairs. Let
$\METArel$ range over binary relations. We sometimes write $x\
\METArel\ y$ instead of $\mypair{x}{y} \in \METArel$. Let
$\MEMdom{\METArel} = \{x$ / $\mypair{x}{y} \in \METArel\}$ and
$\MEMran{\METArel} = \{y$ / $\mypair{x}{y} \in \METArel\}$. A function
is a binary relation $\METAfun$ such that if $\{\mypair{x}{y},
\mypair{x}{z}\} \subseteq \METAfun$ then $y = z$.  Let $\METAfun$
range over functions. Let $\func{s}{s'} = \{\METAfun$ /
$\MEMdom{\METAfun} \subseteq s \wedge \MEMran{\METAfun} \subseteq
s'\}$. We sometimes write $x : s$ instead of $x \in s$.

First, we introduce the set ${\cal L}_{\mathbb N}$ of indexes with an
order relation on indexes.
\begin{definition}
  \begin{enumerate}
  \item An index is a finite sequence of natural numbers $L = (n_i)_{1
      \leq i \leq l}$. We denote ${\cal L}_{\mathbb N}$ the set of
    indexes and $\oslash$ the empty sequence of natural numbers. We
    let $L, K, R$ range over ${\cal L}_{\mathbb N}$.
  \item If $L = (n_i)_{1 \leq i \leq l}$ and $m \in {\mathbb N}$, we
    use $m::L$ to denote the sequence $(r_i)_{1 \leq i \leq l+1}$
    where $r_1 = m$ and for all $i \in \{2, \dots, l+1\}$, $r_i =
    n_{i-1}$.\\ In particular, $k :: \oslash = (k)$.
  \item If $L = (n_i)_{1 \leq i \leq n}$ and $K = (m_i)_{1 \leq i \leq
      m}$, we use $L::K$ to denote the sequence $(r_i)_{1 \leq i \leq
      n+m}$ where for all $i \in \{1, \dots, n\}$, $r_i = n_i$ and for
    all $i \in \{n+1, \dots, n+m\}$, $r_i = m_{i-n}$. In particular,
    $L :: \oslash = \oslash :: L = L$.
  \item We define on ${\cal L}_{\mathbb N}$ a binary relation
    $\preceq$ by:

    $L_1 \preceq L_2$ (or $L_2 \succeq L_1$) if there exists $L_3 \in
    {\cal L}_{\mathbb N}$ such that $L_2 = L_1 :: L_3$.
  \end{enumerate}
\end{definition}

\begin{lemma}
  $\preceq$ is an order relation on ${\cal L}_{\mathbb N}$.
\end{lemma}

The next definition gives the syntax of the indexed calculus and the
notions of reduction.
\begin{definition}
  \begin{enumerate}
  \item Let ${\cal V}$ be a countably infinite set of variables.
    The set of terms ${\cal M}$, the set of free variables $\fv{M}$ of
    a term $M \in {\cal M}$, the degree function $\deg : {\cal M} \f
    {\cal L}_{\mathbb N}$ and the joinability $M \diamond N$ of terms
    $M$ and $N$ are defined by simultaneous induction as follows:
    \begin{itemize}
    \item If $x \in {\cal V}$ and $L \in {\cal L}_{\mathbb N}$, then
      $x^L \in {\cal M}$, $\fv{x^L} = \{x^L\}$ and $\deg(x^L) = L$.
    \item If $M, N \in {\cal M}$, $\deg(M) \preceq \deg(N)$ and $M
      \diamond N$ (see below), then $M \; N \in {\cal M}$, $\fv{MN} =
      \fv{M} \cup \fv{N}$ and $\deg(M \; N) = \deg(M)$.
    \item If $x \in {\cal V}$, $M \in {\cal M}$ and $L \succeq
      \deg(M)$, then $\l x^L.M \in {\cal M}$, $\fv{\l x^L. M} =
      \fv{M}\setminus \{x^L\}$ and $\deg(\l x^L.M) = \deg(M)$.
    \end{itemize}

  \item
    \begin{itemize}
    \item Let $M, N \in {\cal M}$.  We say that $M$ and $N$ are
      joinable and write $M \diamond N$ iff for all $x \in {\cal V}$,
      if $x^L \in \fv{M}$ and $x^K \in \fv{N}$, then $L = K$.
    \item If ${\cal X} \subseteq {\cal M}$ such that for all $M, N \in
      {\cal X}, M \diamond N$, we write, $\diamond{\cal X}$.
    \item If ${\cal X} \subseteq {\cal M}$ and $M \in {\cal M}$ such
      that for all $N \in {\cal X}, M \diamond N$, we write, $ M
      \diamond{\cal X}$.
    \end{itemize}
    The $\diamond$ property ensures that in any term $M$, variables
    have unique degrees.

    We assume the usual definition of subterms and the usual
    convention for parentheses and their omission (see
    Barendregt~\cite{Barendregt:LCSS-1984} and
    Krivine~\cite{krivine-lctm-1990}).  Note that every subterm of $M
    \in {\cal M}$ is also in ${\cal M}$.  We let $x, y, z, etc.$ range
    over ${\cal V}$ and $M, N, P$ range over ${\cal M}$ and use $=$
    for syntactic equality.

  \item The usual simultaneous substitution $M[(x^{L_i}_i:=N_i)_n]$ of
    $N_i \in {\cal M}$ for all free occurrences of $x^{L_i}_i$ in $M
    \in {\cal M}$ is only defined when $\diamond \{M\} \cup \{N_i$ /
    $i \in \{1, \dots, n\}\}$ and for all $i \in \{1, \dots, n\}$,
    $\deg(N_i) = L_i$. In a substitution, we sometimes write
    $x^{L_1}_1:=N_1, \dots, x^{L_n}_{n}:=N_n$ instead of
    $(x^{L_i}_{i}:=N_{i})_n$.  We sometimes write
    $M[(x^{L_i}_i:=N_i)_1$ as $M[x^{L_1}_1 := N_1]$.

  \item We take terms modulo $\alpha$-conversion given by:\\
    $\l x^L. M = \l y^L. (M[x^L:=y^L])$ where for all $L$, $y^L \not \in \fv{M}$.\\
    Moreover, we use the Barendregt convention (BC) where the names of
    bound variables differ from the free ones and where we rewrite
    terms so that not both $\l x^L$ and $\l x^K$ co-occur when $L \not
    = K$.

  \item A relation $\METArel$ on ${\cal M}$ is compatible iff for all
    $M, N, P \in {\cal M}$:
    \begin{itemize}
    \item If $M\ \METArel\ N$ and $\l x^L.M, \l x^L.M \in {\cal M}$
      then $(\l x^L.M)\ \METArel\ (\l x^L.N)$.
    \item If $M\ \METArel\ N$ and $MP, NP \in {\cal M}$ (resp.\ $PM,
      PN \in {\cal M}$), then $(MP)\ \METArel\ (NP)$ (resp.\ $(PM)\
      \METArel\ (PN)$).
    \end{itemize}

  \item The reduction relation $\rhd_\be $ on ${\cal M}$ is defined as
    the least compatible relation closed under the rule: $(\l x^L. M)N
    \rhd_\be M[x^L:=N]$ if $\deg(N) = L$

  \item The reduction relation $\rhd_\eta $ on ${\cal M}$ is defined
    as the least compatible relation closed under the rule: $\l
    x^L. (M \; x^L) \rhd_\eta M$ if $x^L \not \in \fv{M}$

  \item
    The weak head reduction $\rhd_{\wh} $ on ${\cal M}$  is defined by:\\
    $(\l x^L. M) N N_1 \dots N_n \rhd_{\wh} M[x^L:=N]N_1\dots N_n$
    where $n \geq 0$

  \item We let $\rhd_{\beta\eta} = \rhd_{\beta} \cup \rhd_{\eta}$.
    For $r \in \{\beta , \eta , \wh , \beta\eta\}$, we denote by
    $\rhd_r^*$ the reflexive and transitive closure of $\rhd_r$ and by
    $\simeq_{r}$ the equivalence relation induced by $\rhd_{r}^*$.
  \end{enumerate}
\end{definition}
The next theorem whose proof can be found in~\cite{kamnourrahliwells}
states that free variables and degrees are preserved by our notions of reduction.
\begin{theorem}
  \label{deg=}
  Let $M \in {\cal M}$ and $r \in \{\be, {\be\eta}, {\wh}\}$.
  \begin{enumerate}
  \item
    \label{deg=two}
    If $M \rhd_\eta^* N$ then $\fv{N} = \fv{M}$ and $\deg(M) =
    \deg(N)$.
  \item
    \label{deg=one}
    \label{deg=three}
    \label{deg=four}
    If $M \rhd_r^* N$ then $\fv{N} \subseteq \fv{M}$ and $\deg(M) =
    \deg(N)$.
  \end{enumerate}
\end{theorem}

As expansions change the degree of a term, indexes in a term need to
increase/decrease.
\begin{definition}
  Let $i \in {\mathbb N}$ and $M \in {\cal M}$.
  \begin{enumerate}
  \item We define $M^{+i}$ by:\\
    $\bullet (x^L)^{+i} = x^{i::L}\hspace{0.3in}$
    $\bullet (M_1 \; M_2)^{+i} = M_1^{+i} \; M_2^{+i}\hspace{0.3in}$
    $\bullet (\l x^L.M)^{+i} = \l x^{i::L}.M^{+i}$

    Let $M^{+\oslash} = M$ and $M^{+(i::L)} = (M^{+i})^{+L}$.
  \item If $\deg(M) = i::L$, we define $M^{-i}$ by:\\
    $\bullet (x^{i::K})^{-i} = x^{K}\hspace{0.3in}$
    $\bullet (M_1 \; M_2)^{-i} = M_1^{-i} \; M_2^{-i}\hspace{0.3in}$
    $\bullet (\l x^{i::K}.M)^{-i} = \l x^{K}.M^{-i}$

    Let $M^{-\oslash} = M$ and if $\deg(M) \succeq i::L$ then
    $M^{-(i::L)} = (M^{-i})^{-L}$.
  \item Let ${\cal X} \subseteq {\cal M}$. We write ${\cal X}^{+i}$
    for $\{M^{+i}$ / $M \in {\cal X}\}$.
  \end{enumerate}
\end{definition}

Normal forms are defined as usual.
\begin{definition}
  \begin{enumerate}
  \item $M \in {\cal M}$ is in $\be$-normal form ($\be\eta$-normal
    form, $\wh$-normal form resp.) if there is no $N\in {\cal M}$ such
    that $M \rhd_\be N$ ($M \rhd_{\be\eta} N$, $M \rhd_\wh N$ resp.).
  \item $M \in {\cal M}$ is $\be$-normalising ($\be\eta$-normalising,
    $\wh$-normalising resp.) if there is an $N \in {\cal M}$ such that
    $M \rhd^*_\be N$ ($M \rhd_{\be\eta} N$, $M \rhd_\wh N$ resp.) and
    $N$ is in $\be$-normal form ($\be\eta$-normal form, $\wh$-normal
    form resp.).
  \end{enumerate}
\end{definition}
The next theorem states that all of our notions of reduction are
confluent on our indexed calculus.  For a proof see~\cite{kamnourrahliwells}.
\begin{theorem}[Confluence]
  \label{confluenceofbetaeta}
  Let $M, M_1, M_2 \in {\cal M}$ and $r \in \{\be, {\be\eta}, {\wh}\}$.
  \begin{enumerate}
  \item
    \label{confitem}
    If $M \rhd_r^* M_1$ and $M \rhd_r^* M_2$, then there is $M'$ such
    that $M_1 \rhd_r^* M'$ and $M_2 \rhd_r^* M'$.
  \item
    \label{symprop}
    $M_1 \simeq_r M_2$ iff there is a term $M$ such that $M_1 \rhd^*_r
    M$ and $M_2 \rhd^*_r M$.
  \end{enumerate}
\end{theorem}

\section{Typing system}
\label{sectypes}

This paper studies a type system for the indexed $\lambda$-calculus
with the universal type $\o$. In this type system, in order to get
subject reduction and hence completeness, intersections and expansions
cannot occur directly to the right of an arrow (see ${\mathbb U}$
below).

The next two definitions introduce the type system.
\begin{definition}
  \begin{enumerate}
  \item Let $a$ range over a countably infinite set ${\cal A}$ of
    atomic types and let $e$ range over a countably infinite set
    ${\cal E} = \{\overline{e}_0,\overline{e}_1,...\}$ of expansion
    variables. We define sets of types ${\mathbb T}$ and ${\mathbb
      U}$, such that ${\mathbb T} \subseteq {\mathbb U}$, and a
    function $\deg : {\mathbb U} \f {\cal L}_{\mathbb N}$ by:
    \begin{itemize}
    \item If $a \in {\cal A}$, then $a \in {\mathbb T}$ and $\deg(a)
      =\oslash$.
    \item If $U \in {\mathbb U}$ and $T \in {\mathbb T}$, then $U \f T
      \in {\mathbb T}$ and $\deg(U \f T) =\oslash$.
    \item If $L \in {\cal L}_{\mathbb N}$, then $\o^L \in {\mathbb U}$
      and $\deg(\o^L) = L$.
    \item If $U_1,U_2 \in {\mathbb U}$ and $\deg(U_1) = \deg(U_2)$,
      then $U_1 \sqcap U_2 \in {\mathbb U}$ and $\deg(U_1 \sqcap U_2)
      = \deg(U_1) = \deg(U_2)$.
    \item $U \in {\mathbb U}$ and $\overline{e}_i \in {\cal E}$, then $\overline{e}_i U \in
      {\mathbb U}$ and $\deg(\overline{e}_i U) = i::\deg(U)$.
    \end{itemize}
    Note that $d$ remembers the number of the expansion variables
    $\overline{e}_i$ in order to keep a trace of these variables.

    We let $T$ range over ${\mathbb T}$, and $U, V, W$ range over
    ${\mathbb U}$.  We quotient types by taking $\sqcap$ to be
    commutative (i.e. $U_1 \sqcap U_2 = U_2 \sqcap U_1$), associative
    (i.e. $U_1 \sqcap (U_2 \sqcap U_3) = (U_1 \sqcap U_2) \sqcap U_3$)
    and idempotent (i.e. $U \sqcap U = U$), by assuming the
    distributivity of expansion variables over $\sqcap$ (i.e. $e(U_1
    \sqcap U_2) = eU_1 \sqcap eU_2$) and by having $\o^L$ as a neutral
    (i.e. $\o^L \sqcap U = U$). We denote $U_n\sqcap U_{n+1} \dots
    \sqcap U_{m}$ by $\sqcap_{i=n}^{m}U_i$ (when $n \leq m$). We also
    assume that for all $i \geq 0$ and $K \in {\cal L}_{\mathbb N}$,
    $\overline{e}_i\o^K = \o^{i::K}$.

  \item We denote $\overline{e}_{i_1}\dots \overline{e}_{i_n}$ by $\vec{e}_K$, where $K =
    (i_1,\dots,i_n)$ and $U_n\sqcap U_{n+1} \dots \sqcap U_{m}$ by
    $\sqcap_{i=n}^{m}U_i$ (when $n \leq m$).
  \end{enumerate}
\end{definition}

\begin{definition}
  \begin{enumerate}
  \item A type environment is a set $\{x^{L_1}_1:U_1, \dots,
    x^{L_n}_n:U_n\}$ such that for all $i, j \in \{1, \dots, n\}$, if
    $x^{L_i}_i = x^{L_j}_j$ then $U_i = U_j\}$. We let $Env$ be the
    set of environments,  use $\G, \Delta$ to
    range over $Env$ and write $()$ for the empty
    environment. We define $\dom{\G} = \{x^{L}$ / $x^{L} : U \in
    \G\}$.  If $\dom{\G_1} \cap \dom{\G_2} = \emptyset$, we write
    $\G_1, \G_2$ for $\G_1 \cup \G_2$. We write $\G, x^{L} : U$ for
    $\G, \{x^{L} : U\}$ and $x^{L} : U$ for $\{x^{L} : U\}$. We denote
    $x^{L_1}_1:U_1, \dots , x^{L_n}_n:U_n$ by $(x^{L_i}_i:U_i)_n$.
  \item If $M \in {\cal M}$ and $\fv{M} = \{x^{L_1}_1, \dots
    ,x^{L_n}_n\}$, we denote $env^\o_M$ the type environment
    $(x^{L_i}_i : \o^{L_i})_n$.
  \item We say that a type environment $\G$ is OK (and
    write  $\legalenv{\G}$) iff for all $x^{L}:U \in \G$,
    $\deg(U) = L$.
  \item Let $\G_1 = (x^{L_i}_i:U_i)_n,\G'_1$ and $\G_2 =
    (x^{L_i}_i:U'_i)_n,\G'_2$ such that $\dom{\G'_1} \cap \dom{\G'_2}
    = \emptyset$ and for all $i \in \{1, \dots, n\}$, $\deg(U_i) =
    \deg(U'_i)$. We denote $\G_1 \sqcap \G_2$ the type environment
    $(x^{L_i}_i:U_i \sqcap U'_i)_n, \G'_1, \G'_2$. Note that $\G_1
    \sqcap \G_2$ is a type environment, $\dom{\G_1 \sqcap \G_2} =
    \dom{\G_1} \cup \dom{\G_2}$ and that, on environments, $\sqcap$ is
    commutative, associative and idempotent.
  \item Let $\G = (x^{L_i}_i:U_i)_{1 \leq i \leq n}$
    We denote $\overline{e}_j\G = (x^{j::L_i}_i:\overline{e}_jU_i)_{1 \leq i \leq n}$.\\
    Note that $e\G$ is a type environment and $e(G_1\sqcap\G_2) =
    e\G_1\sqcap e\G_2$.
  \item We write $\G_1 \diamond \G_2$ iff $x^L \in \dom{\G_1}$ and $x^K
    \in \dom{\G_2}$ implies $K = L$.
  \item We follow \cite{Car+Wel:ITRS-2004} and write type judgements
    as $M: \<\G \v U\>$ instead of the traditional format of $\G \v M:
    U$, where $\vdash$ is our typing relation.  The typing rules of
    $\v$ are given on the left hand side of
    Figure~\ref{fig:typesystemrules}. In the last clause, the binary
    relation $\sqsubseteq$ is defined on ${\mathbb U}$ by the rules on
    the right hand side of Figure~\ref{fig:typesystemrules}.  We let
    $\Phi$ denote types in ${\mathbb U}$, or environments $\G$ or
    typings $\<\G\v U\>$.  When $\Phi \sqsubseteq \Phi'$, then $\Phi$
    and $\Phi'$ belong to the same set (${\mathbb
      U}$/environments/typings).

    {\footnotesize
      \begin{figure}[t]
        \begin{tabular}{|c|c|}
          \hline
          \begin{footnotesize}
            \begin{tabular}{c}
              \infer[(ax)]{x^{\oslash} : \<(x^{\oslash}:T) \v T\>}{}\\
              \\
              \infer[(\o)]{M : \<env^\o_M \v \o^{\deg(M)}\>}{}\\
              \\
              \infer[(\f_I)]{\l x^L. M : \<\G \v U \f T\>}{M : \<\G,(x^L:U) \v T\>}\\
              \\
              \infer[(\f'_I)]{\l x^L. M : \<\G \v \o^L \f T\>}{M : \<\G \v T\>\;\;\; x^L \not \in \dom{\G}}\\
              \\
              \infer[(\f_E)]{M_1 M_2 : \<\G_1
                \sqcap \G_2 \v T\>}{M_1 : \<\G_1 \v U \f T\> \;\;\;
                \hspace{0.01in} M_2 : \<\G_2 \v U\> \;\;\; \hspace{0.01in} \G_1 \diamond \G_2}\\
              \\
              \infer[(\sqcap_I)]{M : \<\G \v U_1 \sqcap U_2\>}{M: \<\G \v U_1\> \;\;\;\hspace{0.2in} M : \<\G \v U_2\>}\\
              \\
              \infer[(e)]{M^{+j} : \< \overline{e}_j\G \v \overline{e}_jU\>}{M : \<\G \v U\>}\\
              \\
              \infer[(\sqsubseteq)]{M : \<\G'\v U'\>}{M : \<\G\v U\> \;\;\;\hspace{0.2in} \<\G\v U\> \sqsubseteq
                \<\G'\v U'\>}
            \end{tabular}
          \end{footnotesize}
          &
          \begin{footnotesize}
            \begin{tabular}{c}
              \infer[(ref)]{\Phi \sqsubseteq \Phi}{}\\
              \\
              \infer[(tr)]{\Phi_1 \sqsubseteq \Phi_3}{\Phi_1 \sqsubseteq \Phi_2 &
                \Phi_2 \sqsubseteq \Phi_3}\\
              \\
              \infer[(\sqcap_E)]{U_1 \sqcap U_2 \sqsubseteq U_1}{d(U_1) = d(U_2)}\\
              \\
              \infer[(\sqcap)]{U_1 \sqcap U_2 \sqsubseteq V_1 \sqcap V_2}{U_1
                \sqsubseteq V_1 & U_2 \sqsubseteq V_2}\\
              \\
              \infer[(\f)]{U_1 \f T_1 \sqsubseteq U_2 \f T_2}{U_2 \sqsubseteq U_1
                & T_1 \sqsubseteq T_2}\\
              \\
              \infer[(\sqsubseteq_{e})]{eU_1 \sqsubseteq eU_2}{U_1 \sqsubseteq
                U_2}\\
              \\
              \infer[(\sqsubseteq_c)]{\G, y^L : U_1 \sqsubseteq \G, y^L :
                U_2}{U_1 \sqsubseteq U_2}\\
              \\
              \infer[(\sqsubseteq_{\<\>})]{\<\G_1 \v U_1\> \sqsubseteq \<\G_2
                \v U_2\>}{U_1 \sqsubseteq U_2 & \G_2 \sqsubseteq \G_1}
            \end{tabular}
          \end{footnotesize}\\
          \hline
        \end{tabular}
        \label{fig:typesystemrules}
        \caption{Typing rules / Subtyping rules}
      \end{figure}
    }
  \item
    If $L \in {\cal L}_{\mathbb N}$, $U \in {\mathbb U}$ and $\G =
    (x^{L_i}_i : U_i)_n$ is a type environment, we say that:
    \begin{itemize}
    \item $\deg(\G) \succeq L$ if and only if for all $i \in \{1,
      \dots, n\}$, $\deg(U_i) \succeq L$ and $L_i \succeq L$.
    \item $\deg(\<\G \v U\>) \succeq L$ if and only if $\deg(\G) \succeq L$
      and $\deg(U) \succeq L$.
    \end{itemize}
  \end{enumerate}
\end{definition}

To illustrate how our indexed type system works, we give an example:
\begin{example}
  Let $U = \overline{e}_3(\overline{e}_2(\overline{e}_1((\overline{e}_0b \f c) \f (\overline{e}_0(a \sqcap (a \f b)) \f c)) \f
  d) \f (((\overline{e}_2d \f a) \sqcap b) \f a))$ where $a, b, c, d \in {\cal A}$,

  $L_1 = 3 :: \oslash \preceq L_2 = 3 :: 2 :: \oslash \preceq L_3 = 3 :: 2
  :: 1 :: 0 :: \oslash$

  and

  $M = \l x^{L_2}.\l y^{L_1}.(y^{L_1} \, (x^{L_2} \, \l u^{L_3}.
  \l v^{L_3}. (u^{L_3} \, (v^{L_3} \, v^{L_3}))))$.

  We invite the reader to check that $M : \<() \v U\>$.
\end{example}

Just as we did for terms, we decrease the indexes of types,
environments and typings.
\begin{definition}
  \begin{enumerate}
  \item If $\deg(U) \succeq L$, then if $L = \oslash$ then $U^{-L} = U$
    else $L = i::K$ and we inductively define the type $U^{-L}$
    as follows:\\
    $(U_1 \sqcap U_2)^{-i::K} = U_1^{-i::K} \sqcap U_2^{-i::K}$
    \hspace{0.9in} $(\overline{e}_iU)^{-i::K} = U^{-K}$


    We write $U^{-i}$ instead of $U^{-(i)}$.
  \item If $\G = (x^{L_i}_i : U_i)_k$ and $\deg(\G) \succeq L$, then
    for all $i \in \{1, \dots, k\}$, $L_i = L :: L'_i$ and $\deg(U_i)
    \succeq L$ and we denote $\G^{-L} = (x^{L'_i} : U^{-L}_i)_k$.


    We write $\G^{-i}$ instead of $\G^{-(i)}$.
  \item If $U$ is a type and $\G$ is a type environment such that
    $\deg(\G) \succeq K$ and $\deg(U) \succeq K$, then we denote $(\<\G
    \v U\>)^{-K} = \<\G^{-K} \v U^{-K}\>$.
  \end{enumerate}
\end{definition}

The next lemma is informative about types and their degrees.
\begin{lemma}
  \label{degnew-goodeg}
  \begin{enumerate}
  \item If $T \in {\mathbb T}$, then $\deg(T) = \oslash$.
  \item Let $U \in {\mathbb U}$. If $\deg(U) = L = (n_i)_m$, then $U =
    \o^L$ or $U = \vec{e}_L \sqcap_{i=1}^p T_i$ where $p \geq 1$ and
    for all $i \in \{1, \dots, p\}$, $T_i \in {\mathbb T}$.
  \item
    \label{goodegone}
    Let $U_1 \sqsubseteq U_2$.
    \begin{enumerate}
    \item
      \label{goodegoneone}
      $\deg(U_1) = \deg(U_2)$.
    \item
      \label{goodegonetwo}
      If $U_1 = \o^K$ then $U_2 = \o^K$.
    \item
      \label{goodegonethree}
      If $U_1 = \vec{e}_KU$ then $U_2 = \vec{e}_KU'$ and $U
      \sqsubseteq U'$.
    \item
      \label{goodegonethree'}
      If $U_2 = \vec{e}_KU$ then $U_1 = \vec{e}_KU'$ and $U
      \sqsubseteq U'$.
    \item
      \label{goodegonefour}
      If $U_1 = \sqcap_{i=1}^p \vec{e}_K (U_i \f T_i)$ where $p \geq
      1$ then $U_2 = \o^K$ or $U_2 = \sqcap_{j=1}^q \vec{e}_K (U'_j \f
      T'_j)$ where $q \geq 1$ and for all $j \in \{1, \dots, q\}$,
      there exists $i \in \{1, \dots, p\}$ such that $U'_j \sqsubseteq
      U_i$ and $T_i \sqsubseteq T'_j$.
    \end{enumerate}
  \item
    \label{goodegtwo}
    If $U \in {\mathbb U}$ such that $\deg(U) = L$ then $U \sqsubseteq
    \o^L$.
  \item
    \label{goodegthree}
    If $U \sqsubseteq U_1' \sqcap U_2'$ then $U = U_1 \sqcap U_2$
    where $U_1 \sqsubseteq U_1'$ and $U_2 \sqsubseteq U_2'$.
  \item
    \label{goodegfour}
    If $\G \sqsubseteq \G_1' \sqcap \G_2'$ then $\G = \G_1 \sqcap
    \G_2$ where $\G_1 \sqsubseteq \G_1'$ and $\G_2 \sqsubseteq \G_2'$.
  \end{enumerate}
\end{lemma}

The next lemma says how ordering or the decreasing of indexes propagate to
environments.
\begin{lemma}
  \label{env-Phisub}
  \begin{enumerate}
  \item
    \label{lem:legalomegaenv}
    $\legalenv{env^\o_M}$.
  \item
    \label{Phisubone}
    If $\G \sqsubseteq \G'$, $U \sqsubseteq U'$ and $x^L \not \in
    \dom{\G}$ then $\G, (x^L:U)\sqsubseteq \G', (x^L:U')$.
  \item
    \label{Phisubtwo}
    $\G \sqsubseteq \G'$ iff $\G = (x^{L_i}_i : U_i)_n$, $\G' =
    (x^{L_i}_i : U'_i)_n$ and for every $1 \leq i \leq n$, $U_i
    \sqsubseteq U'_i$.
  \item
    \label{Phisubthree}
    $\<\G \v U\> \sqsubseteq \<\G' \v U'\>$ iff $\G' \sqsubseteq \G$
    and $U \sqsubseteq U'$.
  \item
    \label{Phisubfour}
    If $\dom{\G} = \fv{M}$ and $\legalenv{\G}$ then $\G \sqsubseteq
    env^\o_M$
  \item
    \label{Phisubfive}
    If $\G \diamond \D$ and $\deg(\G), \deg(\D) \succeq K$, then
    $\G^{-K} \diamond \D^{-K}$.
  \item
    \label{Phisubnew}
    If $U \sqsubseteq U'$ and $\deg(U) \succeq K$ then $U^{-K}
    \sqsubseteq U'^{-K}$.
  \item
    \label{Phisubnewnew}
    If $\G \sqsubseteq \G'$ and $\deg(\G) \succeq K$ then $\G^{-K}
    \sqsubseteq \G'^{-K}$.
  \item
    \label{lem:interlegalenv}
    If $\legalenv{\G_1}$, $\legalenv{\G_2}$ then $\legalenv{\G_1
      \sqcap \G_2}$.
  \item
    \label{lem:explegalenv}
    If $\legalenv{\G}$ then $\legalenv{e\G}$.
  \item
    \label{lem:deg+legalenvsub}
    If $\G_1 \sqsubseteq \G_2$ then ($\deg(\G_1) \succeq L$ iff
    $\deg(\G_2) \succeq L$) and ($\legalenv{\G_1}$ iff
    $\legalenv{\G_2}$).
  \end{enumerate}
\end{lemma}

The next lemma shows that we do not allow weakening in $\v$.
\begin{lemma}
  \label{structyping}
  \begin{enumerate}
  \item
    \label{structypingtwo}
    For every $\G$ and $M$ such that $\legalenv{\G}$ $\dom{\G} =
    \fv{M}$ and $\deg(M) = K$, we have $M: \<\G \v \o^K\>$.
  \item
    \label{structypingone}
    If $M: \<\G \v U\>$, then $\dom{\G} = \fv{M}$.
  \item
    \label{structypingC}
    If $M_1 : \<\G_1 \v U\>$ and $M_2 : \<\G_2 \v V\>$ then $\G_1
    \diamond \G_2$ iff $M_1 \diamond M_2$.
  \end{enumerate}
\end{lemma}

\begin{proof}
  \ref{structypingtwo}.\ By $\o$, $M: \<env^\o_M \v \o^K\>$. By
  Lemma~\ref{env-Phisub}.\ref{Phisubfour}, $\G \sqsubseteq env^\o_M$.
  Hence, by  $\sqsubseteq$ and $\sqsubseteq_{\<\>}$, $M: \<\G \v \o^K\>$.\\
  \ref{structypingone}.\ By induction on the derivation $M: \<\G \v
  U\>$. \\
  \ref{structypingC}.\ If) Let $x^L \in \dom{\G_1}$ and $x^K \in
  \dom{\G_2}$ then by Lemma~\ref{structyping}.\ref{structypingone},
  $x^L \in \fv{M_1}$ and $x^K \in \fv{M_2}$ so $\G_1 \diamond
  \G_2$. Only if) Let $x^L \in \fv{M_1}$ and $x^K \in \fv{M_2}$ then
  by Lemma~\ref{structyping}.\ref{structypingone}, $x^L \in
  \dom{\G_1}$ and $x^K \in \dom{\G_2}$ so $M_1 \diamond M_2$.
\hfill $\Box$
\end{proof}

The next theorem states that typings are well defined and that within a typing, degrees are well behaved.
\begin{theorem}
  \label{typing'}
    \begin{enumerate}
  \item
    \label{typing'zero}
    The typing relation $\v$ is well defined on ${\cal M}\times Env
    \times {\mathbb U}$.
  \item
    \label{typing'one}
    If $M : \<\G \v  U\>$ then $\legalenv{\G}$, and 
    $\deg(\G) \succeq \deg(U)=\deg(M)$.
  \item
    \label{typing'two'}
    If $M : \<\G \v  U\>$ and
$\deg(U) \succeq K$ then $M^{-K} : \<\G^{-K} \v U^{-K}\>$.
  \end{enumerate}
\end{theorem}

\begin{proof}
  We prove 1.\ and 2.\ simultaneously by induction on the derivation
  $M : \<\G \v U\>$. We prove 3.\ by induction on the derivation $M :
  \<\G \v U\>$. Full details can be found in~\cite{kamnourrahliwells}.
\hfill $\Box$
\end{proof}

Finally, here are two derivable typing rules that we will freely use
in the rest of the article.
\begin{remark}
  \label{remderiv}
  \begin{enumerate}
  \item
    \label{remderivone}
    The rule $\F{M: \<\G_1 \v U_1\> \;\;\;\hspace{0.2in} M : \<\G_2 \v
      U_2\>}{M : \<\G_1 \sqcap \G_2 \v U_1 \sqcap U_2\>} \; \;
    \sqcap'_I$ is derivable.
  \item
    \label{remderivtwo}
    The rule $\F{}{x^{\deg(U)} : \<(x^{\deg(U)} : U) \v U\>} \; \;
    ax'$ is derivable.
  \end{enumerate}
\end{remark}

\section{Subject reduction properties}
\label{srsec}

In this section we show that subject reduction holds for $\v$.  The
proof of subject reduction uses generation and substitution. Hence the
next two lemmas.

\begin{lemma}[Generation for $\v$]
  \label{newgen}
  \begin{enumerate}
  \item
    \label{newgenone}
    If $x^L : \<\G \v U \>$, then $\G = (x^L : V)$ and $V \sqsubseteq
    U$.
  \item
    \label{newgentwo}
    If $\l x^L. M : \<\G \v U\>$, $x^L \in \fv{M}$ and $\deg(U) = K$,
    then $U = \o^K$ or $U = \sqcap_{i=1}^p \vec{e}_K (V_i \f T_i)$
    where $p \geq 1$ and for all $i \in \{1, \dots, p\}$, $M : \<\G ,
    x^L : \vec{e}_K V_i \v \vec{e}_K T_i\>$.
  \item
    \label{newgenthree}
    If $\l x^L. M : \<\G \v U\>$, $x^L \not \in \fv{M}$ and $\deg(U) =
    K$, then $U = \o^K$ or $U = \sqcap_{i=1}^p \vec{e}_K (V_i \f T_i)$
    where $p \geq 1$ and for all $i \in \{1, \dots, p\}$, $M : \<\G \v
    \vec{e}_K T_i\>$.
  \item
    \label{newgenfour'}
    If $M \; x^L : \<\G , (x^L : U) \v T\>$ and $x^L \not \in \fv{M}$,
    then $M : \<\G \v U \f T\>$.
  \end{enumerate}
\end{lemma}

\begin{lemma}[Substitution for $\v$]
  \label{substlem}
  If $M:\<\G,x^L:U \v V\>$, $N :\<\Delta \v U\>$ and $M \diamond N$
  then $M[x^L:=N]:\<\G\sqcap \Delta \v V\>$.
\end{lemma}

Since $\v$ does not allow weakening, we need the next definition since
when a term is reduced, it may lose some of its free variables and
hence will need to be typed in a smaller environment.
\begin{definition}
  If $\G$ is a type environment and ${\cal U} \subseteq \dom{\G}$, then
  we write $\G\r_{\cal U}$ for the restriction of $\G$ on the
  variables of ${\cal U}$.  If ${\cal U} = \fv{M}$ for a term $M$, we
  write $\G \r_M$ instead of $\G \r_{\fv{M}}$.
\end{definition}

Now we are ready to prove the main result of this section:
\begin{theorem}[Subject reduction for $\v$]
  \label{subred}
  \label{subred_be}
  If $M : \<\G \v U\>$ and $M \rhd^*_{\beta\eta} N$, then $N :
  \<\G\r_N \v U\>$.
\end{theorem}

\begin{proof}
  By induction on the length of the derivation $M \rhd^*_{\beta\eta} N$.  
Case $M \rhd_{\beta\eta} N$ is by induction on the derivation
$M : \<\G \v_3 U\>$.  \hfill $\Box$
\end{proof}

\begin{corollary}
  \label{subredC}
  \begin{enumerate}
  \item
    \label{subredC_b}
    If $M : \<\G \v U\>$ and $M \rhd^*_{\beta} N$, then $N : \<\G\r_N
    \v U\>$.
  \item
    \label{subredC_f}
    If $M : \<\G \v U\>$ and $M \rhd^*_{\wh} N$, then $N : \<\G\r_N \v
    U\>$.
  \end{enumerate}
\end{corollary}

\section{Subject expansion properties}
\label{sexpsec}

In this section we show that subject $\beta$-expansion holds for
$\v$ but that  subject $\eta$-expansion fails.

The next lemma is needed for expansion.
\begin{lemma}
  \label{exp1}
  If $M[x^L:= N] : \<\G \v U\>$ and $x^L \in \fv{M}$ then there exist
  a type $V$ and two type environments $\G_1, \G_2$
  such that:\\
  $M : \<\G_1, x^L : V \v U\>\hspace{0.3in}$ $N : \<\G_2 \v
  V\>\hspace{0.3in}$ $\G = \G_1 \sqcap \G_2$
\end{lemma}

Since more free variables might appear in the $\beta$-expansion of a
term, the next definition gives a possible enlargement of an
environment.
\begin{definition}
  Let $m \geq n$, $\G = (x^{L_i}_i : U_i)_n$ and ${\cal U} =
  \{x^{L_1}_1,...,x^{L_m}_m\}$. We write $\G {\u^{\cal U}}$ for
  $x^{L_1}_1 : U_1,...,x^{L_n}_n : U_n,x^{L_{n+1}}_{n+1} :
  \o^{L_{n+1}},...,x^{L_m}_m : \o^{L_m}$.  Note that $\G{\u^{\cal U}}$
  is a type environment.  If $\dom{\G} \subseteq \fv{M}$, we write $\G
  {\u^M}$ instead of $\G {\u^{\fv{M}}}$.
\end{definition}

We are now ready to establish that subject expansion holds for $\beta$
(next theorem) and that it fails for $\eta$ (Lemma~\ref{etaexpfail}).
\begin{theorem}[Subject expansion for $\beta$]
  \label{finalbetaexp}
  If $N : \<\G \v U\>$ and $M \rhd^*_{\beta} N$, then $M : \<\G
  {\u^{M}} \v U\>$.
\end{theorem}
\begin{proof}
  By induction on the length of the derivation $M \rhd^*_{\beta} N$
  using the fact that if $\fv{P} \subseteq \fv{Q}$, then $(\G {\u^{P}})
  {\u^{Q}} = \G {\u^{Q}}$.   \hfill $\Box$
\end{proof}
\begin{corollary}
  \label{finalbetaexpC}
  If $N : \<\G \v U\>$ and $M \rhd^*_{\wh} N$, then $M : \<\G {\u^{M}}
  \v U\>$.
\end{corollary}

\begin{lemma}[Subject expansion fails for $\eta$]
  \label{etaexpfail}
  Let $a$ be an element of ${\cal A}$.  We have:
  \begin{enumerate}
  \item $\l y^\oslash. \l x^\oslash. y^\oslash x^\oslash \rhd_\eta \l
    y^\oslash. y^\oslash$
  \item $\l y^\oslash. y^\oslash : \<() \v a \f a\>$.
  \item It is not possible that

    $\l y^\oslash. \l x^\oslash. y^\oslash x^\oslash : \<() \v a \f a\>$.\\
    Hence, the subject $\eta$-expansion lemmas fail for $\v$.
  \end{enumerate}
\end{lemma}

\begin{proof}
  1. and 2. are easy.  For 3., assume $\l y^\oslash. \l
  x^\oslash. y^\oslash x^\oslash : \<() \v a \f a\>$. \\
  By Lemma~\ref{newgen}.\ref{newgentwo}, $\l x^\oslash. y^\oslash
  x^\oslash : \<(y : a) \v \f a\>$. Again, by
  Lemma~\ref{newgen}.\ref{newgentwo}, $a = \o^\oslash$ or there exists
  $n \geq 1$ such that $a = \sqcap^n_{i=1} (U_i \f T_i)$, absurd.
\hfill $\Box$
\end{proof}

\section{The realisability semantics}
\label{redsec}

In this section we introduce the realisability semantics and show its
soundness for $\v$.

Crucial to a realisability semantics is the notion of a saturated set:
\begin{definition}
  \label{setdefs}
  Let ${\cal X},{\cal Y} \subseteq {\cal M}$.
  \begin{enumerate}
  \item
    We use ${\cal P}({\cal X})$ to denote the powerset of ${\cal X}$,
    i.e.\ $\{{\cal Y} \; / \; {\cal Y} \subseteq {\cal X}\}$.
  \item
    We define ${\cal X}^{+i} = \{M^{+i}$ / $M \in {\cal X}\}$.
  \item
    We define ${\cal X} \fx {\cal Y} = \{M \in {\cal M}$ / $M \;
    N \in {\cal Y}$ for all $N \in {\cal X}$ such that $M \diamond
    N\}$.
  \item We say that ${\cal X} \wr {\cal Y}$ iff for all $M \in {\cal
      X} \fx {\cal Y}$, there exists $N \in {\cal X}$ such that $M
    \diamond N$.
  \item
    For $r \in \{\be, \be\eta, \wh\}$, we say that ${\cal X}$ is
    $r$-saturated if whenever $M \rhd_r^* N$ and $N \in {\cal X}$, then
    $M \in {\cal X}$.
  \end{enumerate}
\end{definition}

Saturation is closed under intersection, lifting and arrows:
\begin{lemma}
  \label{fx+}
  \begin{enumerate}
  \item
    \label{fx+one}
    $({\cal X} \cap {\cal Y})^{+i} = {\cal X}^{+i} \cap {\cal
      Y}^{+i}$.
  \item
    \label{fx+two}
    If ${\cal X},{\cal Y}$ are $r$-saturated sets, then ${\cal X} \cap
    {\cal Y}$ is $r$-saturated.
  \item
    \label{fx+three}
    If ${\cal X}$ is $r$-saturated, then ${\cal X}^{+i}$ is
    $r$-saturated.
  \item
    \label{fx+four}
    If ${\cal Y}$ is $r$-saturated, then, for every set ${\cal X}$,
    ${\cal X} \fx {\cal Y}$ is $r$-saturated.
  \item
    \label{fx+five}
    $({\cal X} \fx {\cal Y})^{+i} \subseteq {\cal X}^{+i} \fx {\cal
      Y}^{+i}$.
  \item
    \label{fx+six}
    If ${\cal X}^{+i} \wr {\cal Y}^{+i}$, then ${\cal X}^{+i} \fx
    {\cal Y}^{+i} \subseteq ({\cal X} \fx {\cal Y})^{+i}$.
  \end{enumerate}
\end{lemma}

We now give the basic step in our realisability semantics: the
interpretations and meanings of types.
\begin{definition}
  \label{intdef}
  Let ${\cal V}_1$, ${\cal V}_2$ be countably infinite, ${\cal V}_1
  \cap {\cal V}_2 = \emptyset$ and ${\cal V} = {\cal V}_1 \cup {\cal
    V}_2$.
  \begin{enumerate}
  \item Let $L \in {\cal L}_{\mathbb N}$. We define ${\cal M}^L = \{M
    \in {\cal M}$ / $\deg(M) = L\}$.
  \item Let $x \in {\cal V}_1$. We define ${\cal N}_x^L = \{x^L \;
    N_1...N_k \in {\cal M}$ / $k \geq 0\}$.
  \item Let $r \in \{\be, \be\eta, \wh\}$. An $r$-interpretation
    ${\cal I}: {\cal A} \mapsto {\cal P} ({\cal M}^{\oslash})$ is a
    function
    such that for all $a \in {\cal A}$:\\
    $\bullet \; {\cal I}(a)$ is $r$-saturated \hspace{0.5in} and
    \hspace{0.5in} $\bullet \; \forall x \in {\cal V}_1.\ {\cal
      N}_x^{\oslash} \subseteq {\cal I}(a)$.

    We extend an $r$-interpretation ${\cal I}$ to ${\mathbb U}$
    as follows:\\
    $\bullet \; {\cal I}(\o^L) = {\cal M}^L$ \hspace{1.4in} $\bullet
    \;
    {\cal I}(\overline{e}_iU) = {\cal I}(U)^{+i}$\\
    $\bullet \; {\cal I}(U_1 \sqcap U_2) = {\cal I}(U_1) \cap {\cal
      I}(U_2)$ \hspace{0.5in} $\bullet \; {\cal I}(U \f T) = {\cal
      I}(U)
    \fx {\cal I}(T)$\\
    Let $r\mbox{-int}= \{{\cal I} \; / \; {\cal I} \mbox{ is an
      $r$-interpretation}\}$.
  \item Let $U \in {\mathbb U}$ and $r \in \{\be, \be\eta,
    \wh\}$. Define $[U]_r$, the $r$-interpretation of $U$ by:

    $[U]_r = \{M \in {\cal M} \; / \; M \; is \; closed \; and \; M
    \in \bigcap_{{\cal I}\in r\mbox{-int}} {\cal I}(U)\}$
  \end{enumerate}
\end{definition}

\begin{lemma}
  \label{interpret-intsub}
  Let $r \in \{\be, \be\eta, \wh\}$.
  \begin{enumerate}
  \item
    \label{interpret}
    \begin{enumerate}
    \item
      \label{interpretone}
      For any $U \in {\mathbb U}$ and ${\cal I}\in r\mbox{-int}$, we
      have ${\cal I}(U)$ is $r$-saturated.
    \item
      \label{interprettwo}
      If $\deg(U) = L$ and ${\cal I}\in r\mbox{-int}$, then for all $x
      \in {\cal V}_1$, ${\cal N}_x^L \subseteq {\cal I}(U) \subseteq
      {\cal M}^L$.
    \end{enumerate}
  \item
    \label{intsub}
    Let $r \in \{\be, \be\eta, \wh\}$. If ${\cal I}\in r\mbox{-int}$
    and $U \sqsubseteq V$, then ${\cal I}(U) \subseteq {\cal I}(V)$.
  \end{enumerate}
\end{lemma}

Here is the soundness lemma.
\begin{lemma}[Soundness]
  \label{adeq}
  Let $r \in \{\be, \be\eta, \wh\}$, $M : \<(x^{L_j}_j:U_j)_n \v U\>$,
  ${\cal I}\in r\mbox{-int}$ and for all $j \in \{1, \dots, n\}$, $N_j
  \in {\cal I}(U_j)$. If $M[(x^{L_j}_j:=N_j)_n] \in {\cal M}$ then
  $M[(x^{L_j}_j:=N_j)_n] \in {\cal I}(U)$.
\end{lemma}
\begin{proof}
  By induction on the derivation $M : \<(x^{L_j}_j:U_j)_n \v U\>$.  \hfill $\Box$
\end{proof}
\begin{corollary}
  \label{finalint}
  Let $r \in \{\be, \be\eta, \wh\}$. If $M : \<() \v U\>$, then $M \in
  [U]_r$. \hfill $\Box$
\end{corollary}

\begin{proof}
  By Lemma~\ref{adeq}, $M \in {\cal I}(U)$ for any $r$-interpretation
  ${\cal I}$. By Lemma~\ref{structyping}.\ref{structypingone}, $\fv{M}
  = \dom{()}= \emptyset$ and hence $M$ is closed. Therefore, $M \in
  [U]_r$.
\hfill $\Box$
\end{proof}

\begin{lemma}[The meaning of types is closed under type operations]
  \label{[]} \mbox{} \hfill\\
  Let $r \in \{\be, \be\eta, \wh\}$. On ${\mathbb U}$, the following
  hold:
  \begin{enumerate}
  \item
    \label{[]one}
    $[\overline{e}_iU]_r = [U]_r^{+i}$
  \item
    \label{[]two}
    $[U \sqcap V]_r = [U]_r \cap [V]_r$
  \item
    \label{[]three}
    If ${\cal I} \in r\mbox{-int}$ and
    $U, V \in {\mathbb U}$, then ${\cal I}(U) \wr {\cal I}(V).$
  \end{enumerate}
\end{lemma}

\begin{proof}
  \ref{[]one}.\ and \ref{[]two}.\ are easy. \ref{[]three}.\ Let
  $\deg(U) = K$, $M \in {\cal I}(U) \fx {\cal I}(V)$ and $x \in {\cal
    V}_1$ such that for all $L$, $x^L \not \in \fv{M}$, then $M
  \diamond x^K$ and by
  lemma~\ref{interpret-intsub}.\ref{interprettwo}, $x^K \in {\cal
    I}(U)$.
\hfill $\Box$
\end{proof}

The next definition and lemma put the realisability semantics in use.
\begin{definition}[Examples]
  Let $a,b \in {\cal A}$ where $a \not = b$.  We define:
  \begin{itemize}
  \item $Id_{0} = a \f a$, $Id_{1} = \overline{e}_1(a \f a)$ and
    $Id'_{1} = \overline{e}_1a \f \overline{e}_1a$.
  \item $D = (a \sqcap (a \f b)) \f b$.
  \item $Nat_{0} = (a \f a) \f (a \f a)$,
    $Nat_{1} = \overline{e}_1((a \f a) \f (a \f a))$,\\
    and $Nat'_{0}=(\overline{e}_1a \f a) \f (\overline{e}_1a \f a)$.
  \end{itemize}
  Moreover, if $M,N$ are terms and $n \in {\mathbb N}$, we define
  $(M)^n \; N$ by induction on $n$: $(M)^0 \; N = N$ and $(M)^{m+1} \;
  N = M \; ((M)^m \; N)$.
\end{definition}

\begin{lemma}
  \label{exlem}
  \label{newexlem}
  \begin{enumerate}
  \item
    \label{exone}
    \label{newexone}
    $[Id_0]_\be = \{ M \in {\cal M}^\oslash$ / $M$ is closed and $M
    \rhd_\be^* \l y^\oslash. y^\oslash\}$.
  \item
    \label{extwo}
    \label{newextwo}
    $[Id_1]_\be = [Id'_1]_\be = \{M \in {\cal M}^{(1)}$ / $M$ is
    closed and $M \rhd_\be^* \l y^{(1)}.y^{(1)}\}$. (Note that $Id'_1
    \not \in {\mathbb U}$.)
  \item
    \label{exthree}
    \label{newexthree}
    $[D]_\be = \{M \in {\cal M}^\oslash$ / $M$ is closed and $M
    \rhd_\be^* \l y^\oslash.y^\oslash y^\oslash\}$.
  \item
    \label{exfour}
    \label{newexfour}
    $[Nat_0]_\be = \{M \in {\cal M}^\oslash$ / $M$ is closed and $M
    \rhd_\be^*\l f^\oslash.f^\oslash$ or $M \rhd_\be^* \l f^\oslash.\l
    y^\oslash.(f^\oslash)^n y^\oslash$ where $n \geq 1\}$.
  \item
    \label{exfive}
    \label{newexfive}
    $[Nat_1]_\be = \{M \in {\cal M}^{(1)}$ / $M$ is closed and $M
    \rhd_\be^*\l f^{(1)}.f^{(1)}$ or $M \rhd_\be^* \l f^{(1)}.\l
    x^{(1)}.(f^{(1)})^n y^{(1)}$ where $n \geq 1\}$. (Note that
    $Nat'_1 \not \in {\mathbb U}$.)
  \item
    \label{exsix}
    \label{newexsix}
    $[Nat'_0]_\be = \{M \in {\cal M}^\oslash$ / $M$ is closed and $M
    \rhd_\be^*\l f^\oslash.f^\oslash$ or $M \rhd_\be^* \l f^\oslash.\l
    y^{(1)}.f^\oslash y^{(1)}\}$.
  \end{enumerate}
\end{lemma}

\section{The completeness theorem}
\label{complesec}

In this section we set out the machinery and prove that completeness
holds for $\v$.

We need the following partition of the set of variables $\{y^L / y \in
{\cal V}_2\}$.
\begin{definition}
  \begin{enumerate}
  \item Let $L \in {\cal L}_{\mathbb N}$. We define ${\mathbb U}^L =
    \{ U \in {\mathbb U} / \deg(U) = L\}$ and ${\cal V}^L = \{x^L / x
    \in {\cal V}_2\}$.
  \item Let $U \in {\mathbb U}$. We inductively define a set of
    variables ${\mathbb V}_U$ as follows:
    \begin{itemize}
    \item If $\deg(U) = \oslash$ then:
      \begin{itemize}
      \item ${\mathbb V}_U$ is an infinite set of variables of degree
        $\oslash$.
      \item If $U \neq V$ and $\deg(U) = \deg(V) = \oslash$, then
        ${\mathbb V}_U \cap {\mathbb V}_V = \emptyset$.
      \item $\bigcup_{U \in {\mathbb U}^{\oslash}}{\mathbb V}_U =
        {\cal V}^{\oslash}$.
      \end{itemize}
    \item If $\deg(U) = L$, then we put ${\mathbb V}_U = \{ y^L$ /
      $y^{\oslash} \in {\mathbb V}_{U^{-L}}\}$.
    \end{itemize}
  \end{enumerate}
\end{definition}

\begin{lemma}
  \label{codeV}
  \begin{enumerate}
  \item If $\deg(U),\deg(V) \succeq L$ and $U^{-L} = V^{-L}$, then $U
    = V$.
  \item If $\deg(U) = L$, then ${\mathbb V}_U$ is an infinite subset
    of ${\cal V}^L$.
  \item If $U \neq V$ and $\deg(U) = \deg(V) = L$, then ${\mathbb V}_U
    \cap {\mathbb V}_V = \emptyset$.
  \item $\bigcup_{U \in {\mathbb U}^L}{\mathbb V}_U = {\cal V}^L$.
  \item If $y^L \in {\mathbb V}_U$, then $y^{i::L} \in {\mathbb
      V}_{\overline{e}_iU}$.
  \item If $y^{i::L} \in {\mathbb V}_U$, then $y^L \in {\mathbb
      V}_{U^{-i}}$.
  \end{enumerate}
\end{lemma}

\begin{proof}
  1.\ If $L = (n_i)_m$, we have $U =
  \overline{e}_{n_1}\dots\overline{e}_{n_m}U'$ and $V =
  \overline{e}_{n_1}\dots\overline{e}_{n_m}V'$. Then $U^{-L} = U'$,
  $V^{-L} = V'$ and $U' = V'$. Thus $U = V$. 2.\, 3.\ and 4.\ By
  induction on $L$ and using 1.\ 5.\ Because $(\overline{e}_i U)^{-i} = U$. 6.\
  By definition.
\hfill $\Box$
\end{proof}

Our partition of the set ${\cal V}_2$ as above will enable us to give
in the next definition useful infinite sets which will contain type
environments that will play a crucial role in one particular type
interpretation.
\begin{definition}
  \label{infsetsenv}
  \begin{enumerate}
  \item Let $L \in {\cal L}_{\mathbb N}$.  We denote ${\mathbb
      G}^L=\{(y^L:U)$ / $U\in{\mathbb U}^L$ and $y^L \in {\mathbb V}_U
    \}$ and ${\mathbb H}^L=\bigcup_{K \succeq L}{\mathbb G}^K$. Note
    that ${\mathbb G}^L$ and ${\mathbb H}^L$ are not type environments
    because they are infinite sets.
  \item Let $L \in {\cal L}_{\mathbb N}$, $M\in{\cal M}$ and
    $U\in{\mathbb U}$, we write:
    \begin{itemize}
    \item $M : \< {\mathbb H}^L \v U\>$ if there is a type environment
      $\G \subset {\mathbb H}^L$ where $M : \< \G \v U\>$
    \item $M : \< {\mathbb H}^L \v^* U\>$ if $M \rhd^*_{\beta\eta} N$
      and $N : \< {\mathbb H}^L \v U\>$
    \end{itemize}
  \end{enumerate}
\end{definition}

\begin{lemma}
  \label{GnHn}
  \begin{enumerate}
  \item If $\G \subset {\mathbb H}^L$ then $\legalenv{\G}$.
  \item If $\G \subset {\mathbb H}^L$ then $\overline{e}_i\G \subset {\mathbb
      H}^{i::L}$.
  \item If $\G \subset {\mathbb H}^{i::L}$ then $\G^{-i} \subset
    {\mathbb H}^L$.
  \item If $\G_1 \subset {\mathbb H}^L$, $\G_2 \subset {\mathbb H}^K$
    and $L \preceq K$ then $\G_1 \sqcap \G_2 \subset {\mathbb H}^L$.
  \end{enumerate}
\end{lemma}

\begin{proof}
  1.\ Let $x^{K} : U \in \G$ then $U \in {\mathbb U}^{K}$ and so
  $\deg(U) = K$.
  2.\ and 3.\ are by lemma \ref{codeV}.
  4.\ First note that by 1., $\G_1 \sqcap \G_2$ is well defined.
  ${\mathbb H}^K \subseteq {\mathbb H}^L$.  Let $(x^R : U_1 \sqcap
  U_2) \in \G_1 \sqcap \G_2$ where $(x^R : U_1) \in \G_1 \subset
  {\mathbb H}^L$ and $(x^R : U_2) \in \G_2 \subset {\mathbb H}^K
  \subseteq {\mathbb H}^L$, then $\deg(U_1) = \deg(U_2) = R$ and $x^R
  \in {\mathbb V}_{U_1} \cap {\mathbb V}_{U_2}$. Hence, by
  lemma~\ref{codeV}, $U_1 = U_2$ and $\G_1 \sqcap \G_2 = \G_1 \cup
  \G_2 \subset {\mathbb H}^L$.
\hfill $\Box$
\end{proof}

For every $L \in {\cal L}_{\mathbb N}$, we define the set of terms of
degree $L$ which contain some free variable $x^K$ where $x \in {\cal
  V}_1$ and $K \succeq L $.
\begin{definition}
  For every $L \in {\cal L}_{\mathbb N}$, let ${\cal O}^L = \{M \in
  {\cal M}^L$ / $x^K \in \fv{M}$, $x \in {\cal V}_1$ and $K \succeq L
  \}$. It is easy to see that, for every $L \in {\cal L}_{\mathbb N}$
  and $x \in {\cal V}_1$, ${\cal N}^L_x \subseteq {\cal O}^L$.
\end{definition}

\begin{lemma}
  \label{Vx}
  \begin{enumerate}
  \item $({\cal O}^L)^{+i} = {\cal O}^{i::L}$.
  \item If $y \in {\cal V}_2$ and $(M y^K) \in {\cal O}^L$, then $M
    \in {\cal O}^L$
  \item If $M \in {\cal O}^L$, $M \diamond N$ and $L \preceq K
    =\deg(N)$, then $M N \in {\cal O}^L$.
  \item If $\deg(M) = L$, $L \preceq K$, $M \diamond N$ and $N \in
    {\cal O}^K$, then $M N \in {\cal O}^L$.
  \end{enumerate}
\end{lemma}

The crucial interpretation ${\mathbb I}$ for the proof of completeness
is given as follows:
\begin{definition}
  \begin{enumerate}
  \item Let ${\mathbb I}_{\be\eta}$ be the $\be\eta$-interpretation
    defined by: for all type variables $a$, ${\mathbb I}_{\be\eta}(a)
    = {\cal O}^{\oslash} \cup \{M \in {\cal M}^{\oslash}$ / $M : \<
    {\mathbb H}^{\oslash} \v^* a\> \}$.
  \item Let ${\mathbb I}_{\be}$ be the $\be$-interpretation defined
    by: for all type variables $a$, ${\mathbb I}_{\be}(a) = {\cal
      O}^{\oslash} \cup \{M \in {\cal M}^{\oslash}$ / $M : \< {\mathbb
      H}^{\oslash} \v a\> \}$.
  \item Let ${\mathbb I}_{\wh}$ be the $\wh$-interpretation defined
    by: for all type variables $a$, ${\mathbb I}_{\wh}(a) = {\cal
      O}^{\oslash} \cup \{M \in {\cal M}^{\oslash}$ / $M : \< {\mathbb
      H}^{\oslash} \v a\> \}$.
  \end{enumerate}
\end{definition}

The next crucial lemma shows that ${\mathbb I}$ is an interpretation
and that the interpretation of a type of order $L$ contains terms of
order $L$ which are typable in these special environments which are
parts of the infinite sets of Definition~\ref{infsetsenv}.
\begin{lemma}
  \label{comp}
  Let $r \in \{\be\eta, \be, \wh\}$ and $r' \in \{\be, \wh\}$
  \begin{enumerate}
  \item If ${\mathbb I}_{r}\in r\mbox{-int}$ and $a \in {\cal A}$ then
    ${\mathbb I}_{r}(a)$ is $r$-saturated and for all $x \in {\cal
      V}_1, {\cal N}_x^{\oslash} \subseteq {\mathbb I}_{r}(a)$.
  \item If $U\in{\mathbb U}$ and $\deg(U) = L$, then ${\mathbb
      I}_{\be\eta}(U) = {\cal O}^L \cup \{M \in {\cal M}^L$ / $M :
    \<{\mathbb H}^L \v^* U\> \}$.
  \item If $U\in{\mathbb U}$ and $\deg(U) = L$, then ${\mathbb
      I}_{r'}(U) = {\cal O}^L \cup \{M \in {\cal M}^L$ / $M :
    \<{\mathbb H}^L \v U\> \}$.
  \end{enumerate}
\end{lemma}

\ifictac
\else
\begin{proof} 1.\ We do two cases:\\
  Case $r = \beta \eta$.
  It is easy to see that $\forall x \in {\cal V}_1, {\cal
    N}_x^\oslash \subseteq {\cal
    O}^\oslash \subseteq {\mathbb I}_{\be\eta}(a)$.
  Now we show that ${\mathbb I}_{\be\eta}(a)$ is $\be\eta$-saturated.
  Let $M \rhd_{\be\eta}^* N$ and $N \in {\mathbb I}_{\be\eta}(a)$.
  \begin{itemize}
  \item If $N \in {\cal O}^{\oslash}$ then $N \in {\cal M}^\oslash$ and
    $\exists L$ and $x \in {\cal V}_1$ such that $x^L \in
    \fv{N}$. By theorem~\ref{deg=}.\ref{deg=three}, $\fv{N} \subseteq
    \fv{M}$ and $\deg(M) = \deg(N)$, hence, $M \in {\cal
      O}^{\oslash}$
  \item If $N \in\{M \in {\cal M}^\oslash$ / $M : \<{\mathbb
      H}^\oslash \v^* a\> \}$ then $N \rhd_{\be\eta}^* N'$ and
    $\exists \G \subset {\mathbb H}^\oslash$, such that $N' : \< \G \v
    a\>$. Hence $M \rhd_{\be\eta}^* N'$ and since by
    theorem~\ref{deg=}.\ref{deg=three}, $\deg(M) = \deg(N')$, $M
    \in\{M \in {\cal M}^\oslash$ / $M : \<{\mathbb H}^\oslash \v^* a\> \}$.
  \end{itemize}
  Case $r = \beta$.
  It is easy to see that $\forall x \in {\cal V}_1, {\cal
    N}_x^\oslash \subseteq {\cal
    O}^\oslash \subseteq {\mathbb I}_{\be}(a)$.
  Now we show that ${\mathbb I}_{\be}(a)$ is $\be$-saturated.
  Let $M \rhd_{\be}^* N$ and $N \in {\mathbb I}_{\be}(a)$.
  \begin{itemize}
  \item If $N \in {\cal O}^{\oslash}$ then $N \in {\cal M}^\oslash$ and
    $\exists L$ and $x \in {\cal V}_1$ such that $x^L \in
    \fv{N}$. By theorem~\ref{deg=}.\ref{deg=one}, $\fv{N} \subseteq
    \fv{M}$ and $\deg(M) = \deg(N)$, hence, $M \in {\cal O}^{\oslash}$
  \item If $N \in\{M \in {\cal M}^\oslash$ / $M : \<{\mathbb
      H}^\oslash \v a\> \}$ then $\exists \G \subset {\mathbb
      H}^\oslash$, such that $N : \< \G \v a\>$. By
    theorem~\ref{finalbetaexp}, $M : \<\G {\u^{M}} \v a\>$. Since by
    theorem~\ref{deg=}.\ref{deg=one}, $\fv{N} \subseteq \fv{M}$, let
    $\fv{N} = \{x_1^{L_1}, \dots, x_n^{L_n}\}$ and $\fv{M} = \fv{N} \cup
    \{x_{n+1}^{L_{n+1}}, \dots, x_{n+m}^{L_{n+m}}\}$. So $\G {\u^{M}}
    = \G, (x_{n+1}^{L_{n+1}} : \o^{L_{n+1}}, \dots, x_{n+m}^{L_{n+m}}
    : \o^{L_{n+m}})$. $\forall n+1 \leq i \leq n+m$, let $U_i$ such
    that $x^{L_i}_i \in {\mathbb V}_{U_i}$. Then $\G, (x_{n+1}^{L_{n+1}} :
    U_{n+1}, \dots, x_{n+m}^{L_{n+m}} : U_{n+m}) \subset {\mathbb
      H}^\oslash$ and by $\sqsubseteq$, $M : \<\G, (x_{n+1}^{L_{n+1}} :
    U_{n+1}, \dots, x_{n+m}^{L_{n+m}} : U_{n+m}) \v a\>$. Thus $M :
    \<{\mathbb H}^\oslash \v a\>$ and since by
    theorem~\ref{deg=}.\ref{deg=three}, $\deg(M) = \deg(N)$, $M \in
    \{M \in {\cal M}^\oslash$ / $M : \<{\mathbb H}^\oslash \v a\>
    \}$.
  \end{itemize}
  2.\  By induction on $U$.
  \begin{itemize}
  \item
    $U = a$: By definition of ${\mathbb I}_{\be\eta}$.

  \item $U = \o^L$: By definition, ${\mathbb I}_{\be\eta}(\o^L) =
    {\cal M}^L$.  Hence, ${\cal O}^L \cup \{M \in {\cal M}^L$ / $M :
    \<{\mathbb
      H}^L \v^* \o^L\> \} \subseteq {\mathbb I}_{\be\eta}(\o^L)$.\\
    Let $M \in {\mathbb I}_{\be\eta}(\o^L)$ where $\fv{M} =
    \{x^{L_1}_1,...,x^{L_n}_n\}$ then $M \in {\cal M}^L$.  $\forall~1
    \leq i \leq n$, let $U_i$ the type such that $x^{L_i}_i \in
    {\mathbb V}_{U_i}$.  Then $\G = (x^{L_i}_i:U_i)_n \subset {\mathbb
      H}^L$. By lemma~\ref{structyping}.\ref{structypingtwo} and
    lemma~\ref{GnHn}, $M : \<\G \v \o^L\>$.  Hence $M : \<{\mathbb
      H}^L \v \o^L\>$. Therefore, ${\mathbb I}(\o^L) \subseteq \{M \in
    {\cal M}^L$ / $M : \< {\mathbb H}^L \v^* \o^L\> \}$.

    We deduce ${\mathbb I}_{\be\eta}(\o^L) = {\cal O}^L \cup \{M \in {\cal
      M}^L$ / $M : \< {\mathbb H}^L \v^* \o^L\> \}$.

  \item
    $U = \overline{e}_iV$: $L = i::K$ and $\deg(V) = K$. By
    IH and lemma~\ref{Vx}, ${\mathbb I}_{\be\eta}(\overline{e}_iV) = ({\mathbb
      I}_{\be\eta}(V))^{+i}
    = $ $({\cal O}^{K} \cup \{M \in {\cal M}^{K}$ / $M : \<
    {\mathbb H}^{K} \v^* V\> \})^{+i} =$ \\${\cal O}^{L} \cup (\{M
    \in {\cal M}^{K}$ / $M : \< {\mathbb H}^{K} \v^* V\> \})^{+i}$.
    \begin{itemize}
    \item If $M \in {\cal M}^{K}$ and $M : \<{\mathbb H}^{K} \v^*
      V\>$, then $M \rhd^*_{\beta\eta} N$ and $N : \<\G \v V\>$ where
      $\G \subset {\mathbb H}^{K}$. By $e$, lemmas~\ref{deg+-f+}
      and~\ref{GnHn}, $N^{+i} : \<\overline{e}_i\G \v \overline{e}_iV\>$, $M^{+i}
      \rhd^*_{\beta\eta} N^{+i}$ and $\overline{e}_i\G \subset {\mathbb
        H}^L$. Thus $M^{+i} \in {\cal M}^L$ and $M^{+i} : \<{\mathbb
        H}^L \v^* U\>$.
    \item If $M \in {\cal M}^L$ and $M : \< {\mathbb H}^L \v^* U\>$,
      then $M \rhd^*_{\beta\eta} N$ and $N : \<\G \v U\>$ where $\G
      \subset {\mathbb H}^L$. By lemmas~\ref{deg+-f+},~\ref{typing'},
      and~\ref{GnHn}, $M^{-i} \rhd^*_{\beta\eta} N^{-i}$, $N^{-i} :
      \<\G^{-i} \v V\>$ and $\G^{-i} \subset {\mathbb H}^{K}$. Thus by
      lemma~\ref{deg+-f+}, $M = (M^{-i})^{+i}$ and $M^{-i} \in \{M \in
      {\cal M}^K$ / $M : \< {\mathbb H}^K \v^* V\> \}$.
    \end{itemize}
    Hence $(\{M \in {\cal M}^K$ / $M : \< {\mathbb H}^K \v^* V\>
    \})^{+i} =  \{M \in {\cal M}^L$ / $M : \< {\mathbb H}^L \v^* U\>
    \}$ and ${\mathbb I}_{\be\eta}(U) = {\cal O}^L \cup \{M \in {\cal M}^L$ /
    $M : \< {\mathbb H}^L \v^* U\> \}$.

  \item
    $U = U_1 \sqcap U_2$:
    By IH, ${\mathbb I}_{\be\eta}(U_1 \sqcap U_2) = {\mathbb I}_{\be\eta}(U_1) \cap
    {\mathbb I}_{\be\eta}(U_2) = $ $({\cal O}^L \cup \{M \in {\cal M}^L$ / $M
    : \< {\mathbb H}^L \v^* U_1\> \}) \cap ({\cal O}^L \cup \{M \in
    {\cal M}^L$ / $M : \< {\mathbb H}^L \v^* U_2\> \}) = $ ${\cal
      O}^L \cup (\{M \in {\cal M}^L$ / $M : \< {\mathbb H}^L \v^*
    U_1\> \} \cap \{M \in {\cal M}^L$ / $M : \< {\mathbb H}^L \v^*
    U_2\> \})$.
    \begin{itemize}
    \item If $M \in {\cal M}^L$, $M : \< {\mathbb H}^L \v^* U_1\>$ and
      $M : \< {\mathbb H}^L \v^* U_2\>$, then $M \rhd^*_{\beta\eta}
      N_1$, $M \rhd^*_{\beta\eta} N_2$, $N_1 : \< \G_1 \v U_1\>$ and
      $N_2 : \< \G_2 \v U_2\>$ where $\G_1, \G_2 \subset {\mathbb
        H}^L$. By confluence theorem~\ref{confluenceofbetaeta} and
      subject reduction theorem~\ref{subred}, $\exists M'$ such that
      $M \rhd^*_{\beta\eta} M'$, $M' : \< \G_1\r_{M'} \v U_1\>$ and
      $M' : \< \G_2\r_{M'} \v U_2\>$. Hence by Remark~\ref{remderiv}
      and lemma~\ref{deg=} and
      lemma~\ref{structyping}.\ref{structypingone} and
      lemma~\ref{restlem}.\ref{restlemtwo}, $M' : \< (\G_1 \sqcap
      \G_2)\r_{M'} \v U_1 \sqcap U_2 \>$ and, by lemma~\ref{GnHn},
      $(\G_1 \sqcap \G_2)\r_{M'} \subseteq \G_1 \sqcap \G_2 \subset
      {\mathbb H}^L$. Thus $M : \< {\mathbb H}^L \v^* U_1 \sqcap
      U_2\>$.
    \item
      If $M \in {\cal M}^L$  and $M : \< {\mathbb H}^L \v^* U_1 \sqcap
      U_2\>$, then $M \rhd^*_{\beta\eta} N$, $N : \< \G \v U_1 \sqcap
      U_2\>$ and $\G \subset {\mathbb H}^L$. By $\sqsubseteq$, $N : \<
      \G \v U_1\>$ and $N : \< \G \v U_2\>$.\\
      Hence, $M : \< {\mathbb H}^L \v^* U_1\>$ and $M : \< {\mathbb H}^L
      \v^* U_2\>$.
    \end{itemize}
    We deduce that ${\mathbb I}_{\be\eta}(U_1 \sqcap T_2) = {\cal O}^L \cup
    \{M \in {\cal M}^L$ / $M : \< {\mathbb H}^L \v^* U_1 \sqcap U_2\>
    \}$.

  \item
    $U = V \f T$: Let $\deg(T) = \oslash \preceq K = \deg(V)$.
    By IH, ${\mathbb I}_{\be\eta}(V) = {\cal O}^K \cup \{M \in {\cal M}^K$ /
    $M : \< {\mathbb H}^K \v^* V\> \}$ and ${\mathbb I}_{\be\eta}(T) = {\cal
      O}^{\oslash} \cup \{M \in {\cal M}^{\oslash}$ / $M : \<
    {\mathbb H}^{\oslash} \v^* T\> \}$. Note that ${\mathbb I}_{\be\eta}(V \f T)
    = {\mathbb I}_{\be\eta}(V) \fx {\mathbb I}_{\be\eta}(T)$.
    \begin{itemize}
    \item
      Let $M \in {\mathbb I}_{\be\eta}(V) \fx {\mathbb I}_{\be\eta}(T)$
      and, by lemma~\ref{codeV},
      let $y^K \in {\mathbb V}_V$ such that $\forall K, y^K \not
      \in \fv{M}$. Then $M \diamond y^K$. By remark~\ref{remderiv}, $y^K
      : \<(y^K : V) \v^*
      V\>$. Hence $y^K : \<{\mathbb H}^K \v^* V\>$. Thus, $y^K \in
      {\mathbb I}_{\be\eta}(V)$ and  $M y^K \in {\mathbb I}_{\be\eta}(T)$.
      \begin{itemize}
      \item
        If $M  y^K \in {\cal O}^{\oslash}$, then since $y \in {\cal V}_2$,
        by lemma~\ref{Vx}, $M \in {\cal O}^{\oslash}$.
      \item If $M y^K \in \{M \in {\cal M}^{\oslash}$ / $M : \<
        {\mathbb H}^{\oslash} \v^* T\> \}$ then $M y^K
        \rhd^*_{\beta\eta} N$ and $N : \<\G \v T\>$ such that $\G
        \subset {\mathbb H}^{\oslash}$, hence, $\l y^K.M y^K
        \rhd_{\beta\eta}^* \l y^K.N$. We have two cases:
        \begin{itemize}
        \item
          If $y^K \in \dom{\G}$,
          then $\G = \Delta , (y^K:V)$ and by $\f_I$, $\l y^K. N : \< \Delta
          \v V \f T\>$.
        \item
          If $y^K \not \in \dom{\G}$, let $\Delta = \G$.
          By $\f'_I$, $\l y^K. N : \< \Delta \v \o^K \f T\>$. By
          $\sqsubseteq$, since $\< \Delta \v \o^K \f T\> \sqsubseteq \<
          \Delta \v V \f T\>$, we have $\l y^K. N : \< \Delta \v  V \f
          T\>$.
        \end{itemize}
        Note that $\Delta \subset {\mathbb H}^{\oslash}$. Since $\l
        y^K.M y^K \rhd_{\beta\eta}^* M$ and $\l y^K.M y^K
        \rhd_{\beta\eta}^* \l y^K.N$, by
        theorem~\ref{confluenceofbetaeta} and theorem~\ref{subred},
        there is $M'$ such that $M \rhd_{\beta\eta}^* M'$, $ \l y^K.N
        \rhd_{\beta\eta}^* M'$, $M' : \< \Delta\r_{M'} \v V \f
        T\>$. Since $\Delta\r_{M'} \subseteq \Delta\subset {\mathbb
          H}^{\oslash}$, $M : \< {\mathbb H}^{\oslash} \v^* V \f T\>$.
      \end{itemize}
    \item Let $M \in {\cal O}^\oslash \cup \{M \in {\cal M}^\oslash$ /
      $M : \< {\mathbb H}^\oslash \v^* V \f T\> \}$ and $N \in
      {\mathbb I}_{\be\eta}(V) = {\cal O}^K \cup \{M \in {\cal M}^K$ /
      $M : \< {\mathbb H}^K \v^* V\> \}$ such that $M \diamond
      N$. Then, $\deg(N) = K \succeq \oslash = \deg(M)$.
      \begin{itemize}
      \item
        If $M \in {\cal O}^{\oslash}$, then, by
        lemma \ref{Vx}, $M N \in {\cal O}^{\oslash}$.
      \item
        If $M \in \{M \in {\cal M}^{\oslash}$ /
        $M : \< {\mathbb H}^{\oslash} \v^* V \f T\> \}$, then
        \begin{itemize}
        \item
          If $N \in {\cal O}^K$, then, by lemma~\ref{Vx}, $M  N \in
          {\cal O}^{\oslash}$.
        \item If $N \in \{M \in {\cal M}^K$ / $M : \< {\mathbb H}^K
          \v^* V\> \}$ then $M \rhd^*_{\beta\eta} M_1$, $N
          \rhd^*_{\beta \eta} N_1$, $M_1 : \< \G_1 \v V \f T\>$ and
          $N_1 : \< \G_2 \v V\>$ where $\G_1 \subset {\mathbb
            H}^{\oslash}$ and $\G_2 \subset {\mathbb H}^K$. By
          lemma~\ref{deg+-f+} and theorem~\ref{deg=}, $M N
          \rhd^*_{\beta \eta} M_1 N_1$ and, by $\f_E$ and
          lemma~\ref{structyping}.\ref{structypingC}, $M_1 N_1 : \<
          \G_1 \sqcap \G_2 \v T\>$. By lemma~\ref{GnHn}, $\G_1 \sqcap
          \G_2 \subset {\mathbb H}^{\oslash}$. Therefore $M N : \<
          {\mathbb H}^{\oslash} \v^* T\>$.
        \end{itemize}
      \end{itemize}
    \end{itemize}
    We deduce  that ${\mathbb I}_{\be\eta}(V \f T) = {\cal O}^{\oslash} \cup
    \{M \in {\cal M}^{\oslash}$ / $M : \< {\mathbb H}^{\oslash} \v^* V
    \f T\> \}$.
  \end{itemize}
  3.\ We only do the
  case $r = \beta$.
  By induction on $U$.

  \begin{itemize}
  \item
    $U = a$: By definition of ${\mathbb I}_{\be}$.

  \item $U = \o^L$: By definition, ${\mathbb I}_{\be}(\o^L) = {\cal
      M}^L$.  Hence, ${\cal O}^L \cup \{M \in {\cal M}^L$ / $M :
    \<{\mathbb
      H}^L \v \o^L\> \} \subseteq {\mathbb I}_{\be}(\o^L)$.\\
    Let $M \in {\mathbb I}_{\be}(\o^L)$ where $\fv{M} =
    \{x^{L_1}_1,...,x^{L_n}_n\}$ then $M \in {\cal M}^L$.  $\forall~1
    \leq i \leq n$, let $U_i$ the type such that $x^{L_i}_i \in
    {\mathbb V}_{U_i}$.  Then $\G = (x^{L_i}_i:U_i)_n \subset {\mathbb
      H}^L$. By lemma~\ref{structyping}.\ref{structypingtwo} and
    lemma~\ref{GnHn}, $M : \<\G \v \o^L\>$.  Hence $M : \<{\mathbb
      H}^L \v \o^L\>$. Therefore, ${\mathbb I}(\o^L) \subseteq \{M \in
    {\cal M}^L$ / $M : \< {\mathbb H}^L \v \o^L\> \}$.

    We deduce ${\mathbb I}_{\be}(\o^L) = {\cal O}^L \cup \{M \in {\cal
      M}^L$ / $M : \< {\mathbb H}^L \v \o^L\> \}$.

  \item
    $U = \overline{e}_iV$: $L = i::K$ and $\deg(V) = K$. By
    IH and lemma~\ref{Vx}, ${\mathbb I}_{\be}(\overline{e}_iV) = ({\mathbb
      I}_{\be}(V))^{+i}
    = $ $({\cal O}^{K} \cup \{M \in {\cal M}^{K}$ / $M : \<
    {\mathbb H}^{K} \v V\> \})^{+i} =$ \\${\cal O}^{L} \cup (\{M
    \in {\cal M}^{K}$ / $M : \< {\mathbb H}^{K} \v V\> \})^{+i}$.
    \begin{itemize}
    \item If $M \in {\cal M}^{K}$ and $M : \<{\mathbb H}^{K} \v V\>$,
      then $M : \<\G \v V\>$ where $\G \subset {\mathbb H}^{K}$. By
      $e$ and~\ref{GnHn}, $M^{+i} : \<\overline{e}_i\G \v \overline{e}_iV\>$ and $\overline{e}_i\G
      \subset {\mathbb H}^L$. Thus $M^{+i} \in {\cal M}^L$ and $M^{+i}
      : \<{\mathbb H}^L \v U\>$.
    \item If $M \in {\cal M}^L$ and $M : \< {\mathbb H}^L \v U\>$,
      then $M : \<\G \v U\>$ where $\G \subset {\mathbb H}^L$. By
      lemmas~\ref{typing'}, and~\ref{GnHn}, $M^{-i} : \<\G^{-i} \v
      V\>$ and $\G^{-i} \subset {\mathbb H}^{K}$. Thus by
      lemma~\ref{deg+-f+}, $M = (M^{-i})^{+i}$ and $M^{-i} \in \{M \in
      {\cal M}^K$ / $M : \< {\mathbb H}^K \v V\> \}$.
    \end{itemize}
    Hence $(\{M \in {\cal M}^K$ / $M : \< {\mathbb H}^K \v V\>
    \})^{+i} =  \{M \in {\cal M}^L$ / $M : \< {\mathbb H}^L \v U\>
    \}$ and ${\mathbb I}_{\be}(U) = {\cal O}^L \cup \{M \in {\cal M}^L$ /
    $M : \< {\mathbb H}^L \v U\> \}$.

  \item
    $U = U_1 \sqcap U_2$:
    By IH, ${\mathbb I}_{\be}(U_1 \sqcap U_2) = {\mathbb I}_{\be}(U_1) \cap
    {\mathbb I}_{\be}(U_2) = $ $({\cal O}^L \cup \{M \in {\cal M}^L$ / $M
    : \< {\mathbb H}^L \v U_1\> \}) \cap ({\cal O}^L \cup \{M \in
    {\cal M}^L$ / $M : \< {\mathbb H}^L \v U_2\> \}) = $ ${\cal
      O}^L \cup (\{M \in {\cal M}^L$ / $M : \< {\mathbb H}^L \v
    U_1\> \} \cap \{M \in {\cal M}^L$ / $M : \< {\mathbb H}^L \v
    U_2\> \})$.
    \begin{itemize}
    \item
      If $M \in {\cal M}^L$, $M : \< {\mathbb H}^L
      \v U_1\>$ and $M : \< {\mathbb H}^L \v U_2\>$, then $M : \<
      \G_1 \v U_1\>$ and $M : \< \G_2 \v U_2\>$ where
      $\G_1, \G_2 \subset {\mathbb H}^L$. Hence by
      Remark~\ref{remderiv}, $M : \<\G_1 \sqcap
      \G_2 \v U_1 \sqcap U_2\>$ and, by lemma~\ref{GnHn}, $\G_1 \sqcap
      \G_2 \subset {\mathbb H}^L$. Thus $M : \< {\mathbb H}^L \v U_1
      \sqcap U_2\>$.
    \item
      If $M \in {\cal M}^L$  and $M : \< {\mathbb H}^L \v U_1 \sqcap
      U_2\>$, then $M : \< \G \v U_1 \sqcap
      U_2\>$ and $\G \subset {\mathbb H}^L$. By $\sqsubseteq$, $M : \<
      \G \v U_1\>$ and $M : \< \G \v U_2\>$.
      Hence, $M : \< {\mathbb H}^L \v U_1\>$ and $M : \< {\mathbb H}^L
      \v U_2\>$.
    \end{itemize}
    We deduce that ${\mathbb I}_{\be}(U_1 \sqcap T_2) = {\cal O}^L \cup
    \{M \in {\cal M}^L$ / $M : \< {\mathbb H}^L \v U_1 \sqcap U_2\>
    \}$.

  \item
    $U = V \f T$: Let $\deg(T) = \oslash \preceq K = \deg(V)$.
    By IH, ${\mathbb I}_{\be}(V) = {\cal O}^K \cup \{M \in {\cal M}^K$ /
    $M : \< {\mathbb H}^K \v V\> \}$ and ${\mathbb I}_{\be}(T) = {\cal
      O}^{\oslash} \cup \{M \in {\cal M}^{\oslash}$ / $M : \<
    {\mathbb H}^{\oslash} \v T\> \}$. Note that ${\mathbb I}_{\be}(V \f T)
    = {\mathbb I}_{\be}(V) \fx {\mathbb I}_{\be}(T)$.
    \begin{itemize}
    \item
      Let $M \in {\mathbb I}_{\be}(V) \fx {\mathbb I}_{\be}(T)$
      and, by lemma~\ref{codeV},
      let $y^K \in {\mathbb V}_V$ such that $\forall K, y^K \not
      \in \fv{M}$. Then $M \diamond y^K$. By remark~\ref{remderiv}, $y^K
      : \<(y^K : V) \v^*
      V\>$. Hence $y^K : \<{\mathbb H}^K \v V\>$. Thus, $y^K \in
      {\mathbb I}_{\be}(V)$ and  $M y^K \in {\mathbb I}_{\be}(T)$.
      \begin{itemize}
      \item
        If $M  y^K \in {\cal O}^{\oslash}$, then since $y \in {\cal V}_2$,
        by lemma~\ref{Vx}, $M \in {\cal O}^{\oslash}$.
      \item If $M y^K \in \{M \in {\cal M}^{\oslash}$ / $M : \<
        {\mathbb H}^{\oslash} \v T\> \}$ then $My^K : \<\G \v T\>$
        such that $\G \subset {\mathbb H}^{\oslash}$. Since by
        lemma~\ref{structyping}.\ref{structypingone}, $\dom{\G} =
        \fv{My^K}$ and $y^K \in \fv{My^K}$, $\G = \Delta,
        (y^K:V')$. Since $(y^K:V') \in {\mathbb H}^{\oslash}$, by
        lemma~\ref{codeV}, $V = V'$. So $My^K : \<\D, (y^K:V) \v T\>$
        and by lemma~\ref{newgen} $M : \< \D \v V \f T\>$.  Note that
        $\D \subset {\mathbb H}^\oslash$, hence $M : \< {\mathbb
          H}^{\oslash} \v V \f T\>$.
      \end{itemize}
    \item Let $M \in {\cal O}^\oslash \cup \{M \in {\cal M}^\oslash$ /
      $M : \< {\mathbb H}^\oslash \v V \f T\> \}$ and $N \in {\mathbb
        I}_{\be\eta}(V) = {\cal O}^K \cup \{M \in {\cal M}^K$ / $M :
      \< {\mathbb H}^K \v V\> \}$ such that $M \diamond N$. Then,
      $\deg(N) = K \succeq \oslash = \deg(M)$.
      \begin{itemize}
      \item
        If $M \in {\cal O}^{\oslash}$, then, by
        lemma~\ref{Vx}, $M N \in {\cal O}^{\oslash}$.
      \item
        If $M \in \{M \in {\cal M}^{\oslash}$ /
        $M : \< {\mathbb H}^{\oslash} \v V \f T\> \}$, then
        \begin{itemize}
        \item
          If $N \in {\cal O}^K$, then, by lemma~\ref{Vx}, $M  N \in
          {\cal O}^{\oslash}$.
        \item If $N \in \{M \in {\cal M}^K$ / $M : \< {\mathbb H}^K \v
          V\> \}$ then $M : \< \G_1 \v V \f T\>$ and $N : \< \G_2 \v
          V\>$ where $\G_1 \subset {\mathbb H}^{\oslash}$ and $\G_2
          \subset {\mathbb H}^K$. By $\f_E$ and
          lemma~\ref{structyping}.\ref{structypingC}, $MN : \< \G_1
          \sqcap \G_2 \v T\>$. By lemma~\ref{GnHn}, $\G_1 \sqcap \G_2
          \subset {\mathbb H}^{\oslash}$. Therefore $MN : \< {\mathbb
            H}^{\oslash} \v T\>$.
        \end{itemize}
      \end{itemize}
    \end{itemize}
    We deduce  that ${\mathbb I}_{\be}(V \f T) = {\cal O}^{\oslash} \cup
    \{M \in {\cal M}^{\oslash}$ / $M : \< {\mathbb H}^{\oslash} \v V
    \f T\> \}$. \hfill $\Box$
  \end{itemize}
\end{proof}

\fi

Now, we use this crucial  ${\mathbb I}$ to establish completeness of our semantics.
\begin{theorem}[Completeness of $\v$]\label{comple}
  Let $U\in {\mathbb U}$ such that $\deg(U) = L$.
  \begin{enumerate}
  \item\label{compleone}
    $[U]_{\be\eta} = \{M \in {\cal M}^L$ / $M$ closed, $M
    \rhd^*_{\beta \eta} N$ and $N : \< () \v U\> \}$.
  \item\label{completwo}
    $[U]_{\be} = [U]_{\wh} = \{M \in {\cal M}^L$ / $M : \< () \v U\>
    \}$.
  \item\label{complethree}
    $[U]_{\be\eta}$ is stable by reduction. I.e., If $M \in
    [U]_{\be\eta}$ and $M \rhd_{\beta\eta}^* N$ then $N \in
    [U]_{\be\eta}$.
  \end{enumerate}
\end{theorem}

\begin{proof}
  Let $r \in \{\be, \wh, \be\eta\}$.
  \begin{itemize}
  \item[\ref{compleone}.]  Let $M \in [U]_{\be\eta}$. Then $M$ is a
    closed term and $M \in {\mathbb I}_{\be\eta}(U)$. Hence, by Lemma
    \ref{comp}, $M \in {\cal O}^L \cup \{M \in {\cal M}^L$ / $M : \<
    {\mathbb H}^L \v^* U\> \}$. Since $M$ is closed, $M \not \in {\cal
      O}^L$. Hence, $M \in \{M \in {\cal M}^L$ / $M : \< {\mathbb H}^L
    \v^* U\>\}$ and so, $M \rhd^*_{\beta\eta} N$ and $N :\<\G \v U\>$
    where $\G \subset {\mathbb H}^L$. By Theorem~\ref{deg=}, $N$ is
    closed and, by Lemma~\ref{structyping}.\ref{structypingone}, $N :
    \< () \v U\>$.

    Conversely, take $M$ closed such that $M \rhd^*_{\beta} N$ and $N
    : \< () \v U\>$. Let ${\cal I}\in \be\eta\mbox{-int}$.  By
    Lemma~\ref{adeq}, $N \in {\cal I}(U)$. By
    Lemma~\ref{interpret-intsub}.\ref{interpret}, ${\cal I}(U)$ is
    $\be\eta$-saturated.  Hence, $M \in {\cal I}(U)$. Thus $M \in
    [U]$.

  \item[\ref{completwo}.]  Let $M \in [U]_{\be}$. Then $M$ is a closed
    term and $M \in {\mathbb I}_{\be}(U)$.  Hence, by Lemma
    \ref{comp}, $M \in {\cal O}^L \cup \{M \in {\cal M}^L$ / $M : \<
    {\mathbb H}^L \v U\>\}$. Since $M$ is closed, $M \not \in {\cal
      O}^L$. Hence, $M \in \{M \in {\cal M}^L$ / $M : \< {\mathbb H}^L
    \v U\>\}$ and so, $M :\< \G \v U\>$ where $\G \subset {\mathbb
      H}^L$. By Lemma~\ref{structyping}.\ref{structypingone}, $M : \<
    () \v U\>$.

    Conversely, take $M$ such that $M : \< () \v U\>$. By
    Lemma~\ref{structyping}.\ref{structypingone}, $M$ is closed. Let
    ${\cal I}\in\be\mbox{-int}$.  By Lemma~\ref{adeq}, $M \in {\cal
      I}(U)$. Thus $M \in [U]_{\be}$.

    It is easy to see that $[U]_{\be} = [U]_{\wh}$.

  \item[\ref{complethree}.]  Let $M \in [U]_{\be\eta}$ and $M
    \rhd_{\beta\eta}^* N$.  By~\ref{compleone}, $M$ is closed, $M
    \rhd^*_{\beta\eta} P$ and $P : \< () \v U\>$. By confluence
    Theorem~\ref{confluenceofbetaeta}, there is $Q$ such that $P
    \rhd^*_{\beta\eta} Q$ and $N \rhd^*_{\beta\eta} Q$. By subject
    reduction Theorem~\ref{subred}, $Q : \< () \v U\>$. By
    Theorem~\ref{deg=}, $N$ is closed and, by~\ref{compleone}, $N \in
    [U]_{\be\eta}$. \hfill $\Box$
  \end{itemize}
\end{proof}

\section{Conclusion}
\label{concsec}

Expansion may be viewed to work like a multi-layered simultaneous
substitution.  Moreover, expansion is a crucial part of a procedure
for calculating principal typings and helps support compositional type
inference.  Because the early definitions of expansion were
complicated, expansion variables (E-variables) were introduced to
simplify and mechanise expansion. The aim of this paper is to give a
complete semantics for intersection type systems with expansion
variables.

The only earlier attempt (see Kamareddine, Nour, Rahli and
Wells~\cite{report}) at giving a semantics for expansion variables
could only handle the $\lambda I$-calculus, did not allow a universal
type, and was incomplete in the presence of more than one expansion
variable. This paper overcomes these difficulties and gives a complete
semantics for an intersection type system with an arbitrary (possibly
infinite) number of expansion variables using a calculus indexed with
finite sequences of natural numbers.

\bibliographystyle{jbwc}
{\bibliography{literature,macros,bibliography,conferences}}

\ifictac
\else
\newpage
\appendix

\section{Proofs of Section~\ref{secpure}}

The next lemma is needed in the proofs.
\begin{lemma}
  \label{degsub}
  Let $M, M', N, N_1, \dots, N_n \in {\cal M}$.
  \begin{enumerate}
  \item
    \label{lem:refl+symm}
    $M \diamond M$ and if $M \diamond N$ then $N \diamond M$.
  \item
    \label{lem:degsubfvsub}
    If $\fv{M} \subseteq \fv{M'}$ and $M' \diamond N$ then $M \diamond
    N$.
  \item
    \label{degsubone}
    If $M \diamond N$ and $M'$ is a subterm of $M$
    then $M' \diamond N$.
  \item
    \label{degsubtwo}
    If $\deg(M) = L$ and $x^K$ occurs in $M$, then
    $K \succeq L$.
  \item
    \label{degsubthree}
    If ${\cal X} = \{M\}\cup \{N_i / 1 \leq i \leq n\}$, for all $i
    \in \{1, \dots, n\}$, $\deg(N_i) = L_i$ and $\diamond {\cal X}$
    then $M[(x^{L_i}_i := N_i)_n] \in {\cal M}$ and $\deg(M[(x^{L_i}_i
    := N_i)_n]) = \deg(M)$.
  \item
    \label{degsubfour}
    If ${\cal X} = \{M,N\}\cup \{N_i / 1 \leq i \leq n\}$, for all $i
    \in \{1, \dots, n\}$, $\deg(N_i) = L_i$ and $\diamond {\cal X}$
    then $M[(x^{L_i}_i := N_i)_n] \diamond N[(x^{L_i}_i := N_i)_n]$
  \end{enumerate}
\end{lemma}

\begin{proof}
  \begin{enumerate}
  \item First, we prove $M \diamond M$ by induction on $M$.
    \begin{itemize}
    \item Let $M = x^{L}$ then it is trivial.
    \item Let $M = \l x^{L}. N$ such that $N \in {\cal M}$ and $L
      \succeq \deg(N)$. Let $y^{K}, y^{K'} \in \fv{M}$ then $y^{K},
      y^{K'} \in \fv{N}$ and we conclude using IH on $N$.
    \item Let $M = M_1M_2$ such that $M_1, M_2 \in {\cal M}$,
      $\deg(M_1) \preceq \deg(M_2)$ and $M_1 \diamond M_2$. Let
      $x^{L}, x^{K} \in \fv{M}$ then either $x^{L}, x^{K} \in
      \fv{M_1}$ and we conclude using IH on $M_1$. Or $x^{L}, x^{K}
      \in \fv{M_2}$ and we conclude using IH on $M_2$. Or $x^{L} \in
      \fv{M_1}$ and $x^{K} \in \fv{M_2}$ and we conclude using $M_1
      \diamond M_2$.
    \end{itemize}

    Let $M \diamond N$, we prove $N \diamond M$. It is trivial by
    definition.

  \item Let $x^{L} \in \fv{M} \subseteq \fv{M'}$ and $x^{K} \in
    \fv{N}$ then by hypothesis $K = L$.

  \item
    By induction on $M$.
    \begin{itemize}
    \item Case $M = x^L$ is trivial.
    \item Case $M = \l x^L. P$ where $\forall K \in {\cal L}_{\mathbb
        N}, x^K \not \in \fv{N}$. If $M' = M$ then nothing to
      prove. Else $M'$ is a subterm of $P$. If we prove that
      $P\diamond N$ then we can use IH to get $M' \diamond N$.  Hence,
      now we prove $P\diamond N$.  Let $y \in {\cal V}$ such that $y^K
      \in \fv{P}$ and $y^{K'} \in \fv{N}$. Since $x^{K'} \not \in
      \fv{N}$, then $x \not = y$ and $y^K \not = x^L$.  Hence $y^K \in
      \fv{M}$ and since $M \diamond N$ then $K = K'$. Hence, $P
      \diamond N$.
    \item Case $M = M_1M_2$. Let $i \in \{1, 2\}$. First we prove that
      $M_i \diamond N$: let $ x \in {\cal V}$, such that $x^L \in
      \fv{M_i}$ and $x^K \in \fv{N}$, then $x^L \in \fv{M}$ and so $L =
      K$.  Now, if $M' = M$ then nothing to prove.  Else
      \begin{itemize}
      \item Either $M'$ is a subterm of $M_1$ and so by IH, since $M_1
        \diamond N$, $M' \diamond N$.
      \item Or $M'$ is a subterm of $M_2$ and so by IH, since $M_2
        \diamond N$, $M' \diamond N$.
      \end{itemize}
    \end{itemize}

  \item By induction on $M$.
    \begin{itemize}
    \item If $M = x^{K}$ then $\deg(M) = K$ and since $\succeq$ is an
      order relation, $K \succeq K$.
    \item If $M = M_1M_2$ then $\deg(M) = \deg(M_1)$. Let $L' =
      \deg(M_2)$ so $L' \succeq L$. By IH, if $x^{K}$ occurs in $M_1$
      then $K \succeq L$ and if $x^{K}$ occurs in $M_2$ then $K
      \succeq L'$. Since $x^{K}$ occurs in $M$, $K \succeq L$.
    \item If $M = \l x^{L_1}. M_1$ then $L_1 \succeq \deg(M_1) =
      \deg(\l x^{L_1}. M_1) = L$. If $x^{K}$ occurs in $M$, then
      $x^{K} = x^{L_1}$ or $x^{K}$ occurs in $M_1$. By IH, if $x^{K}$
      occurs in $M_1$ then $K \succeq L$.
    \end{itemize}

  \item By induction on $M$.
    \begin{itemize}
    \item If $M = y^K$ then if $y^K = x^{L_i}_i$, for $1 \leq i \leq
      n$, then $M[(x^{L_i}_i := N_i)_n] = N_i \in {\cal M}$ and
      $\deg(M[(x^{L_i}_i := N_i)_n]) = \deg(N_i) = L_i = K$. Else,
      $M[(x^{L_i}_i := N_i)_n] = y^K \in {\cal M}$ and
      $\deg(M[(x^{L_i}_i := N_i)_n]) = \deg(y^{K})$.
    \item If $M = M_1M_2$ then $\deg(M) = \deg(M_1)$ and $M[(x^{L_i}_i
      := N_i)_n] = M_1[(x^{L_i}_i := N_i)_n]M_2[(x^{L_i}_i :=
      N_i)_n]$. Since $\forall N \in {\cal X}, M \diamond N$, by
      \ref{degsubone}., $\forall N \in {\cal X}, M_1 \diamond N$ and
      $M_2 \diamond N$. Since $M_1, M_2 \in {\cal M}$, by IH,
      $M_1[(x^{L_i}_i := N_i)_n], M_2[(x^{L_i}_i := N_i)_n] \in {\cal
        M}$, $\deg(M_1[(x^{L_i}_i := N_i)_n]) = \deg(M_1)$ and
      $\deg(M_2[(x^{L_i}_i := N_i)_n]) = \deg(M_2)$. Let $x^K \in
      \fv{M_1[(x^{L_i}_i := N_i}_n])$ and $x^{K'} \in \fv{M_2[(x^{L_i}_i
      := N_i}_n])$. If $x^K \in \fv{M_1}$ then by \ref{degsubone}.,
      $\diamond (\{M_1,M_2\}\cup \{N_i / 1 \leq i \leq n\})$ hence $K
      = K'$. Let $1 \leq i \leq n$. If $x^K \in \fv{N_i}$ then by
      \ref{degsubone}., $\diamond (\{M_2\}\cup \{N_i / 1 \leq i \leq
      n\})$ hence $K = K'$. So $M_1[(x^{L_i}_i := N_i)_n] \diamond
      M_2[(x^{L_i}_i := N_i)_n]$. Furthermore, $\deg(M_2[(x^{L_i}_i :=
      N_i)_n]) = \deg(M_2) \succeq \deg(M_1) = \deg(M_1[(x^{L_i}_i :=
      N_i)_n])$ hence $M_1[(x^{L_i}_i := N_i)_n]M_2[(x^{L_i}_i :=
      N_i)_n] \in {\cal M}$ and $\deg(M_1[(x^{L_i}_i :=
      N_i)_n]M_2[(x^{L_i}_i := N_i)_n]) = \deg(M_1[(x^{L_i}_i :=
      N_i)_n]) = \deg(M_1) = \deg(M)$.
    \item If $M = \l y^K. M_1$ where $K \succeq \deg(M_1)$ and
      $\forall 1 \leq i \leq n$, $y \not = x_i$ and $\forall K' \in
      {\cal L}_{\mathbb N}$, $y^{K'} \not \in \fv{N_i} \cup
      \{x^{L_i}_i\}$ then $M[(x^{L_i}_i := N_i)_n] = \l
      y^K. M_1[(x^{L_i}_i := N_i)_n]$.  Since $M_1 \in {\cal M}$, then
      by \ref{degsubone}.\ and IH $M_1[(x^{L_i}_i := N_i)_n] \in {\cal
        M}$ and $\deg(M_1[(x^{L_i}_i := N_i)_n]) = \deg(M_1)$. So $\l
      y^K. M_1[(x^{L_i}_i := N_i)_n] \in {\cal M}$ and $\deg(\l
      y^K. M_1[(x^{L_i}_i := N_i)_n]) = \deg(M_1[(x^{L_i}_i :=
      N_i)_n]) = \deg(M_1) = \deg(M)$.
    \end{itemize}

  \item By \ref{degsubthree}., $M[(x^{L_i}_i := N_i)_n],N[(x^{L_i}_i
    := N_i)_n] \in {\cal M}$. Let $x^L \in \fv{M[(x^{L_i}_i :=
    N_i}_n])$ and $x^K \in \fv{N[(x^{L_i}_i := N_i}_n])$. So $x^L \in
    \fv{M}\cup \fv{N_1} \cup ... \cup \fv{N_n}$ and $x^K \in \fv{N}\cup
    \fv{N_1} \cup ... \cup \fv{N_n}$. Since $\diamond {\cal X}$, then $K
    = L$. Hence, $M[(x^{L_i}_i := N_i)_n] \diamond N[(x^{L_i}_i :=
    N_i)_n]$ \hfill $\Box$
  \end{enumerate}
\end{proof}

\begin{proof}[Of Theorem~\ref{deg=}]
  \begin{enumerate}
  \item By induction on $M \rhd_\eta^* N$, we only do the base step:
    \begin{itemize}
    \item $M = \l x^L. Nx^L \rhd_\eta N$ and $x^L \not \in \fv{N}$. By
      definition $\fv{M} = \fv{Nx^L}\setminus\{x^L\} = \fv{N}$ and
      $\deg(M) = \deg(Nx^L) = \deg(N)$.
    \item $M = \l x^L. M_1 \rhd_\eta \l x^L. N_1 = N$ and $M_1
      \rhd_\eta N_1$.  By IH, $\fv{N_1} = \fv{M_1}$ and $\deg(M_1) =
      \deg(N_1)$.  Hence, $\deg(M) = \deg(M_1) = \deg(N_1) = \deg(N)$
      and $\fv{N} = \fv{N_1} \setminus \{x^L\} = \fv{M_1} \setminus
      \{x^L\} = \fv{M}$.
    \item $M = M_1M_2 \rhd_\eta N_1M_2 = N$ such that $M_1 \rhd_\eta
      N_1$.  By IH, $\fv{N_1} = \fv{M_1}$ and $\deg(M_1) = \deg(N_1)$.
      By definition, $\fv{N} = \fv{N_1} \cup \fv{M_2} = \fv{M_1} \cup
      \fv{M_2} = \fv{M}$ and $\deg(M) = \deg(M_1) = \deg(N_1) =
      \deg(N)$.
    \item $M = M_1M_2 \rhd_\eta M_1N_2 = N$ such that $M_2 \rhd_\eta
      N_2$.  By IH, $\fv{N_2} = \fv{M_2}$ and $\deg(M_2) = \deg(N_2)$.
      By definition, $\fv{N} = \fv{M_1} \cup \fv{N_2} = \fv{M_1} \cup
      \fv{M_2} = \fv{M}$ and $\deg(M) = \deg(M_1) = \deg(N)$.
    \end{itemize}
  \item Case $r = \be$.  By induction on $M \rhd_\be^* N$, we only do
    the base step:
    \begin{itemize}
    \item $M = (\l x^L. M_1)M_2 \rhd_\be M_1[x^L:=M_2] = N$ such that
      $\deg(M_2) = L$.  If $x^L \in \fv{M_1}$ then $\fv{N} =
      (\fv{M_1}\setminus\{x^L\})\cup \fv{M_2} = \fv{M}$.  If $x^L \not
      \in \fv{M_1}$ then $\fv{N} = \fv{M_1} = \fv{M_1}\setminus
      \{x^L\}\subseteq \fv{M}$.  By definition, $\deg(M) =
      \deg(M_1)$. Because $N \in {\cal M}$ then $M_1 \diamond M_2$ and
      $\deg(M_2) = L$. So, by lemma~\ref{degsub}.\ref{degsubthree},
      $\deg(N) = \deg(M_1)$.
    \item $M = \l x^L. M_1 \rhd_\be \l x^L. N_1 = N$ such that $M_1
      \rhd_\be N_1$.  By IH, $\fv{N_1} \subseteq \fv{M_1}$ and
      $\deg(M_1) = \deg(N_1)$.  By definition $\deg(M) = \deg(M_1) =
      \deg(N_1) = \deg(N)$ and $\fv{N} = \fv{N_1}\setminus\{x^L\}
      \subseteq \fv{M_1}\setminus\{x^L\} = \fv{M}$.
    \item $M = M_1M_2 \rhd_\be N_1M_2 = N$ such that $M_1 \rhd_\be
      N_1$.  By IH, $\fv{N_1} \subseteq \fv{M_1}$ and $\deg(M_1) =
      \deg(N_1)$.  By definition, $\fv{N} = \fv{N_1} \cup \fv{M_2}
      \subseteq \fv{M_1} \cup \fv{M_2} = \fv{M}$ and $\deg(M) =
      \deg(M_1) = \deg(N_1) = \deg(N)$.
    \item $M = M_1M_2 \rhd_\be M_1N_2 = N$ such that $M_2 \rhd_\be
      N_2$.  By IH, $\fv{N_2} \subseteq \fv{M_2}$ and $\deg(M_2) =
      \deg(N_2)$.  By definition, $\fv{N} = \fv{M_1} \cup \fv{N_2}
      \subseteq \fv{M_1} \cup \fv{M_2} = \fv{M}$ and $\deg(M) =
      \deg(M_1) = \deg(N)$.
    \end{itemize}
    Case $r = \beta\eta$, by the $\beta$ and $\eta$ cases.  Case $r =
    {\wh}$, by the $\be$ case. \hfill $\Box$
  \end{enumerate}
\end{proof}

The next lemma is again needed in the proofs.
\begin{lemma}
  \label{deg+-f+}
  Let $i,p \geq 0$, $M, N, N_1, N_2, \dots, N_p \in {\cal M}$,
  $\su'\in \{\rhd_{\beta}^*,\rhd_{\eta}^*,\rhd_{\beta\eta}^*\}$ and
  $\su \in \{\rhd_{\beta},\rhd_{\eta},\rhd_{\beta\eta},\rhd_\wh,
  \rhd_{\beta}^*,\rhd_{\eta}^*,\rhd_{\beta\eta}^*,\rhd_{\wh}^*\}$.  We
  have:
  \begin{enumerate}
  \item
    \label{one'}
    $M^{+i}\in {\cal M}$ and $\deg(M^{+i}) = i::\deg(M)$ and $x^K$
    occurs in $M^{+i}$ iff $K = i::L$ and $x^L$ occurs in $M$.
  \item
    \label{one''}
    $M \diamond N$ iff $M^{+i} \diamond N^{+i}$.
  \item
    \label{one'''}
    Let ${\cal X} \subseteq {\cal M}$ then $\diamond {\cal X}$ iff
    $\diamond {\cal X}^{+i}$.
  \item
    \label{one}
    $(M^{+i})^{-i} = M$.
  \item
    \label{two}
    If $\diamond \{M\} \cup \{N_j$ / $j \in \{1, \dots, p\}\}$ then
    $(M[(x^{L_j}_j:=N_j)_p])^{+i} =
    M^{+i}[(x^{i::L_j}_j:=N_j^{+i})_p]$.
  \item
    \label{three}
    If $M \su N$, then $M^{+i} \su N^{+i}$.
  \item
    \label{four}
    If $\deg(M) = i :: L$, then:
    \begin{enumerate}
    \item $M = P^{+i}$ for some $P \in {\cal M}$, $\deg(M^{-i}) = L$
      and $(M^{-i})^{+i} = M$.
    \item If $\forall 1 \leq j \leq p, \deg(N_j) = i :: K_j$ and
      $\diamond \{M\} \cup \{N_j$ / $j \in \{1, \dots, p\}\}$ then
      $(M[(x^{i :: K_j}_j:=N_j)_p])^{-i} =
      M^{-i}[(x^{K_j}_j:=N^{-i}_j)_p]$.
    \item If $M \su N$ then $M^{-i}\su N^{-i}$.
    \end{enumerate}
  \item
    \label{five'}
    If $M \su N$, $P \su Q$ and $M \diamond P$ then $N \diamond Q$
  \item
    \label{five}
    If $M \su N^{+i}$, then there is $P \in {\cal M}$ such that $M=
    P^{+i}$ and $P \su N$.
  \item
    \label{six}
    If $M^{+i} \su N$, then there is $P \in {\cal M}$ such that $N=
    P^{+i}$ and $M \su P$.
  \item
    \label{lem:substlemma}
    If $y^{K} \not \in \fv{N} \cup \{x^{L}\}$, $\deg(P) = K$, $\deg(N)
    = L$, $\diamond \{M, N, P\}$ then\\
    $M[y^{K}:=P][x^{L}:=N] = M[x^{L}:=N][y^{K}:=P[x^{L}:=N]]$.
  \item
    \label{eight}
    If $M \su N$ and $\deg(P) = L$ and $\diamond \{M, N, P\}$, then
    $M[x^L:=P] \su N[x^L:=P]$.
  \item
    \label{nine}
    If $N \su' P$ and $\deg(N) = L = \deg(P)$ and $\diamond \{M, N, P\}$, then
    $M[x^L:=N] \su' M[x^L:=P]$.
  \item
    \label{ten}
    If $M \su' M'$, $P \su' P'$ and $\deg(P) = L$ and $\diamond \{M,
    M', P, P'\}$, then $M[x^L:=P] \su' M'[x^L:=P']$.
  \end{enumerate}
\end{lemma}

\begin{proof}
  \begin{enumerate}
  \item[\ref{one'}] We only prove the lemma by induction on $M$:
    \begin{itemize}
    \item If $M = x^L$ then $M^{+i} = x^{i::L} \in {\cal M}$ and
      $\deg(x^{i::L}) = i::L = i::\deg(x^{L})$.

    \item If $M = \l x^L. M_1$ then $M_1 \in {\cal M}$, $L \succeq
      \deg(M_1)$ and $M^{+i} = \l x^{i::L}. M_1^{+i}$.  By IH,
      $M_1^{+i} \in {\cal M}$ and $\deg(M_1^{+i}) = i::\deg(M_1)$ and
      $x^{K}$ occurs in $M_1^{+i}$ iff $K = i::K'$ and $y^{K'}$ occurs
      in $M_1$.  So $i::L \succeq i::\deg(M_1) =
      \deg(M_1^{+i})$. Hence, $\l x^{i::L}. M_1^{+i} \in {\cal
        M}$. Moreover, $\deg(M^{+i}) = \deg(M_1^{+i}) = i::\deg(M_1) =
      i::\deg(M)$.  If $y^{K}$ occurs in $M^{+i}$ then either $y^{K} =
      x^{i::L}$, so it is done because $x^{L}$ occurs in $M$. Or
      $y^{K}$ occurs in $M_1^{+i}$. By IH, $K = i::K'$ and $y^{K'}$
      occurs in $M_1$. So $y^{K'}$ occurs in $M$.  If $y^{K}$ occurs
      in $M$ then either $y^{K} = x^{L}$ and then $y^{i::K}$ occurs in
      $M^{+i}$. Or $y^{K}$ occurs in $M_1$. Then by IH, $y^{i::K}$
      occurs in $M_1^{+i}$. So, $y^{i::K}$ occurs in $M^{+i}$.

    \item If $M = M_1M_2$ then $M_1, M_2 \in {\cal M}$, $\deg(M_1)
      \preceq \deg(M_2)$, $M_1 \diamond M_2$ and $M^{+i} =
      M_1^{+i}M_2^{+i}$.  By IH, $M_1^{+i}, M_2^{+i} \in {\cal M}$,
      $\deg(M_1^{+i}) = i::\deg(M_1)$, $\deg(M_2^{+i}) =
      i::\deg(M_2)$, $y^{K}$ occurs in $M_1^{+i}$ iff $K = i::K'$ and
      $y^{K'}$ occurs in $M_1$, and $y^{K}$ occurs in $M_2^{+i}$ iff
      $K = i::K'$ and $y^{K'}$ occurs in $M_2$.  Let $x^{L} \in
      \fv{M_1^{+i}}$ and $x^{K} \in \fv{M_2^{+i}}$ then, using IH, $L
      = i::L'$, $K = i::K'$, $x^{L'}$ occurs in $M_1$ and $x^{K'}$
      occurs in $M_2$. Using $M_1 \diamond M_2$, we obtain $L' = K'$,
      so $L = K$. Hence, $M_1^{+i} \diamond M_2^{+i}$.  Because
      $\deg(M_1) \preceq \deg(M_2)$, then $\deg(M_1^{+i}) =
      i::\deg(M_1) \preceq i::\deg(M_2) = \deg(M_2^{+i})$.  So,
      $M^{+i}\in {\cal M}$.  Moreover, $\deg(M^{+1}) = \deg(M_1^{+i})
      = i::\deg(M_1) = i::\deg(M)$.  If $x^{L}$ occurs in $M^{+i}$
      then either $x^{L}$ occurs in $M_1^{+i}$ and using IH, $L =
      i::L'$ and $x^{L'}$ occurs in $M_1$, so $x^{L'}$ occurs in
      $M$. Or $x^{L}$ occurs in $M_2^{+i}$ and using IH, $L = i::L'$
      and $x^{L'}$ occurs in $M_2$, so $x^{L'}$ occurs in $M$.  If
      $x^{L}$ occurs in $M$ then either $x^{L}$ occurs in $M_1$ so by
      IH $x^{i::L}$ occurs in $M_1^{+i}$, hence $x^{i::L}$ occurs in
      $M^{+i}$. Or $x^{L}$ occurs in $M_2$ so by IH $x^{i::L}$ occurs
      in $M_2^{+i}$, hence $x^{i::L}$ occurs in $M^{+i}$.
    \end{itemize}

  \item[\ref{one''}] Assume $M \diamond N$. Let $x^{L} \in \fv{M^{+i}}$ and $x^{K}
    \in \fv{N^{+i}}$ then by lemma~\ref{deg+-f+}.\ref{one'}, $L =
    i::L'$, $K = i::K'$, $x^{L'} \in \fv{M}$ and $x^{K'} \in
    \fv{N}$. Using $M \diamond N$ we obtain $K' = L'$ and so $K = L$.

    Assume $M^{+i} \diamond N^{+i}$. Let $x^{L} \in \fv{M}$ and $x^{K}
    \in \fv{N}$, then by lemma~\ref{deg+-f+}.\ref{one'}, $x^{i::L} \in
    \fv{M^{+i}}$ and $x^{i::K} \in \fv{N^{+i}}$. Using $M^{+i}
    \diamond N^{+i}$ we obtain $i::K = i::L$ and so $K = L$.

  \item[\ref{one'''}] Let ${\cal X} \subseteq {\cal M}$.

    Assume $\diamond {\cal X}$. Let $M, N \in {\cal X}^{+i}$. Then by
    definition, $M = P^{+i}$ and $N = Q^{+i}$ such that $P, Q \in
    {\cal X}$. Because by hypothesis $P \diamond Q$ then by
    lemma~\ref{deg+-f+}.\ref{one''}, $M \diamond N$.

    Assume $\diamond {\cal X}^{+i}$. Let $M, N \in {\cal X}$ then
    $M^{+i}, N^{+i} \in {\cal X}^{+i}$. Because by hypothesis $M^{+i}
    \diamond N^{+i}$ then by lemma~\ref{deg+-f+}.\ref{one''}, $M
    \diamond N$.

  \item[\ref{one}] By lemma~\ref{deg+-f+}.\ref{one'}, $M^{+i} \in {\cal M}$ and
    $\deg(M^{+i}) = i::\deg(M)$. We prove the lemma by induction on
    $M$.
    \begin{itemize}
    \item Let $M = x^{L}$ then $M^{+i} = x^{i::L}$ and $(M^{+i})^{-i}
      = x^{L}$.
    \item Let $M = \l x^{L}. M_1$ such that $M_1 \in {\cal M}$ and $L
      \succeq \deg(M_1)$. Then, $(M^{+i})^{-i} = (\l
      x^{i::L}. M_1^{+i})^{-i} = \l x^{L}. (M_1^{+i})^{-i} =^{IH} \l
      x^{L}. M_1$.
    \item Let $M = M_1M_2$ such that $M_1, M_2 \in {\cal M}$, $M_1
      \diamond M_2$ and $\deg(M_1) \preceq \deg(M_2)$. Then,
      $(M^{+i})^{-i} = (M_1^{+i}M_2^{+i})^{-i} =
      (M_1^{+i})^{-i}(M_2^{+i})^{-i} =^{IH} M_1M_2$.
    \end{itemize}

  \item[\ref{two}] By~\ref{one'''}, $\diamond \{M^{+i}\} \cup \{N_j^{+i}$ / $j
    \in \{1, \dots, p\}\}$.  By lemma~\ref{degsub}.\ref{degsubthree},
    $M[(x^{L_j}_j:=N_j)_p]$ and $M^{+i}[(x^{i::L_j}_j:=N_j^{+i})_p]
    \in {\cal M}$.  By induction on $M$:
    \begin{itemize}
    \item Let $M = y^K$. If $\forall 1 \leq j \leq p, y^K \neq
      x^{L_j}_j$ then $y^K[(x^{L_j}_j:=N_j)_p] = y^K$. Hence
      $(y^K[(x^{L_j}_j:=N_j)_p])^{+i} = y^{i::K} =
      y^{i::K}[(x^{i::L_j}_j:=N_j^{+i})_p]$. If $\exists 1 \leq j \leq
      p, y^K = x^{L_j}_j$ then $y^K[(x^{L_j}_j:=N_j)_p] = N_j$. Hence
      $(y^K[(x^{L_j}_j:=N_j)_p])^{+i} = N^{+i}_j =
      y^{i::K}[(x^{i::L_j}_j:=N_j^{+i})_p]$.
    \item Let $M = \l y^K. M_1$. Then $M[(x^{L_j}_j:=N_j)_p] = \l
      y^K. M_1[(x^{L_j}_j:=N_j)_p]$ where $\forall 1 \leq j \leq p,
      y^K \not \in \fv{N_j} \cup \{x^{L_j}_j\}$. By
      lemma~\ref{degsub}.\ref{degsubone}, $\diamond \{M_1\} \cup
      \{N_j$ / $j \in \{1, \dots, p\}\}$.  By IH,
      $(M_1[(x^{L_j}_j:=N_j)_p])^{+i} =
      M_1^{+i}[(x^{i::L_j}_j:=N_j^{+i})_p]$. \\Hence,
      $(M[(x^{L_j}_j:=N_j)_p])^{+i} = \l y^{i::K}. (
      M_1[(x^{L_j}_j:=N_j)_p])^{+i} = $\\$ \l
      y^{i::K}. M_1^{+i}[(x^{i::L_j}_j:=N_j^{+i})_p] = (\l
      y^K. M_1)^{+i}[(x^{i::L_j}_j:=N_j^{+i})_p]$.
    \item Let $M = M_1M_2$. $M[(x^{L_j}_j:=N_j)_p] =
      M_1[(x^{L_j}_j:=N_j)_p]M_2[(x^{L_j}_j:=N_j)_p]$.  By
      lemma~\ref{degsub}.\ref{degsubone}, $\diamond \{M_1\} \cup
      \{N_j$ / $j \in \{1, \dots, p\}\}$ and $\diamond \{M_2\} \cup
      \{N_j$ / $j \in \{1, \dots, p\}\}$.  By IH,
      $(M_1[(x^{L_j}_j:=N_j)_p])^{+i} =
      M_1^{+i}[(x^{i::L_j}_j:=N_j^{+i})_p]$ and
      $(M_2[(x^{L_j}_j:=N_j)_p])^{+i} =
      M_2^{+i}[(x^{i::L_j}_j:=N_j^{+i})_p]$. \\Hence
      $(M[(x^{L_j}_j:=N_j)_p])^{+i} =
      (M_1[(x^{L_j}_j:=N_j)_p])^{+i}(M_2[(x^{L_j}_j:=N_j)_p])^{+i} =
      M_1^{+i}[(x^{i::L_j}_j:=N_j^{+i})_p]M_2^{+i}[(x^{i::L_j}_j:=N_j^{+i})_p]
      = M^{+i}[(x^{i::L_j}_j:=N_j^{+i})_p]$.
    \end{itemize}

  \item[\ref{three}] By lemma~\ref{deg+-f+}.\ref{one'}, if $M, N \in {\cal M}$ then
    $M^{+i}, N^{+i} \in {\cal M}$.
    \begin{itemize}
    \item Let $\su$ be $\rhd_\be$. By induction on $M \rhd_\be N$.
      \begin{itemize}
      \item Let $M = (\l x^L. M_1)M_2 \rhd_\be M_1[x^L := M_2] = N$
        where $\deg(M_2) = L$, then by lemma~\ref{deg+-f+}.\ref{one'},
        $\deg(M_2^{+i}) = i::L$ and $M^{+i} = (\l
        x^{i::L}. M_1^{+i})M_2^{+i} \rhd_\be M_1^{+i}[x^{i::L} :=
        M_2^{+i}] = (M_1[x^L := M_2])^{+i}$.
      \item Let $M = \l x^L. M_1 \rhd_\be \l x^L. N_1 = N$ such that
        $M_1 \rhd_\be N_1$.
        By IH, $M_1^{+i} \rhd_\be N_1^{+i}$, hence $M^{+i} = \l
        x^{i::L}. M_1^{+i} \rhd_\be \l x^{i::L} N_1^{+i} = N^{+i}$.
      \item Let $M = M_1 M_2 \rhd_\be N_1M_2 = N$ such that $M_1
        \rhd_\be N_1$. By IH, $M_1^{+i} \rhd_\be N_1^{+i}$, hence
        $M^{+i} = M_1 ^{+i}M_2^{+i} \rhd_\be N_1^{+i}M_2^{+i} =
        N^{+i}$.
      \item Let $M = M_1 M_2 \rhd_\be M_1N_2 = N$ such that $M_2
        \rhd_\be N_2$. By IH, $M_2^{+i} \rhd_\be N_2^{+i}$, hence
        $M^{+i} = M_1 ^{+i}M_2^{+i} \rhd_\be N_1^{+i}M_2^{+i} =
        N^{+i}$.
      \end{itemize}
    \item Let $\su$ be $\rhd_\be^*$. By induction on $\rhd_\be^*$
      using $\rhd_\be$.
    \item Let $\su$ be $\rhd_\eta$. We only do the base case. The
      inductive cases are as for $\rhd_\be$. Let $M = \l x^L. Nx^L
      \rhd_\eta N$ where $x^L \not \in \fv{N}$. By
      lemma~\ref{deg+-f+}.\ref{one'}, $x^{i::L} \not \in \fv{N^{+i}}$
      Then $M^{+i} = \l x^{i::L}. N^{+i}x^{i::L} \rhd_\eta N^{+i}$.
    \item Let $\su$ be $\rhd_\eta^*$. By induction on $\rhd_\eta^*$
      using $\rhd_\eta$.
    \item Let $\su$ be $\rhd_{\be\eta}$, $\rhd_{\be\eta}^*$,
      $\rhd_{\wh}$ or $\rhd_{\wh}^*$. By the previous items.
    \end{itemize}

  \item[\ref{four}]
    \begin{enumerate}
    \item By induction on $M$:
      \begin{itemize}
      \item Let $M = y^{i::L}$ then $y^L \in {\cal M}$ and
        $\deg((y^{i::L})^{-i}) = \deg(y^{L}) = L$ and
        $((y^{i::L})^{-i})^{+i} = y^{i::L}$.

      \item Let $M = \l y^K. M_1$ such that $M_1 \in {\cal M}$ and $K
        \succeq \deg(M_1)$.  Because $\deg(M_1) = \deg(M) = i::L$, by
        IH, $M_1 = P^{+i}$ for some $P \in {\cal M}$, $\deg(M_1^{-i})
        = L$ and $(M_1^{-i})^{+i} = M_1$.  Because, $K \succeq i::L$
        then $K = i::L::K'$ for some $K'$. Let $Q = \l y^{L::K'}. P$.
        Because $P =^{\ref{deg+-f+}.\ref{one}} (P^{+i})^{-i} =
        M_1^{-i}$, then $\deg(P) = L$.  Because $L \preceq L::K'$,
        then $Q \in {\cal M}$ and $Q^{+i} = M$.  Moreover,
        $\deg(M^{-i}) =^{\ref{deg+-f+}.\ref{one}} \deg(Q) = \deg(P) =
        L$ and $(M^{-i})^{+i} = P^{+i} = M$.

      \item Let $M = M_1M_2$ such that $M_1, M_2 \in {\cal M}$, $M_1
        \diamond M_2$ and $\deg(M_1) \preceq \deg(M_2)$.  Then
        $\deg(M) = \deg(M_1) \preceq \deg(M_2)$, so $\deg(M_2) =
        i::L::L'$ for some $L'$. By IH $M_1 = P_1^{+i}$ for some $P_1
        \in {\cal M}$, $\deg(M_1^{-i}) = L$ and $(M_1^{-i})^{+i} =
        M_1$. Again by IH, $M_2 = P_2^{+i}$ for some $P_2 \in {\cal
          M}$, $\deg(M_2^{-i}) = L::L'$ and $(M_2^{-i})^{+i} = M_2$.
        If $y^{K_1} \in \fv{P_1}$ and $y^{K_2} \in \fv{P_2}$, then by
        lemma~\ref{deg+-f+}.\ref{one'}, $K'_1 = i::K_1$, $K'_2 =
        i::K_2$, $x^{K'_1} \in \fv{M_1}$ and $x^{K'_2} \in
        \fv{M_2}$. Thus $K'_1 = K'_2$, so $K_1 = K_2$ and $P_1
        \diamond P_2$.  Because $\deg(P_1) = \deg(M_1^{-i}) = L
        \preceq L::L' = \deg(M_2^{-i}) = \deg(P_2)$ then $Q = P_1P_2
        \in {\cal M}$ and $Q^{+i} = (P_1P_2)^{+i} = P_1^{+i}P_2^{+i} =
        M$.  Moreover, $\deg(M^{-i}) =^{\ref{deg+-f+}.\ref{one}}
        \deg(Q) = \deg(P_1) = L$ and $(M^{-i})^{+i} = Q^{+i} = M$.
      \end{itemize}

    \item By the previous item, there exist $M', N'_1,
      \dots, N'_n \in {\cal M}$ such that $M = M'^{+i}$ and for all $j
      \in \{1, \dots, p\}$, $N_j = N'^{+i}_j$. So by
      lemma~\ref{deg+-f+}.\ref{one'''}, $\diamond \{M'\} \cup \{N'_j$
      / $j \in \{1, \dots, p\}\}$.  By lemma~\ref{deg+-f+}.\ref{one},
      $M^{-i} = M'$ and for all $j \in \{1, \dots, p\}$, $N_j^{-i} =
      N'_j$. So, $\diamond \{M^{-i}\} \cup \{N_j^{-i}$ / $j \in \{1,
      \dots, p\}\}$. By lemma~\ref{degsub}.\ref{degsubthree},
      $M[(x^{i::K_j}_j:=N_j)_p], M^{-i}[(x^{K_j}_j:=N^{-i}_j)_p] \in
      {\cal M}$ and $\deg(M[(x^{i::K_j}_j:=N_j)_p]) = \deg(M) = i::L$.
      We prove the result by induction on $M$:
      \begin{itemize}
      \item Let $M = y^{i::L}$. If $\forall 1 \leq j \leq p, y^{i::L} \neq
        x^{i::K_j}_j$ then $y^{i::L}[(x^{i::K_j}_j:=N_j)_p] = y^{i::L}$. Hence
        $(y^{i::L}[(x^{i::K_j}_j:=N_j)_p])^{-i} = y^L =
        y^L[(x^{K_j}_j:=N_j^{-i})_p]$. If $\exists
        1 \leq j \leq p, y^{i::L} = x^{i::K_j}_j$ then
        $y^{i::L}[(x^{i::K_j}_j:=N_j)_p] = N_j$. Hence
        $(y^{i::L}[(x^{i::K_j}_j:=N_j)_p])^{-i} = N^{-i}_j =
        y^L[(x^{K_j}_j:=N_j^{-i})_p]$.

      \item Let $M = \l y^K. M_1$ such that $M_1 \in {\cal M}$ and $K
        \succeq \deg(M_1)$. Then, $M[(x^{i::K_j}_j:=N_j)_p] = \l
        y^K. M_1[(x^{i::K_j}_j:=N_j)_p]$ where $\forall 1 \leq j \leq
        p, y^K \not \in \fv{N_j} \cup \{x^{i::K_j}_j\}$.  By
        lemma~\ref{degsub}.\ref{degsubone}, $\diamond \{M_1\} \cup
        \{N_j$ / $j \in \{1, \dots, p\}\}$. By definition $\deg(M) =
        \deg(M_1)$.  By IH, $(M_1[(x^{i::K_j}_j:=N_j)_p])^{-i} =
        M_1^{-i}[(x^{K_j}_j:=N_j^{-i})_p]$. Because $\deg(M_1) = i::L
        \preceq K$, $K = i::L::K'$ for some $K'$.\\
        Hence, $(M[(x^{i::K_j}_j:=N_j)_p])^{-i} = \l
        y^{L::K'}. (M_1[(x^{i::K_j}_j:=N_j)_p])^{-i} = \l
        y^{L::K'}. M_1^{-i}[(x^{K_j}_j:=N_j^{-i})_p] = (\l
        y^K. M_1)^{-i}[(x^{K_j}_j:=N_j^{-i})_p]$.

      \item Let $M = M_1M_2$ such that $M_1, M_2 \in {\cal M}$, $M_1
        \diamond M_2$ and $\deg(M_1) \preceq \deg(M_2)$. Then,
        $M[(x^{i::K_j}_j:=N_j)_p] =
        M_1[(x^{i::K_j}_j:=N_j)_p]M_2[(x^{i::K_j}_j:=N_j)_p]$. By
        lemma~\ref{degsub}.\ref{degsubone}, $\diamond \{M_1\} \cup
        \{N_j$ / $j \in \{1, \dots, p\}\}$ and $\diamond \{M_2\} \cup
        \{N_j$ / $j \in \{1, \dots, p\}\}$. By definition $\deg(M) =
        \deg(M_1) \preceq \deg(M_2)$. So $\deg(M_2) = i::L::L'$ for
        some $L'$.
        By IH, $(M_1[(x^{i::K_j}_j:=N_j)_p])^{-i} =
        M_1^{-i}[(x^{K_j}_j:=N_j^{-i})_p]$ and
        $(M_2[(x^{i::K_j}_j:=N_j)_p])^{-i} =
        M_2^{-i}[(x^{K_j}_j:=N_j^{-i})_p]$. Hence\\
        $(M[(x^{i::K_j}_j:=N_j)_p])^{-i} =
        (M_1[(x^{i::K_j}_j:=N_j)_p])^{-i}(M_2[(x^{i::K_j}_j:=N_j)_p])^{-i}
        $\\$=
        M_1^{-i}[(x^{K_j}_j:=N_j^{-i})_p]M_2^{-i}[(x^{K_j}_j:=N_j^{-i})_p]
        = M^{-i}[(x^{K_j}_j:=N_j^{-i})_p]$.
      \end{itemize}

    \item Using lemma~\ref{deg+-f+}.\ref{one},
      lemma~\ref{deg=} and the first item, we prove that $M^{-i},
      N^{-i} \in {\cal M}$.
      \begin{itemize}
      \item Let $\su$ be $\rhd_\be$. By induction on $M \rhd_\be N$.
        \begin{itemize}
        \item Let $M = (\l x^K. M_1)M_2 \rhd_\be M_1[x^K := M_2] = N$
          where $\deg(M_2) = K$.  Because $M \in {\cal M}$ then $M_1
          \in {\cal M}$.  Because $i::L = \deg(M) = \deg(M_1) \preceq
          K$, then $K = i::L::K'$.  By lemma~\ref{deg+-f+}.\ref{four},
          $\deg(M_2^{-i}) = L::K'$.  So $M^{-i} = (\l
          x^{L::K'}. M_1^{-i})M_2^{-i} \rhd_\be M_1^{-i}[x^{L::K'} :=
          M_2^{-i}] = (M_1[x^K := M_2])^{-i}$.

        \item Let $M = \l x^K. M_1 \rhd_\be \l x^K. N_1 = N$ such that
          $M_1 \rhd_\be N_1$.  Because $M \in {\cal M}$, $M_1 \in
          {\cal M}$ and $K \succeq \deg(M_1)$.  By definition $\deg(M)
          = \deg(M_1)$.  Because $i::L = \deg(M_1) \preceq K$, $K =
          i::L::K'$ for some $K'$.  By IH, $M_1^{-i} \rhd_\be
          N_1^{-i}$, hence $M^{-i} = \l x^{L::K'}. M_1^{-i} \rhd_\be
          \l x^{L::K'} N_1^{-i} = N^{-i}$.

        \item Let $M = M_1 M_2 \rhd_\be N_1M_2 = N$ such that $M_1
          \rhd_\be N_1$. Because $M \in {\cal M}$ then $M_1 \in {\cal
            M}$. By definition $\deg(M) = \deg(M_1) = i::L$.  By IH,
          $M_1^{-i} \rhd_\be N_1^{-i}$, hence $M^{-i} = M_1
          ^{-i}M_2^{-i} \rhd_\be N_1^{-i}M_2^{-i} = N^{-i}$.

        \item Let $M = M_1 M_2 \rhd_\be M_1N_2 = N$ such that $M_2
          \rhd_\be N_2$. Because $M \in {\cal M}$ then $M_2 \in {\cal
            M}$. By definition $\deg(M_2) \succeq \deg(M_1) = \deg(M)
          = i::L$. So $\deg(M_2) = i::L::L'$ for some $L'$.  By IH,
          $M_2^{-i} \rhd_\be N_2^{-i}$, hence $M^{-i} = M_1
          ^{-i}M_2^{-i} \rhd_\be N_1^{-i}M_2^{-i} = N^{-i}$.
        \end{itemize}

      \item Let $\su$ be $\rhd_\be^*$. By induction on $\rhd_\be^*$.
        using $\rhd_\be$.

      \item Let $\su$ be $\rhd_\eta$. We only do the base case. The
        inductive cases are as for $\rhd_\be$. Let $M = \l x^K. Nx^K
        \rhd_\eta N$ where $x^K \not \in \fv{N}$.  Because $i::L =
        \deg(M) = \deg(N) \preceq K$, then $K = i::L::K'$ for some
        $K'$.  By lemma~\ref{deg+-f+}.\ref{four}, $N = N'^{+i}$ for
        some $N' \in {\cal M}$. By lemma~\ref{deg+-f+}.\ref{four}, $N'
        = N^{-i}$.  By lemma~\ref{deg+-f+}.\ref{one'}, $x^{L::K'} \not
        \in \fv{N^{-i}}$.  Then $M^{-i} = \l
        x^{L::K'}. N^{-i}x^{L::K'} \rhd_\eta N^{-i}$.

      \item Let $\su$ be $\rhd_\eta^*$. By induction on $\rhd_\eta^*$
        using $\rhd_\eta$.

      \item Let $\su$ be $\rhd_{\be\eta}$, $\rhd_{\be\eta}^*$,
        $\rhd_{\wh}$ or $\rhd_{\wh}^*$.  By the previous items.
      \end{itemize}
    \end{enumerate}

  \item[\ref{five'}] Let $x^L \in \fv{N} \subseteq^{\ref{deg=}}
    \fv{M}$ and $X^K \in \fv{Q} \subseteq^{\ref{deg=}} \fv{P}$, since
    $M \diamond P$, $L = K$. Hence $N \diamond Q$.

  \item[\ref{five}] By lemma~\ref{deg+-f+}.\ref{one'}, $\deg(N^{+i}) =
    i::\deg(N)$. By lemma~\ref{deg=}, $\deg(M) = \deg(N^{+i})$. By
    lemma~\ref{deg+-f+}.\ref{four}, $M = M'^{+i}$ such that $M' \in
    {\cal M}$. By lemma~\ref{deg+-f+}.\ref{four}, $M'
    =^{\ref{deg+-f+}.\ref{one}} (M'^{+i})^{-i} = M^{-i} \su
    (N^{+i})^{-i} =^{\ref{deg+-f+}.\ref{one}} N$.

  \item[\ref{six}] By lemma~\ref{deg+-f+}.\ref{one'}, $\deg(M^{+i}) =
    i::\deg(M)$. By lemma~\ref{deg=}, $\deg(M^{+i}) = \deg(N)$. By
    lemma~\ref{deg+-f+}.\ref{four}, $N = N'^{+i}$ such that $N' \in
    {\cal M}$. By lemma~\ref{deg+-f+}.\ref{four}, $M
    =^{\ref{deg+-f+}.\ref{one}} (M^{+i})^{-i} \su N^{-i} =
    (N'^{+i})^{-i} =^{\ref{deg+-f+}.\ref{one}} N'$.

  \item[\ref{lem:substlemma}] By lemma~\ref{degsub}.\ref{degsubthree},
    $M[y^{K}:=P] \in {\cal M}$. Let us now prove $\diamond
    \{M[y^{K}:=P], N\}$. Let $z^{R} \in \fv{M[y^{K}:=P]}$ and $z^{R'}
    \in \fv{N}$ then $z^{R} \in \fv{M}$ or $z^{R} \in \fv{P}$. In both
    cases, because $M \diamond N$ and $P \diamond N$, we obtain $R =
    R'$. So by lemma~\ref{degsub}.\ref{degsubthree},
    $M[y^{K}:=P][x^{L}:=N] \in {\cal M}$.

    By lemma~\ref{degsub}.\ref{degsubthree}, $M[x^{L}:=N], P[x^{L}:=N]
    \in {\cal M}$ and $\deg(P[x^{L}:=N]) = \deg(P) = K$. Let us now
    prove that $\diamond \{M[x^{L}:=N], P[x^{L}:=N]\}$. Let $z^{R} \in
    \fv{M[x^{L}:=N}$ and $z^{R'} \in \fv{P[x^{L}:=N]}$ then either
    $z^{R} \in \fv{M}$ or $z^{R} \in \fv{N}$ and either $z^{R'} \in
    \fv{P}$ or $z^{R'} \in \fv{N}$. In all of the four cases, because
    by hypotheses and lemma~\ref{degsub}.\ref{lem:refl+symm}, $M
    \diamond P$, $M \diamond N$, $N \diamond P$ and $N \diamond N$, we
    obtain $R = R'$. So by lemma~\ref{degsub}.\ref{degsubthree},
    $M[x^{L}:=N][y^{K}:=P[x^{L}:=N]] \in {\cal M}$.

    We prove this lemma by induction on the structure of $M$.
    \begin{itemize}
    \item Let $M = z^{R}$.
      \begin{itemize}
      \item If $z^{R} = y^{K}$ then $M[y^{K}:=P][x^{L}:=N] =
        P[x^{L}:=N] = M[y^{K}:=P[x^{L}:=N]] =
        M[x^{L}:=N][y^{K}:=P[x^{L}:=N]]$.
      \item Else
        \begin{itemize}
        \item If $M = x^{L}$ then $M[y^{K}:=P][x^{L}:=N] = M[x^{L}:=N]
          = N = N[y^{K}:=P[x^{L}:=N]] =
          M[x^{L}:=N][y^{K}:=P[x^{L}:=N]]$.
        \item Else $M[y^{K}:=P][x^{L}:=N] = M[x^{L}:=N] = M =
          M[y^{K}:=P[x^{L}:=N]] = M[x^{L}:=N][y^{K}:=P[x^{L}:=N]]$.
        \end{itemize}
      \end{itemize}
    \item Let $M = \l z^{R} M_1$ such that $R \succeq \deg(M_1)$ and
      $M_1 \in {\cal M}$. By lemma~\ref{degsub}.\ref{degsubone},
      $\diamond \{M_1, N, P\}$. Then, $M[y^{K}:=P][x^{L}:=N] = \l
      z^{R}. M_1[y^{K}:=P][x^{L}:=N] =^{IH} \l
      z^{R}. M_1[x^{L}:=N][y^{K}:=P[x^{L}:=N]] =
      M[x^{L}:=N][y^{K}:=P[x^{L}:=N]]$ such that $z^{R} \not \in
      \fv{N} \cup \fv{P} \cup \{y^{K}, x^{L}\}$.
    \item Let $M = M_1M_2$ such that $M_1, M_2 \in {\cal M}$,
      $\deg(M_1) \preceq \deg(M_2)$ and $M_1 \diamond M_2$. By
      lemma~\ref{degsub}.\ref{degsubone}, $\diamond \{M_1, N, P\}$ and
      $\diamond \{M_2, N, P\}$. Then, $M[y^{K}:=P][x^{L}:=N] =
      M_1[y^{K}:=P][x^{L}:=N]M_2[y^{K}:=P][x^{L}:=N] =^{IH}
      M_1[x^{L}:=N][y^{K}:=P[x^{L}:=N]]M_2[x^{L}:=N][y^{K}:=P[x^{L}:=N]]
      = M[x^{L}:=N][y^{K}:=P[x^{L}:=N]]$.
    \end{itemize}

  \item[\ref{eight}] By lemma~\ref{degsub}.\ref{degsubthree} and using
    the hypothesis, we obtain $M[x^L:=P], N[x^L:=P] \in {\cal M}$.
    \begin{itemize}
    \item Let $\su = \rhd_{\beta}$. We prove the result by induction
      on $M \rhd_{\beta} N$.
      \begin{itemize}
      \item Let $M = (\l y^{K}. M_1)M_2 \rhd_{\beta} M_1[y^{K}:=M_2] =
        N$ such that $\deg(M_2) = K$. Then $M[x^{L}:=P] = (\l
        y^{K}. M_1[x^{L}:=P])M_2[x^{L}:=P]$ and $N[x^{L}:=P]
        =^{\ref{deg+-f+}.\ref{lem:substlemma}}
        M_1[x^{L}:=P][y^{K}:=M_2[x^{L}:=P]]$ such that $y^{K} \not \in
        \fv{P} \cup \{x^{L}\}$. By
        lemma~\ref{degsub}.\ref{degsubthree}, $\deg(M_2[x^{L}:=P]) =
        \deg(M_2) = K$. So, $M[x^{L}:=P] \rhd_{\beta} N[x^{L}:=P]$.
      \item Let $M = \l y^{K}. M_1 \rhd_{\beta} \l y^{K}. N_1 = N$
        such that $M_1 \rhd_{\beta} N_1$.  Then $M[x^{L}:=P] = \l
        y^{K}. M_1[x^{L}:=P]$ and $N[x^{L}:=P] = \l
        y^{K}. N_1[x^{L}:=P]$ such that $y^{K} \not \in \fv{P} \cup
        \{x^{L}\}$.  By lemma~\ref{degsub}.\ref{degsubone}, $\diamond
        \{M_1, N_1, P\}$.  By IH, $M_1[x^{L}:=P] \rhd_{\beta}
        N_1[x^{L}:=P]$. So, $M[x^{L}:=P] \rhd_{\beta} N[x^{L}:=P]$.
      \item Let $M = M_1M_2 \rhd_{\beta} N_1M_2 = N$ such that $M_1
        \rhd_{\beta} N_1$.  By lemma~\ref{degsub}.\ref{degsubone},
        $\diamond \{M_1, N_1, P\}$. By IH, $M_1[x^{L}:=P] \rhd_{\beta}
        N_1[x^{L}:=P]$.  So, $M[x^{L}:=P] \rhd_{\beta} N[x^{L}:=P]$.
      \item Let $M = M_1M_2 \rhd_{\beta} M_1N_2 = N$ such that $M_2
        \rhd_{\beta} N_2$.  By lemma~\ref{degsub}.\ref{degsubone},
        $\diamond \{M_2, N_2, P\}$. By IH, $M_2[x^{L}:=P] \rhd_{\beta}
        N_2[x^{L}:=P]$.  So, $M[x^{L}:=P] \rhd_{\beta} N[x^{L}:=P]$.
      \end{itemize}

    \item Let $\su = \rhd_{\eta}$. We only prove the base case. The
      other cases are similar as the ones for $\rhd_{\beta}$. Let $M =
      \l y^{K}. Ny^{K} \rhd_{\eta} N$ such that $y^{K} \not \in
      \fv{N}$. Then $M[x^{L}:=P] = \l y^{K}. N[x^{L}:=P]y^{K}$ such
      that $y^{K} \not \in \fv{P} \cup \{x^{L}\}$. So $y^{K} \not \in
      \fv{N[x^{L}:=P]}$. Hence, $M[x^{L}:=P] \rhd_{\eta} N[x^{L}:=P]$.

    \item The other cases are based on the two previous ones.
    \end{itemize}

  \item[\ref{nine}] By lemma~\ref{degsub}.\ref{degsubthree} and using
    the hypothesis, we obtain $M[x^L:=P], M[x^L:=N] \in {\cal M}$. We
    prove the result by induction on the structure of $M$.
    \begin{itemize}
    \item Let $M = y^{K}$.
      \begin{itemize}
      \item If $y^{K} = x^{L}$ then $M[x^{L}:=P] = P \su' N =
        M[x^{L}:=N]$.
      \item Else, $M[x^{L}:=P] = M \su' M = M[x^{L}:=N]$.
      \end{itemize}
    \item Let $M = \l y^{K}. M_1$ such that $K \succeq \deg(M_1)$ and
      $M_1 \in {\cal M}$. Then $M[x^{L}:=P] = \l y^{K}. M_1[x^{L}:=P]$
      and $M[x^{L}:=N] = \l y^{K}. M_1[x^{L}:=N]$ such that $y^{K}
      \not \in \fv{P} \cup \fv{N} \cup \{x^{L}\}$. By
      lemma~\ref{degsub}.\ref{degsubone}, $\diamond \{M_1, N,
      P\}$. By IH, $M_1[x^L:=N] \su' M_1[x^L:=P]$. So, $M[x^L:=N] \su'
      M[x^L:=P]$.
    \item Let $M = M_1M_2$ such that $M_1, M_2 \in {\cal M}$, $M_1
      \diamond M_2$ and $\deg(M_1) \preceq \deg(M_2)$.  By
      lemma~\ref{degsub}.\ref{degsubone}, $\diamond \{M_1, N, P\}$ and
      $\diamond \{M_2, N, P\}$. By IH, $M_1[x^L:=N] \su' M_1[x^L:=P]$
      and $M_2[x^L:=N] \su' M_2[x^L:=P]$.  By
      lemma~\ref{degsub}.\ref{degsubthree}, $M_1[x^L:=N], M_2[x^L:=N],
      M_1[x^L:=P], M_2[x^L:=P] \in {\cal M}$ and $\deg(M_1[x^L:=N]) =
      \deg(M_1) \preceq \deg(M_2) = \deg(M_2[x^L:=N])$ and
      $\deg(M_1[x^L:=P]) = \deg(M_1) \preceq \deg(M_2) =
      \deg(M_2[x^L:=N])$ and $\deg(M_1[x^L:=P]) = \deg(M_1) \preceq
      \deg(M_2) = \deg(M_2[x^L:=P])$.  By
      lemma~\ref{degsub}.\ref{degsubfour}, $M_1[x^L:=N] \diamond
      M_2[x^L:=N]$ and $M_1[x^L:=P] \diamond M_2[x^L:=N]$ and
      $M_1[x^L:=P] \diamond M_2[x^L:=P]$. So $M_1[x^L:=N]M_2[x^L:=N],
      M_1[x^L:=P]M_2[x^L:=N], M_1[x^L:=P]M_2[x^L:=P] \in {\cal M}$.
      So $M_1[x^L:=N]M_2[x^L:=N] \su' M_1[x^L:=P]M_2[x^L:=N]$ and
      $M_1[x^L:=P]M_2[x^L:=N] \su' M_1[x^L:=P]M_2[x^L:=P]$.  Hence,
      $M[x^L:=N] \su' M[x^L:=P]$.
    \end{itemize}

  \item[\ref{ten}] By lemma~\ref{deg+-f+}.\ref{eight}, $M[x^L:=P] \su'
    M'[x^L:=P]$. By lemma~\ref{deg+-f+}.\ref{nine}, $M'[x^L:=P] \su'
    M'[x^L:=P']$. So, $M[x^L:=P] \su' M'[x^L:=P']$.
\hfill $\Box$
   \end{enumerate}
\end{proof}

Next we give a lemma that will be used in the rest of the article.
\begin{lemma}
  \label{trivial1}
  \begin{enumerate}
  \item
    \label{mynewme}
    If $M[y^L:=x^L] \rhd_\be N$ then $M \rhd_\be N'$ where $N =
    N'[y^L:=x^L]$.
  \item
    \label{newme}
    If $M[y^L:=x^L]$ is $\be$-normalising then $M$ is
    $\be$-normalising.
  \item
    \label{trivial1two}
    Let $k \geq 1$. If $M x^{L_1}_1...x^{L_k}_k$ is $\be$-normalising,
    then $M$ is $\be$-normalising.
  \item
    \label{trivial1three}
    Let $k \geq 1$, $1 \leq i \leq k$, $l \geq 0$, $x^{L_i}_i
    N_1...N_l$ be in normal form and $M$ be closed. If $M
    x^{L_1}_1...x^{L_k}_k \rhd_\be^* x^{L_i}_i N_1...N_l$, then for
    some $m \geq i$ and $n \leq l$, $M \rhd_\be^* \l x^{L_1}_1. ...\l
    x^{L_m}_m. x^{L_i}_i M_1...M_n$ where $n + k = m + l$, $M_j
    \simeq_\be N_j$ for every $1 \leq j \leq n$ and $N_{n+j}
    \simeq_\be x^{L_{m+j}}_{m+j}$ for every $1 \leq j \leq k-m$.
  \end{enumerate}
\end{lemma}

\begin{proof}
  \begin{itemize}
  \item[\ref{mynewme}.]  By induction on $M[y^L:=x^L] \rhd_\be N$.
  \item[\ref{newme}.]  Immediate by \ref{mynewme}.
  \item[\ref{trivial1two}.]  By induction on $k \geq 1$. We only prove
    the basic case. The proof is by cases.
    \begin{itemize}
    \item If $M\, x^{L_1}_1\rhd^*_\be M'\, x^{L_1}_1$ where $M'\,
      x^{L_1}_1$ is in $\be$-normal form and $M \rhd^*_\be M'$ then
      $M'$ is in $\be$-normal form and $M$ is $\be$-normalising.
    \item If $M\, x^{L_1}_1 \rhd^*_\be (\l y^{L_1}.N)\, x^{L_1}_1
      \rhd_\be N[y^{L_1}:=x^{L_1}_1] \rhd^*_\be P$ where $P$ is in
      $\be$-normal form and $M \rhd^*_\be \l y^{L_1}.N$ then
      by~\ref{newme}, $N$ has a $\be$-normal form and so, $\l
      y^{L_1}.N$ has a $\be$-normal form. Hence, $M$ has a
      $\be$-normal form.
    \end{itemize}
  \item[\ref{trivial1three}.]  By \ref{trivial1two}, $M$ is
    $\be$-normalising and, since $M$ is
    closed, its $\be$-normal form is \\
    $\l x^{L_1}_1. ...\l x^{L_m}_m. x^{L_p}_p M_1 ... M_n$ for $n, m
    \geq 0$ and $1 \leq p \leq m$.\\
    Since by theorem~\ref{confluenceofbetaeta}, $x^{L_i}_i N_1...N_l
    \simeq_{\beta} (\l x^{L_1}_1. ...\l x^{L_m}_m. x^{L_p}_p M_1
    ... M_n)x^{L_1}_1...x^{L_k}_k$ then $m \leq k$, $x^{L_i}_i \;
    N_1...N_l \simeq_\be x^{L_p}_p M_1 ... M_n
    x^{L_{m+1}}_{m+1}...x^{L_k}_k$. Hence, $n \leq l$, $i = p \leq m$,
    $l = n + k - m$, for every $1 \leq j \leq n$, $M_j \simeq_\be N_j$
    and for every $1 \leq j \leq k-m$, $N_{n+j} \simeq_\be
    x^{n_{m+j}}_{m+j}$.
\hfill $\Box$
  \end{itemize}
\end{proof}

\subsection{Confluence of $\rhd_{\be}^*$, $\rhd_{\wh}^*$ and
  $\rhd_{\be\eta}^*$}
\label{confproof}

In this section we establish the confluence of $\rhd_{\be}^*$,
$\rhd_{\wh}^*$ and $\rhd_{\be\eta}^*$ using the standard parallel
reduction method for $\rhd_{\be}^*$ and $\rhd_{\be\eta}^*$ .

\begin{definition}\label{newdefrho}
Let $r \in \{\be, \beta\eta\}$. We define on ${\cal M}$ the binary
relation $\myrhor$ by:
\begin{itemize}
\item $M \myrhor M$
\item If $M \myrhor M'$ then $\l x^L.M \myrhor \l x^L.M'$.
\item If $M \myrhor M'$, $N \myrhor N'$ and $M \diamond N$ and
  $\deg(M) \succeq \deg(N)$ then $M N \myrhor M' N'$
\item If $M \myrhor M'$, $N \myrhor N'$, $\deg(N) = L \succeq \deg(M)$
  and $M \diamond N$, then $(\l x^L. M) N \myrhor M'[x^n := N']$
\item If $M \myrhobet M'$, $x^L \diamond M$ and $L \succeq \deg(M)$
  then $\l x^L. M x^L \myrhobet M'$
\end{itemize}
We denote the transitive closure of $\myrhor$ by $\myrrhor$.  When $M
\myrhor N$ (resp.\ $M \myrrhor N$), we can also write $N \leftrhor M$
(resp.\ $N \lleftrhor M$).  If $R, R' \in \{\myrhor, \myrrhor,
\leftrhor, \lleftrhor\}$, we write $M_1\ R\ M_2\ R'\ M_3$ instead of
$M_1\ R\ M_2$ and $M_2\ R'\ M_3$.
\end{definition}

\begin{lemma}
  \label{newrho2}
  Let $M \in {\cal M}$.
  \begin{enumerate}
  \item
    \label{newrho2one}
    If $M \rhd_r M'$, then $M \myrhor M'$.
  \item
    \label{newrho2two}
    If $M \myrhor M'$, then $M' \in {\cal M}$, $M \rhd_r^* M'$,
    $\fv{M'} \subseteq \fv{M}$ and $\deg(M) = \deg(M')$.
  \item
    \label{newrho2three}
    If $M \myrhor M'$, $N \myrhor N'$ and $M\diamond N$ then
    $M'\diamond N'$
  \end{enumerate}
\end{lemma}

\begin{proof}
\ref{newrho2one}.\ By induction on the derivation $M\rhd_r M'$.
\ref{newrho2two}.\ By induction on the derivation of $M \myrhor
M'$ using  theorem~\ref{deg=} and lemma~\ref{deg+-f+}.
\ref{newrho2three}.\ Let $x^L\in \fv{M'}$ and $x^K\in \fv{N'}$. By
\ref{newrho2two}.,  $\fv{M'} \subseteq \fv{M}$ and $\fv{N'} \subseteq
\fv{N}$. Hence, since $M\diamond N$, $L =K$, so $M'\diamond N'$.
\hfill $\Box$
\end{proof}

\begin{lemma}\label{newrho1}
 Let $M, N \in {\cal M}$, $M\diamond N$ and $N \myrhor N'$.  We have:
\begin{enumerate}
\item\label{newrho1one}
  $M[x^L := N] \myrhor M[x^L :=N']$.
\item\label{newrho1two}
  If $M \myrhor M'$ and $\deg(N) = L$, then $M[x^L := N] \myrhor
  M'[x^L :=N']$.
\end{enumerate}
\end{lemma}

\begin{proof}
\ref{newrho1one}.\
By induction on $M$:
\begin{itemize}
\item Let $M = y^K$. If $y^K = x^L$, then $M[x^L := N] = N$, $M[x^L
  :=N'] = N'$ and by hypothesis, $N \myrhor N'$. If $y^K \not = x^L$,
  then $M[x^L := N] = M$, $M[x^L :=N'] = M$ and by definition, $M
  \myrhor M$.
\item Let $M = \l y^K. M_1$. $M[x^L := N] = \l y^K. M_1[x^L := N]$ and
  since $M_1 \diamond N$, by IH, $M_1[x^L := N] \myrhor M_1[x^L :=
  N']$ and so $\l y^K. M_1[x^L := N] \myrhor \l y^K. M_1[x^L := N']$
\item Let $M = M_1M_2$. $M[x^L := N] = M_1[x^L := N]M_2[x^L := N]$ and
  since $M_1 \diamond N$ and $M_2 \diamond N$, by IH, $M_1[x^L := N]
  \myrhor M_1[x^L := N']$ and $M_2[x^L := N] \myrhor M_2[x^L :=
  N']$. By lemma~\ref{degsub}.\ref{degsubfour}, $M_1[x^L := N]
  \diamond M_2[x^L := N]$, so $M_1[x^L := N]M_2[x^L := N] \myrhor
  M_1[x^L := N']M_2[x^L := N']$.
\end{itemize}
\ref{newrho1two}.\ By induction on  $M \myrhor M'$.
\begin{itemize}
\item If $M = M'$, then \ref{newrho1one}..
\item If $\l y^K.M \myrhor \l y^K.M'$ where $M \myrhor M'$, then by
  IH, $M[x^L := N] \myrhor M'[x^L := N']$. Hence $(\l y^K. M)[x^L :=
  N] = \l y^K. M[x^L := N] \myrhor \l y^K. M'[x^L := N'] = (\l
  y^K. M')[x^L := N']$ where $y^K \not \in \fv{N'} \subseteq \fv{N}$.
\item  If $P Q \myrhor P' Q'$ where $P \myrhor P'$, $Q \myrhor
  Q'$ and $P \diamond Q$, then by IH, $P[x^L := N] \myrhor P'[x^L :=
  N']$ and $Q[x^L := N] \myrhor Q'[x^L := N']$. By
  lemma~\ref{degsub}.\ref{degsubfour}, $P[x^L := N]
  \diamond Q[x^L := N]$, so $P[x^L := N]Q[x^L := N] \myrhor
  P'[x^L := N']Q'[x^L := N']$.
\item $(\l y^K. P)Q \myrhor P'[y^K := Q']$ where $P \myrhor P'$,
  $Q \myrhor Q'$, $P \diamond Q$ and $\deg(Q) = K$, then by IH, $P[x^L
  := N] \myrhor P'[x^L := N']$,
  $Q[x^L := N] \myrhor Q'[x^L := N']$. Moreover, $((\l y^K. P)Q)[x^L := N] =
  (\l y^K. P)[x^L := N]Q[x^L := N] = \l y^K. P[x^L := N]Q[x^L := N]$
  where $y^K \not \in \fv{N'} \subseteq \fv{N}$. By
  lemma~\ref{degsub}.\ref{degsubfour}, $P[x^L := N]
  \diamond Q[x^L := N]$ and by lemma~\ref{degsub}.\ref{degsubthree}
  $\deg(Q) =\deg(Q[x^L := N])$ so $\l y^K. P[x^L :=
  N]Q[x^L := N] \myrhor P'[x^L := N'][y^K := Q'[x^L := N']] = P'[y^K :=
  Q'][x^L := N']$.
\item If  $\l y^K. M y^K \myrhobet M'$ where $M \myrhobet M'$, $K
  \succeq \deg(M)$ and $\forall K \in {{\cal L}_{\mathbb N}}, y^K
  \not \in \fv{M}$, then by IH $M[x^L := N]
  \myrhobet M'[x^L :=N']$. Moreover, $(\l y^K. M y^K)[x^L := N] = \l
  y^K. M[x^L := N] y^K[x^L := N] = \l y^K. M[x^L := N] y^K$ where
  $\forall K \in {{\cal L}_{\mathbb N}}, y^K
  \not \in \fv{N'} \subseteq \fv{N}$. Since by
  lemma~\ref{degsub}.\ref{degsubthree} $\deg(M) =\deg(M[x^L := N])$,
  $\l y^K. M[x^L := N] y^K \myrhobet M'[x^L := N']$.
\hfill $\Box$
\end{itemize}
\end{proof}

\begin{lemma}\label{rhorem}
\begin{enumerate}
\item\label{rhoremone}
  If $x^L \myrhor N$, then $N = x^L$.
\item\label{rhoremthree}
  If $\l x^L.P \myrhobet N$ then one of the following holds:
  \begin{itemize}
  \item
    $N = \l x^L.P'$ where $P \myrhobet  P'$.
  \item
    $P = P' x^L$ where $\forall L \in {{\cal L}_{\mathbb N}}, x^L
  \not \in \fv{P'}$, $L \succeq \deg(P')$ and $P' \myrhobet  N$.
  \end{itemize}
\item\label{rhoremthree'}
  If $\l x^L.P \myrhobe N$ then $N = \l x^L.P'$ where $P \myrhobe
  P'$.
\item\label{rhoremfour}
  If $PQ \myrhor  N$, then one of the following holds:
  \begin{itemize}
  \item
    $N = P'Q'$, $P \myrhor P'$, $Q \myrhor Q'$ and $P\diamond Q$.
  \item
    $P = \l x^L.P'$, $N = P''[x^L:=Q']$, $P' \myrhor P''$, $Q  \myrhor
    Q'$, $P'\diamond Q$ and $\deg(Q) = L$.
  \end{itemize}
\end{enumerate}
\end{lemma}

\begin{proof}
\ref{rhoremone}.\ By induction on the derivation $x^L \myrhor N$.\\
\ref{rhoremthree}.\ By induction on the derivation $\l x^L.P \myrhobet N$.\\
\ref{rhoremthree'}.\ By induction on the derivation $\l x^L.P \myrhobe N$.\\
\ref{rhoremfour}.\ By induction on the derivation  $PQ \myrhor N$.
\hfill $\Box$
\end{proof}

\begin{lemma}\label{newrho3}
Let $M, M_1, M_2 \in {\cal M}$.
\begin{enumerate}
\item\label{newrho3one}
If $M_2 \leftrhor M \myrhor M_1$, then there is $M' \in {\cal M}$ such that
$M_2 \myrhor M'  \leftrhor M_1$.
\item\label{newrho3two}
If $M_2 \lleftrhor M \myrrhor M_1$, then there is $M' \in {\cal M}$ such that
$M_2 \myrrhor M'  \lleftrhor M_1$.
\end{enumerate}\end{lemma}

\begin{proof}
\ref{newrho3one}.\ By induction on $M$:
\begin{itemize}
\item Let $r = \be\eta$:
  \begin{itemize}
  \item
    If $M = x^L$, by lemma~\ref{rhorem}, $M_1 = M_2 = x^L$.  Take $M' =
    x^L$.

  \item If $N_2P_2 \leftrhobet NP \myrhobet N_1P_1$ where $N_2
    \leftrhobet N \myrhobet N_1$, $P_2 \leftrhobet P \myrhobet P_1$ and
    $N \diamond P$ then, by IH, $\exists N', P'$ such that $N_2
    \myrhobet N'\leftrhobet N_1$ and $P_2 \myrhobet P'\leftrhobet P_1$.
    By lemma~\ref{newrho2}.\ref{newrho2three}, $N_1 \diamond P_1$ and
    $N_2 \diamond P_2$, hence $N_2P_2 \myrhobet N'P' \leftrhobet
    N_1P_1$.

  \item If $(\l x^L.P_1) Q_1 \leftrhobet (\l x^L.P) Q \myrhobet
    P_2[x^L:=Q_2]$ where $\l x^L. P\myrhobet \l x^L.P_1$, $P \myrhobet
    P_2$, $Q_1 \leftrhobet Q \myrhobet Q_2$, $\deg(Q) = L$, $(\l
    x^L. P)\diamond Q$ and $P\diamond Q$ then, by
    lemma~\ref{rhorem}, $P \myrhobet P_1$. By IH, $\exists
    P', Q'$ such that $P_1 \myrhobet P'\leftrhobet P_2$ and $Q_1
    \myrhobet Q'\leftrhobet Q_2$. By
    lemma~\ref{newrho2}.\ref{newrho2two}, $\deg(Q_1) = \deg(Q_2) = \deg(Q) =
    L$. By lemma~\ref{newrho2}.\ref{newrho2three}, $P_1 \diamond
    Q_1$. Hence, $(\l x^L.P_1)Q_1  \myrhobet P'[x^L:=Q']$.

    Moreover, since $P_2 \myrhobet P'$, $Q_2 \myrhobet Q'$,
    $\deg(Q_2) = L$ and by lemma~\ref{newrho2}.\ref{newrho2three},
    $P_2 \diamond Q_2$, then, by lemma~\ref{newrho1}.\ref{newrho1two},
    $P_2[x^L:=Q_2] \myrhobet P'[x^L:=Q']$.

  \item If  $P_1[x^L:= Q_1] \leftrhobet (\l x^L. P) Q \myrhobet
    P_2[x^L:=Q_2]$ where $P_1 \leftrhobet P \myrhobet P_2$, $Q_1
    \leftrhobet Q \myrhobet Q_2$, $\deg(Q) = L$ and $P\diamond Q$,
    then, by IH,  $\exists P', Q'$ where $P_1
    \myrhobet P'\leftrhobet P_2$ and $Q_1 \myrhobet Q'\leftrhobet
    Q_2$. By lemma~\ref{newrho2}.\ref{newrho2two}, $\deg(Q_1) =
    \deg(Q_2) = \deg(Q) = L$. By
    lemma~\ref{newrho2}.\ref{newrho2three}, $P_1 \diamond
    Q_1$ and $P_2 \diamond Q_2$. Hence, by
    lemma~\ref{newrho1}.\ref{newrho1two}, $P_1[x^L:=Q_1]\myrhobet
    P'[x^L:=Q'] \leftrhobet P_2[x^L:=Q_2]$.

  \item
    If $\l x^L. N_2 \leftrhobet \l x^L. N \myrhobet\l x^L. N_1$ where $N_2
    \leftrhobet N \myrhobet N_1$, by IH, there is $N'$ such that
    $N_2 \myrhobet N'  \leftrhobet N_1$. Hence, $\l x^L. N_2
    \myrhobet \l x^L. N' \leftrhobet \l x^L. N_1$.

  \item
    If $M_1 \leftrhobet \l x^L. P x^L \myrhobet  M_2$ where $\forall L
    \in {{\cal L}_{\mathbb N}}, x^L \not \in \fv{P}$, $L \succeq
    \deg(P)$ and $M_1 \leftrhobet P \myrhobet M_2$, then, by IH, there
    is $M'$ such that $M_2 \myrhobet M' \leftrhobet M_1$.

  \item
    If $M_1 \leftrhobet \l x^L. P x^L \myrhobet \l x^L.P'$, where
    $P \myrhobet M_1$, $Px^L \myrhobet P'$ and $\forall L \in {{\cal
        L}_{\mathbb N}}, x^L \not \in \fv{P}$ and $L \succeq \deg(P)$. By
    lemma~\ref{rhorem} there are two cases:
    \begin{itemize}
    \item $P' = P''x^L$ and $P  \myrhobet P''$. By IH, there is $M'$
      such that  $P'' \myrhobet M' \leftrhobet M_1$. By
      lemma~\ref{newrho2}.\ref{newrho2two}, $\forall L
      \in {{\cal L}_{\mathbb N}}, x^L \not \in \fv{P''}$ and $L\succeq
      \deg(P'')$, hence, $\l x^L. P' = \l x^L.P''x^L  \myrhobet M'
      \leftrhobet M_1$.
    \item $P = \l y^L. Q$, $Q  \myrhobet Q'$, $Q \diamond x^L$ and $P'
      = Q'[y^L := x^L]$.
      So we have $M_1 \leftrhobet \l x^L. (\l y^L. Q)x^L \myrhobet \l
      x^L. Q'[y^L := x^L]$ where $M_1 \leftrhobet \l y^L. Q = \l
      x^L. Q[y^L := x^L]$ since $\forall L
      \in {{\cal L}_{\mathbb N}}, x^L \not \in \fv{P}$.\\ By
      lemma~\ref{newrho1}.\ref{newrho1two}, $\l x^L. Q[y^L := x^L]
      \myrhobet \l x^L. Q'[y^L := x^L]$. Hence by IH, there is $M'$
      such that  $M_1 \myrhobet M' \leftrhobet \l x^L. Q'[y^L := x^L]$.
    \end{itemize}
  \end{itemize}

\item Let $r = \be$:
  \begin{itemize}
  \item
    If $M = x^L$, by lemma~\ref{rhorem}, $M_1 = M_2 = x^L$.  Take $M' =
    x^L$.

  \item If $N_2P_2 \leftrhobe NP \myrhobe N_1P_1$ where $N_2
    \leftrhobe N \myrhobe N_1$, $P_2 \leftrhobe P
    \myrhobe P_1$ and $N \diamond P$, then, by IH, $\exists N', P'$
    such that $N_2
    \myrhobe N'\leftrhobe N_1$ and $P_2 \myrhobe P'\leftrhobe P_1$. By
    lemma~\ref{newrho2}.\ref{newrho2three}, $N_1 \diamond P_1$ and
    $N_2 \diamond P_2$. Hence, $N_2P_2 \myrhobe N'P' \leftrhobe N_1P_1$.

  \item If $(\l x^L.P_1) Q_1 \leftrhobe (\l x^L.P) Q \myrhobe
    P_2[x^L:=Q_2]$ where $\l x^L. P\myrhobe \l x^L.P_1$, $P \myrhobe
    P_2$, $Q_1 \leftrhobe Q \myrhobe Q_2$, $\deg(Q) = L$, $P \diamond
    Q$ and $(\l x^L.P) \diamond Q$, then, by  lemma~\ref{rhorem}, $P
    \myrhobe P_1$. By IH, $\exists
    P', Q'$ such that $P_1 \myrhobe P'\leftrhobe P_2$ and $Q_1
    \myrhobe Q'\leftrhobe Q_2$. By
    lemma~\ref{newrho2}.\ref{newrho2two}, $\deg(Q_1) = \deg(Q_2) =
    \deg(Q) = L$. By lemma~\ref{newrho2}.\ref{newrho2three}, $P_1
    \diamond Q_1$. Hence, $(\l x^L.P_1)Q_1  \myrhobe P'[x^L:=Q']$.

    Moreover, since $P_2 \myrhobe P'$, $Q_2 \myrhobe Q'$,
    $\deg(Q_2) = L$ and by lemma~\ref{newrho2}.\ref{newrho2three}, $P_2
    \diamond Q_2$., then, by lemma~\ref{newrho1}.\ref{newrho1two},
    $P_2[x^L:=Q_2] \myrhobe P'[x^L:=Q']$.

  \item If  $P_1[x^L:= Q_1] \leftrhobe (\l x^L. P) Q \myrhobe
    P_2[x^L:=Q_2]$ where
    $P_1 \leftrhobe P \myrhobe P_2$, $Q_1 \leftrhobe Q \myrhobe Q_2$,
    $\deg(Q) = L$ and $P \diamond Q$ then by IH, $\exists P', Q'$ where $P_1
    \myrhobe P'\leftrhobe P_2$ and $Q_1 \myrhobe Q'\leftrhobe
    Q_2$. By lemma~\ref{newrho2}.\ref{newrho2two}, $\deg(Q_1) =
    \deg(Q_2) = \deg(Q) = L$. By
    lemma~\ref{newrho2}.\ref{newrho2three}, $P_1 \diamond Q_1$ and
    $P_2 \diamond Q_2$. Hence, by
    lemma~\ref{newrho1}.\ref{newrho1two}, $P_1[x^L:=Q_1]\myrhobe
    P'[x^L:=Q'] \leftrhobe P_2[x^L:=Q_2]$.

  \item
    If $\l x^L. N_2 \leftrhobe \l x^L. N \myrhobe\l x^L. N_1$ where $N_2
    \leftrhobe N \myrhobe N_1$, by IH, there is $N'$ such that
    $N_2 \myrhobe N'  \leftrhobe N_1$. Hence, $\l x^L. N_2
    \myrhobe \l x^L. N' \leftrhobe \l x^L. N_1$.
  \end{itemize}
\end{itemize}
\ref{newrho3two}.\ First show by induction on $M \myrrhor M_1$
(and using \ref{newrho3one}) that if $M_2 \leftrhor M \myrrhor
M_1$, then there is $M'$ such that $M_2 \myrrhor M'  \leftrhor
M_1$. Then use this to show \ref{newrho3two} by induction on  $M
\myrrhor M_2$.
\hfill $\Box$
\end{proof}

\begin{proof}[Of Theorem~\ref{confluenceofbetaeta}]
\begin{itemize}
\item[\ref{confitem}.]
For $r \in \{\be, \be\eta\}$, by lemma~\ref{newrho3}.\ref{newrho3two},
$\myrrhor$ is confluent. by lemma~\ref{newrho2}.\ref{newrho2one}
and~\ref{newrho2}.\ref{newrho2two}, $M \myrrhor N$ iff $M \rhd_r^* N$. Then
$\rhd_r^*$ is confluent.\\
For $r = \wh$, since if $M \rhd_r^* M_1$ and $M \rhd_r^* M_2$, $M_1 =
M_2$, we take $M' = M_1$.
\item[\ref{symprop}.]
If) is by definition of $\simeq_r$.  Only if) is by induction on $M_1
\simeq_r M_2$ using \ref{confitem}.
\hfill $\Box$
\end{itemize}
\end{proof}

\section{Proofs of section~\ref{sectypes}}

\begin{proof}[Of lemma~\ref{degnew-goodeg}]
\begin{enumerate}
\item By definition.
\item By induction on $U$.
\begin{itemize}
\item
If $U = a$ ($\deg(U) = \oslash$), nothing to prove.
\item
If $U = V \f T$ ($\deg(U) = \oslash$), nothing to prove.
\item
If $U = \o^L$, nothing to prove.
\item
If $U = U_1 \sqcap U_2$ ($\deg(U) = \deg(U_1) = \deg(U_2) = L$), by IH
we have four cases:
\begin{itemize}
\item
If $U_1 = U_2 = \o^L$ then $U = \o^L$.
\item
If $U_1 = \o^L$ and $U_2 = \vec{e}_L \sqcap_{i=1}^k T_i$ where $k
\geq 1$ and $\forall 1 \leq i \leq k$, $T_i \in {\mathbb T}$ then $U =
U_2$ (since $\o^L$ is a neutral).
\item
If $U_2 = \o^L$ and $U_1 = \vec{e}_L \sqcap_{i=1}^k T_i$ where $k
\geq 1$ and $\forall 1 \leq i \leq k$, $T_i \in {\mathbb T}$ then $U =
U_1$ (since $\o^L$ is a neutral).
\item
If $U_1 = \vec{e}_L \sqcap_{i=1}^p T_i$ and $U_2 =
\vec{e}_L \sqcap_{i=p+1}^{p+q} T_i$ where $p,q \geq 1$,
$\forall 1 \leq i \leq p+q$, $T_i \in {\mathbb T}$ then $U =
\vec{e}_L \sqcap_{i=1}^{p+q} T_i$.
\end{itemize}
\item
If $U = \overline{e}_{n_1}V$ ($L = \deg(U) = n_1::\deg(V) = n_1::K$), by IH we have two
cases:
\begin{itemize}
\item
If $V = \o^K$, $U = \overline{e}_{n_1}\o^K = \o^L$.
\item
If $V = \vec{e}_K \sqcap_{i=1}^p T_i$ where $p \geq 1$ and
$\forall 1 \leq i \leq p$, $T_i \in {\mathbb T}$ then $U =
\vec{e}_L \sqcap_{i=1}^p T_i$ where $p \geq 1$ and $\forall 1
\leq i \leq p$, $T_i \in {\mathbb T}$.
\end{itemize}
\end{itemize}
\item
  \begin{enumerate}
  \item By induction on $U_1 \sqsubseteq U_2$.
  \item By induction on $U_1 \sqsubseteq U_2$.
  \item By induction on $K$. We do the induction step. Let $U_1 =
    \overline{e}_iU$. By induction on $\overline{e}_iU \sqsubseteq U_2$ we obtain $U_2 =
    \overline{e}_iU'$ and $U \sqsubseteq U'$.
  \item same proof as in the previous item.
  \item By induction on $U_1 \sqsubseteq U_2$:
    \begin{itemize}
    \item By $ref$, $U_1 = U_2$.
    \item If $\F{\sqcap_{i=1}^p \vec{e}_K (U_i \f T_i) \sqsubseteq U
    \;\;\; \; \;\;\; U \sqsubseteq U_2}{\sqcap_{i=1}^p
    \vec{e}_K (U_i \f T_i) \sqsubseteq U_2}$. If $U = \o^K$ then by
    (b), $U_2 = \o^K$. If $U = \sqcap_{j=1}^q \vec{e}_K (U'_j \f
    T'_j)$ where $q \geq 1$ and $\forall 1 \leq j \leq q$, $\exists 1
    \leq i \leq p$ such that $U'_j \sqsubseteq U_i$ and $T_i
    \sqsubseteq T'_j$ then by IH, $U_2 = \o^K$ or $U_2 =
    \sqcap_{k=1}^r \vec{e}_K (U''_k \f T''_k)$ where $r \geq 1$ and
    $\forall 1 \leq k \leq r$, $\exists 1 \leq j \leq q$ such that
    $U''_k \sqsubseteq U'_j$ and $T'_j \sqsubseteq T''_k$. Hence, by
    $tr$, $\forall 1 \leq k \leq r$, $\exists 1 \leq i \leq p$ such
    that $U''_k \sqsubseteq U_i$ and $T_i \sqsubseteq T''_k$.
    \item By $\sqcap_E$, $U_2 = \o^K$ or $U_2 = \sqcap_{j=1}^q \vec{e}_K (U'_j \f
    T'_j)$ where $1 \leq q \leq p$ and $\forall 1 \leq j \leq q$,
    $\exists 1 \leq i \leq p$ such that $U_i = U'_j$ and $T_i =
    T'_j$.
  \item Case $\sqcap$ is by IH.
  \item Case $\f$ is trivial.
  \item If $\F{\sqcap_{i=1}^p \vec{e}_L (U_i \f T_i) \sqsubseteq
      U_2}{\sqcap_{i=1}^p \vec{e}_K (U_i \f T_i) \sqsubseteq \overline{e}_iU_2}$
    where $K = i::L$ then by IH, $U_2 = \o^L$ and so $\overline{e}_iU_2 = \o^K$
    or $U_2 = \sqcap_{j=1}^q \vec{e}_L (U'_j \f T'_j)$ so $\overline{e}_iU_2 =
      \sqcap_{j=1}^q \vec{e}_K (U'_j \f T'_j)$ where $q
    \geq 1$ and $\forall 1 \leq j \leq q$, $\exists 1 \leq i \leq p$
    such that $U'_j \sqsubseteq U_i$ and $T_i \sqsubseteq T'_j$.
 \end{itemize}
\end{enumerate}
\item
  By $\sqcap_E$ and since $\o^{L}$ is a neutral.
\item
  By induction on  $U \sqsubseteq U_1' \sqcap U_2'$.
  \begin{itemize}
  \item Let $\F{}{ U_1' \sqcap U_2' \sqsubseteq U_1' \sqcap U_2'}$.
    By $ref$, $U_1' \sqsubseteq U_1'$ and $U_2'
    \sqsubseteq U_2'$.
  \item Let $\F{U \sqsubseteq U'' \;\;\; \; \;\;\; U'' \sqsubseteq U_1'
      \sqcap U_2'} {U \sqsubseteq U_1' \sqcap U_2'}$. By IH,
    $U'' = U_1'' \sqcap U_2''$ such that  $U_1'' \sqsubseteq
    U_1'$ and $U_2'' \sqsubseteq U_2'$. Again by IH, $U = U_1
    \sqcap U_2$ such that  $U_1 \sqsubseteq U_1''$ and $U_2
    \sqsubseteq U_2''$. So by $tr$, $U_1 \sqsubseteq U_1'$ and $U_2
    \sqsubseteq U_2'$.
  \item Let $\F{}{(U_1' \sqcap U_2') \sqcap U \sqsubseteq U_1' \sqcap
      U_2'}$.  By  $ref$,   $U_1' \sqsubseteq U_1'$  and $U_2'
      \sqsubseteq U_2'$.  Moreover $\deg(U) =\deg(U_1' \sqcap U_2') =
      \deg(U_1')$ then by $\sqcap_E$, $U_1' \sqcap U \sqsubseteq
      U_1'$.
  \item If $\F{U_1 \sqsubseteq U_1' \;\;\; \& \;\;\; U_2 \sqsubseteq
      U_2'}{U_1 \sqcap U_2 \sqsubseteq U_1' \sqcap U_2'}$ there is
    nothing to prove.
  \item $\F{V_2 \sqsubseteq V_1 \;\;\; \& \;\;\; T_1 \sqsubseteq
      T_2}{V_1 \f T_1 \sqsubseteq V_2 \f T_2}$ then $U_1' = U_2' = V_2
      \f T_2$ and $U = U_1 \sqcap U_2$ such that $U_1 = U_2 = V_1 \f
      T_1$ and we are done.
    \item If $\F{U \sqsubseteq U_1' \sqcap U_2'}{eU \sqsubseteq eU_1'
        \sqcap eU_2'}$ then by IH $U = U_1 \sqcap U_2$ such that $U_1
      \sqsubseteq U_1'$ and $U_2 \sqsubseteq U_2'$. So, $eU = eU_1
      \sqcap eU_2$ and by $\sqsubseteq_e$, $eU_1 \sqsubseteq eU_1'$
      and $eU_2 \sqsubseteq eU_2'$.
  \end{itemize}
\item
  By induction on  $\G \sqsubseteq \G_1' \sqcap \G_2'$.
  \begin{itemize}
  \item Let $\F{}{\G_1' \sqcap \G_2' \sqsubseteq \G_1' \sqcap \G_2'}$.  By
    $ref$, $\G_1' \sqsubseteq \G_1'$ and $\G_2'
    \sqsubseteq \G_2'$.
  \item Let $\F{\G \sqsubseteq \G'' \;\;\; \; \;\;\; \G'' \sqsubseteq \G_1'
      \sqcap \G_2'} {\G \sqsubseteq \G_1' \sqcap \G_2'}$. By IH,
    $\G'' = \G_1'' \sqcap \G_2''$ such that  $\G_1'' \sqsubseteq
    \G_1'$ and $\G_2'' \sqsubseteq \G_2'$. Again by IH, $\G = \G_1
    \sqcap \G_2$ such that  $\G_1 \sqsubseteq \G_1''$ and $\G_2
    \sqsubseteq \G_2''$. So by $tr$, $\G_1 \sqsubseteq \G_1'$ and $\G_2
    \sqsubseteq \G_2'$.
  \item Let $\F{U_1 \sqsubseteq U_2}{\G, (y^n : U_1) \sqsubseteq  \G, (y^n :
      U_2)}$ where $\G, (y^n : U_2) = \G_1' \sqcap
    \G_2'$. 
    \begin{itemize}
    \item If $\G_1' = \G_1'', (y^n : U_2')$ and $\G_2' = \G_2'', (y^n
      : U_2'')$ such that $U_2 = U_2' \sqcap U_2''$, then
      by~\ref{goodegthree}, $U_1 = U_1' \sqcap U_1''$ such
      that $U_1' \sqsubseteq U_2'$ and $U_1'' \sqsubseteq U_2''$. Hence $\G
      = \G_1'' \sqcap \G_2''$ and $\G, (y^n : U_1) = \G_1 \sqcap \G_2$
      where $\G_1 = \G_1'', (y^n : U'_1)$ and $\G_2 = \G_2'', (y^n :
      U''_1)$ such that $\G_1 \sqsubseteq \G_1'$ and $\G_2 \sqsubseteq
      \G_2'$ by $\sqsubseteq_c$.
    \item If $y^n  \not \in \dom{\G_1'}$ then $\G = \G_1'\sqcap \G_2''$
      where $\G_2'', (y^n:U_2) = \G_2'$. Hence, $\G, (y^n:U_1) =
      \G_1'\sqcap\G_2$ where $\G_2 = \G_2'',(y^n:U_1)$. By $ref$ and
      $\sqsubseteq_c$, $\G_1' \sqsubseteq \G_1'$ and $\G_2 \sqsubseteq
      \G_2'$.
    \item If  $y^n  \not \in \dom{\G_2'}$ then similar to the above case.
\hfill $\Box$
    \end{itemize}
  \end{itemize}
\end{enumerate}
\end{proof}

\begin{proof}[Of lemma~\ref{env-Phisub}]
\ref{lem:legalomegaenv}.\ By definition, if $\fv{M} = \{x^{L_1}_1,
\dots, x^{L_n}_n\}$ then $env^{\o}_{M} = (x^{L_i}_i : \o^{L_i})_n$
and by definition, for all $i \in \{1, \dots, n\}$, $\deg(\o^{L_i}) =
L_i$. Moreover, if $x^{L} : U, x^{L} : V \in env^{\o}_{M}$, then $U =
\o^{L} = V$.\\
\ref{Phisubone}.\ First show by induction on the derivation $\G
\sqsubseteq \G'$ that if  $\G \sqsubseteq \G'$ and $\G, (x^L:U)$ is an
environment, then $\G, (x^L:U)\sqsubseteq \G',  (x^L:U)$. Then use
$(tr)$ and $(\sqsubseteq_c)$.\\
\ref{Phisubtwo}.\ Only if) By induction on the derivation $\G
\sqsubseteq \G'$. If) By induction on $n$ using  \ref{Phisubone}.\\
\ref{Phisubthree}.\ Only if) By induction on the derivation $\<\G \v
U\> \sqsubseteq \<\G' \v U'\>$. If) By $\sqsubseteq_{\<\>}$.\\
\ref{Phisubfour}.\
Let $\fv{M} =\{x^{L_1}_1, \dots, x^{L_n}_n\}$ and  $\G = (x^{L_i}_i
:U_i)_n$. By definition, $env^\o_M = (x^{L_i}_i : \o^{L_i})_n$.
Because $\legalenv{\G}$, then for all $i \in \{1, \dots, n\}$,
$\deg(U_i) = L_i$.  Hence,
by lemma~\ref{degnew-goodeg}.\ref{goodegtwo} and \ref{Phisubtwo},  $\G
\sqsubseteq env^\o_M$.\\
\ref{Phisubfive}.\
Let $x^{L_1} \in \dom{\G^{-K}}$ and $x^{L_2} \in \dom{\D^{-K}}$, then
$x^{K::L_1} \in \dom{\G}$ and $x^{K::L_2} \in \dom{\D}$, hence $K::L_1 =
K::L_2$ and so $L_1 = L_2$.\\
\ref{Phisubnew}.\ Let $\deg(U) = L = K::K'$. By lemma~\ref{degnew-goodeg}:
  \begin{itemize}
    \item If $U = \o^L$ then by
    lemma~\ref{degnew-goodeg}.\ref{goodegonetwo}, $U' = \o^L$
    and by $ref$, $U^{-K} = \o^{K'} \sqsubseteq \o^{K'} = U'^{-K}$.
    \item If $U = \vec{e}_L \sqcap_{i=1}^p T_i$ where $p \geq 1$ and
    $\forall~1 \leq i \leq p$, $T_i \in {\mathbb T}$ then by
    lemma~\ref{degnew-goodeg}.\ref{goodegonethree}, $U' =
    \vec{e}_L V$ and $\sqcap_{i=1}^p T_i \sqsubseteq V$. Hence, by
    $\sqsubseteq_e$,  $U^{-K} = \vec{e}_{K'} \sqcap_{i=1}^p T_i
    \sqsubseteq \vec{e}_{K'} V = U'^{-K}$.
  \end{itemize}
\ref{Phisubnewnew}.\ Let $\G = (x^{L_i}_i : U_i)_n$, so
  by lemma~\ref{env-Phisub}.\ref{Phisubtwo}, $\G' = (x^{L_i}_i :
  U'_i)_n$ and $\forall 1 \leq i \leq n$, $U_i \sqsubseteq
  U'_i$. Because $\deg(\G) \succeq K$, then by definition $\forall 1
  \leq i \leq n$,
  $\deg(U_i)  \succeq K$. By lemma~\ref{env-Phisub}.\ref{Phisubnew},
  $\forall i \in \{1, \dots,  n\}$, $U_i^{-K} \sqsubseteq
  {U'}_i^{-K}$ and by lemma~\ref{env-Phisub}.\ref{Phisubtwo}, $\G^{-K}
  \sqsubseteq \G'^{-K}$.\\
\ref{lem:interlegalenv}.\ Let $\G_1 = (x^{L_i}_i : U_i)_n, \G'_1$ and
$\G_2 = (x^{L_i}_i : U'_i)_n, \G'_2$ such that $\dom{\G'_1} \cap
\dom{\G'_2}$. Then, by hypotheses, for all $i \in \{1, \dots, n\}$,
$\deg(U_i) = L_i = \deg(U'_i)$.  Then $\G_1 \sqcap \G_2 = (x^{L_i}_i :
U_i \sqcap U'_i)_n, \G'_1, \G'_2$ is well defined. Moreover, for all
$x^{L} : U \in \G'_1$, $\deg(U) = L$ and for all $x^{L} : U \in
\G'_2$, $\deg(U) = L$ and for all $i \in \{1, \dots, n\}$, $\deg(U_i
\sqcap U'_i) = \deg(U_i) = L_i = \deg(U'_i)$.\\
\ref{lem:explegalenv}.\ Let $\G = (x^{L_j}_j : U_j)_n$ then by
hypothesis, for all $j \in \{1, \dots, n\}$, $\deg(U_j) = L_j$ and
$\overline{e}_i\G = (x^{i::L_j}_j : \overline{e}_iU_j)$. So, for all $j \in \{1, \dots,
n\}$, $\deg(\overline{e}_iU_j) = i::\deg(U_j) = i::L_j$.\\
\ref{lem:deg+legalenvsub}.\  By
lemma~\ref{env-Phisub}.\ref{Phisubtwo}, $\G_1 = (x^{L_i}_i : U_i)_n$
and $\G_2 = (x^{L_i}_i : U'_i)_n$ and for all $i \in \{1, \dots, n\}$,
$U_i \sqsubseteq U'_i$. By
lemma~\ref{degnew-goodeg}.\ref{goodegoneone}, for all $i \in \{1,
\dots, n\}$, $\deg(U_i) = \deg(U'_i)$.
Assume $\deg(\G_1) \succeq K$ then for all $i \in \{1, \dots, n\}$,
$\deg(U_i) = \deg(U'_i) \succeq K$ and $L_i \succeq K$, so $\deg(\G_2)
\succeq K$.
Assume $\deg(\G_2) \succeq K$ then for all $i \in \{1, \dots, n\}$,
$\deg(U_i) = \deg(U'_i) \succeq K$ and $L_i \succeq K$, so $\deg(\G_1)
\succeq K$.
Assume $\legalenv{\G_1}$ then for all $i \in \{1, \dots, n\}$, $L_i =
\deg(U_i) = \deg(U'_i)$, so $\legalenv{\G_2}$.
Assume $\legalenv{\G_2}$ then for all $i \in \{1, \dots, n\}$, $L_i =
\deg(U'_i) = \deg(U_i)$, so $\legalenv{\G_1}$.
\hfill $\Box$
\end{proof}

\begin{proof}[Of theorem~\ref{typing'}]
  \begin{enumerate}
  \item[\ref{typing'zero}.\ and \ref{typing'one}.] By
    lemma~\ref{structyping}.\ref{structypingone} and
    lemma~\ref{structypingC}.\ref{structypingC}, $\G \diamond \G$.
    \begin{itemize}
    \item If $\F{}{x^{\oslash} : \<(x^{\oslash}:T) \v T\>}$ then, by
      hypothesis, $T \in {\mathbb T} \subseteq {\mathbb U}$ and
      $\deg(T) = \oslash = \deg(x^{\oslash})$. So,
      $\legalenv{(x^{\oslash}:T)}$ and $x^{\oslash} \in {\cal M}$.

    \item If $\F{}{M : \<env^\o_M \v \o^{\deg(M)}\>}$. By definition,
      $M$ is defined to range over ${\cal M}$ and
      $\legalenv{env^\o_M}$ by
      lemma~\ref{env-Phisub}.\ref{lem:legalomegaenv}.  By definition,
      $\o^{\deg(M)} \in {\mathbb U}$.  Let $\fv{M} =
      \{x^{L_1},\dots,x^{L_n}\}$, so $env^\o_M = (x^{L_i}_i :
      \o^{L_i})_n$ and by lemma~\ref{degsub}.\ref{degsubtwo}, $\forall
      1 \leq i \leq n, L_i \succeq \deg(M) = \deg(\o^{\deg(M)})$.

    \item If $\F{M : \<\G,(x^L:U) \v T\>}{\l x^L. M : \<\G \v U \f
        T\>}$ then by IH, $M \in {\cal M}$, $T \in {\mathbb U}$
      $\G,(x^L:U) \in Env$, $\legalenv{\G,(x^L:U)}$ and
      $\deg(\G,(x^L:U)) \succeq \deg(T) = \deg(M)$.  By hypothesis, $T
      \in {\mathbb T}$. Because $\G,(x^L:U) \in Env$, then $U \in
      {\mathbb U}$. So $U \f T \in {\mathbb T} \subset {\mathbb U}$.
      Let $\G = (x^{L_i}_i : U_i)_n$, then for all $i \in \{1, \dots,
      n\}$, $L_i = \deg(U_i) \succeq \deg(T) = \deg(U \f T)$ and
      $\deg(U) = L \succeq \deg(M)$. Hence, $\l x^{L}. M \in {\cal M}$
      and $\legalenv{\G}$.  So, $\deg(\l x^L. M) = \deg(M) = \deg(T) =
      \deg(U \f T)$.

    \item If $\F{M : \<\G \v T\>\;\;\; x^L \not \in \dom{\G}}{\l
        x^L. M : \<\G \v \o^L \f T\>}$ then by IH, $M \in {\cal M}$,
      $T \in {\mathbb U}$, $\G \in Env$, $\legalenv{\G}$ and $\deg(\G)
      \succeq \deg(T) = \deg(M)$.  By hypothesis, $T \in {\mathbb T}$.
      So $\deg(T) = \oslash = \deg(M) \preceq L$.  By definition,
      $\o^{L} \in {\mathbb U}$.  So, $\o^{L} \f T \in {\mathbb T}
      \subset {\mathbb U}$.
      So, $\l x^{L}. M \in {\cal M}$ and
      $\deg(\l x^L. M) = \deg(M) = \deg(T) = \deg(\o^{L} \f T)$.

    \item If $\F{M_1 : \<\G_1 \v U \f T\> \;\;\; \hspace{0.2in} M_2 :
        \<\G_2 \v U\> \;\;\; \hspace{0.2in} \G_1 \diamond \G_2}{M_1
        M_2 : \<\G_1 \sqcap \G_2 \v T\>}$ then by IH, $M_1, M_2 \in
      {\cal M}$, $\G_1, \G_1 \in Env$, $U \f T, U \in {\mathbb U}$,
      $\legalenv{\G_1}$, $\legalenv{\G_2}$ and $\deg(\G_1) \succeq
      \deg(U \f T) = \deg(M_1)$ and $\deg(\G_2) \succeq \deg(U) =
      \deg(M_2)$.  By definition, $\G_1 \sqcap \G_2$ is a type
      environment.  By hypothesis, $T \in {\mathbb T} \subset {\mathbb
        U}$.  By lemma~\ref{env-Phisub}.\ref{lem:interlegalenv} and
      lemma~\ref{structyping}.\ref{structypingC}, $\legalenv{\G_1
        \sqcap \G_2}$ and $M_1 \diamond M_2$.  Because $\deg(M_2) =
      \deg(U) \succeq \oslash = \deg(U \f T) = \deg(M_1)$, then
      $M_1M_2 \in {\cal M}$.  We have trivially, $\deg(\G_1 \sqcap
      \G_1) \succeq \oslash$.  Moreover $\deg(M_1M_2) = \deg(M_1) =
      \deg(U \f T) = \deg(T)$.

    \item If $\F{M: \<\G \v U_1\> \;\;\;\hspace{0.2in} M : \<\G \v
        U_2\>}{M : \<\G \v U_1 \sqcap U_2\>}$ then by IH, $M \in {\cal
        M}$, $\G \in Env$, $U_1, U_2 \in {\mathbb U}$, $\legalenv{\G}$
      and $\deg(\G) \succeq \deg(U_1) = \deg(M)$ and $\deg(\G) \succeq
      \deg(U_2) = \deg(M)$.  So $\deg(U_1) = \deg(M) = \deg(U_2)$.
      Hence, $U_1 \sqcap U_2 \in {\mathbb U}$.  Moreover, $\deg(\G)
      \succeq \deg(U_1) = \deg(U_1 \sqcap U_2) = \deg(M)$.

    \item If $\F{M : \<\G \v U\>}{M^{+k} : \< \overline{e}_k\G \v
        \overline{e}_kU\>}$ then by IH, $M \in {\cal M}$, $\G \in
      Env$, $U \in {\mathbb U}$, $\legalenv{\G}$ and $\deg(\G) \succeq
      \deg(U) = \deg(M)$.  Then, by definition, $\overline{e}_{k}U \in
      {\mathbb U}$.  By definition, $\overline{e}_k\G \in Env$.  Then,
      by lemma~\ref{deg+-f+}.\ref{one'} and
      lemma~\ref{env-Phisub}.\ref{lem:explegalenv}, $M^{+i} \in {\cal
        M}$ and $\legalenv{\overline{e}_{k}\G}$.  Let $\G = (x^{L_j}_j
      : U_j)_n$ so $\overline{e}_k\G= (x^{k::L_j}_j :
      \overline{e}_kU_j)_n$ and for all $j \in \{1, \dots, n\}$,
      because $\deg(U_j) = L_j \succeq \deg(U)$ then
      $\deg(\overline{e}_kU_j) = k::\deg(U_j) = k::L_j \succeq
      k::\deg(U) = \deg(\overline{e}_kU) = k::\deg(M)
      =^{\ref{deg+-f+}.\ref{one'}} \deg(M^{+k})$.

    \item If $\F{M : \<\G\v U\> \;\;\;\hspace{0.2in} \<\G\v U\>
        \sqsubseteq \<\G'\v U'\>}{M : \<\G'\v U'\>}$ then by IH, $M
      \in {\cal M}$, $U \in {\mathbb U}$, $\G \in Env$,
      $\legalenv{\G}$ and $\deg(\G) \succeq \deg(U) = \deg(M)$.  By
      lemma~\ref{env-Phisub}.\ref{Phisubthree}, $\G' \sqsubseteq \G$,
      hence, $\G' \in Env$.  By
      lemma~\ref{env-Phisub}.\ref{lem:deg+legalenvsub},
      $\legalenv{\G'}$.  Let $\G= (x^{L_i}_i : U_i)_n$, so $\forall 1
      \leq i \leq n, \deg(U_i) = L_i \succeq \deg(U)$. By
      lemma~\ref{env-Phisub}.\ref{Phisubtwo}, $\G' = (x^{L_i}_i :
      U'_i)_n$ and $\forall 1 \leq i \leq n$, $U_i \sqsubseteq U'_i$
      so by lemma~\ref{degnew-goodeg}.\ref{goodegoneone}, $\deg(U_i) =
      \deg(U'_i)$. By lemma~\ref{env-Phisub}.\ref{Phisubthree}, $U
      \sqsubseteq U'$ so by
      lemma~\ref{degnew-goodeg}.\ref{goodegoneone}, $\deg(U) =
      \deg(U')$. Hence $\forall 1 \leq i \leq n, \deg(U'_i) = L_i
      \succeq \deg(U') = \deg(M)$.
    \end{itemize}
  \item[\ref{typing'two'}.]
    By induction on $M : \<\G \v U\>$. Case $K = \oslash$ is
    trivial, consider $K = i::K'$. Let $\deg(U) = K::L$. Since
    $\deg(U) \succeq K$, $U^{-K}$ is well defined. Since by
    1. $\deg(\G) \succeq \deg(U) = \deg(M)$, $M^{-K}$ and
    $\G^{-K}$ are well defined too.
    \begin{itemize}
    \item
      If $\F{}{M : \<env^\o_M \v \o^{\deg(M)}\>}$. 
      By $\o$, $M^{-K} : \<env^\o_{M^{-K}} \v \o^L\>$.
    \item $\sqcap_I$ is by IH.
    \item If $\F{M : \<\G \v U\>}{M^{+i} : \< \overline{e}_i\G \v \overline{e}_iU\>}$. Since
      $\deg(\overline{e}_iU) = i::K'::L$, $\deg(U) = K'::L$, so $\deg(U) \succeq
      K'$ and by IH, $M^{-K'} : \<\G^{-K'} \v U^{-K'}\>$, so by $e$
      and lemma~\ref{deg+-f+}.\ref{one}, $(M^{+i})^{-K} :
      \<(\overline{e}_i\G)^{-K} \v (\overline{e}_iU)^{-K}\>$.
    \item If $\F{M : \<\G\v U\> \;\;\;\hspace{0.2in} \<\G\v U\>
        \sqsubseteq \<\G'\v U'\>}{M : \<\G'\v U'\>}$ then by
      lemma~\ref{env-Phisub}.\ref{Phisubthree}, $\G' \sqsubseteq \G$
      and $U \sqsubseteq U'$. By
      lemma~\ref{degnew-goodeg}.\ref{goodegoneone}, $\deg(U) =
      \deg(U') \succeq K$. By IH, $M^{-K} : \<\G^{-K} \v
      U^{-K}\>$. Hence by
      lemma~\ref{env-Phisub}.\ref{lem:deg+legalenvsub},
      lemma~\ref{env-Phisub}.\ref{Phisubnew},
      lemma~\ref{env-Phisub}.\ref{Phisubnewnew} and $\sqsubseteq$,
      $M^{-K} : \<\G'^{-K} \v U'^{-K}\>$.
\hfill $\Box$
   \end{itemize}
 \end{enumerate}
\end{proof}

\begin{proof}[Of  remark~\ref{remderiv}]
  \begin{enumerate}
  \item Let $M: \<\G_1 \v U_1\>$ and $M : \<\G_2 \v U_2\>$. By
    lemma~\ref{structyping}.\ref{structypingone}, $\dom{\G_1} = \fv{M}
    = \dom{\G_2}$. Let $\G_1 = (x_i^{L_i} : V_i)_n$ and $\G_2 =
    (x_i^{L_i} : V'_i)_n$. Then, by
    lemma~\ref{typing'}.\ref{typing'one}, $\forall 1 \leq i \leq n$,
    $\deg(V_i) = \deg(V'_i) = L_i$. By $\sqcap_E$, $V_i \sqcap V'_i
    \sqsubseteq V_i$ and $V_i \sqcap V'_i \sqsubseteq V'_i$. Hence, by
    lemma~\ref{env-Phisub}.\ref{Phisubtwo}, $\G_1 \sqcap \G_2
    \sqsubseteq\G_1$ and $\G_1 \sqcap \G_2 \sqsubseteq\G_2$ and by
    $\sqsubseteq$ and $\sqsubseteq_{\<\>}$, $M :\<\G_1 \sqcap \G_2 \v
    U_1\>$ and $M :\<\G_1 \sqcap \G_2 \v U_2\>$.  Finally, by
    $\sqcap_I$, $M :\<\G_1 \sqcap \G_2 \v U_1\sqcap U_2\>$.
  \item
    By lemma \ref{degnew-goodeg},  either $U = \o^L$ so by $\o$, $x^L :
    \<(x^L : \o^L) \v \o^L\>$.Or $U = \sqcap_{i=1}^p \vec{e}_L T_i$
    where $p \geq 1$, and $\forall 1 \leq i \leq p$,  $T_i \in
    {\mathbb T}$. Let $1 \leq i \leq p$.By $ax$, $x^\oslash :
    \<(x^\oslash : T_i) \v T_i\>$, hence by $e$, $x^L : \<(x^L :
    \vec{e}_L T_i) \v \vec{e}_L T_i\>$. Now, by $\sqcap_I'$, $x^L :
    \<(x^L : U) \v U\>$.
\hfill $\Box$
  \end{enumerate}
\end{proof}

\section{Proofs of section~\ref{srsec}}

\begin{proof}[Of lemma~\ref{newgen}]
  \ref{newgenone}.\ By induction on  the derivation $x^L : \<\G \v  U
  \>$. We have fives cases:
  \begin{itemize}
  \item If $\F{}{x^{\oslash} : \<(x^{\oslash}:T) \v T\>}$ then it is
    done using (ref).
  \item If $\F{}{x^L : \<(x^L :\o^L) \v \o^L\>}$ then it is done using
    (ref).
  \item If $\F{x^L: \<\G \v U_1\> \;\;\;\hspace{0.2in} x^L : \<\G \v U_2\>}
      {x^L : \<\G \v U_1 \sqcap U_2\>}$. By IH, $\G = (x^L : V)$, $V
    \sqsubseteq U_1$ and $V \sqsubseteq U_2$, then by rule $\sqcap$, $ V
    \sqsubseteq U_1 \sqcap U_2$.
  \item If $\F{x^L : \<\G \v U\>}{x^{i::L} : \< \overline{e}_i\G \v \overline{e}_iU\>}$. Then
    by IH, $\G = (x^L : V)$ and $V \sqsubseteq U$, so  $\overline{e}_i\G =
    (x^{i::L} : \overline{e}_iV)$ and by $\sqsubseteq_e$, $\overline{e}_iV \sqsubseteq \overline{e}_iU$,
  \item If $\F{x^L:\<\G' \v U'\>\;\;\;\<\G' \v U'\> \sqsubseteq \<\G
      \v U\> }{x^L:\< \G \v U\>}$. By
    lemma~\ref{env-Phisub}.\ref{Phisubthree}, $\G \sqsubseteq \G'$ and $U'
    \sqsubseteq U$ and, by IH, $\G' = (x^L : V')$ and $V' \sqsubseteq
    U'$. Then, by lemma~\ref{env-Phisub}.\ref{Phisubtwo}, $\G = (x^L : V)$, $V
    \sqsubseteq V'$ and, by rule $tr$, $V \sqsubseteq U$.
  \end{itemize}
  \ref{newgentwo}.\ By induction on  the derivation $\l x^L. M : \<\G
  \v  U\>$. We have five cases:
  \begin{itemize}
  \item If $\F{}{\l x^L. M : \<env^\o_{\l x^L. M} \v  \o^{\deg(\l
        x^L. M)}\>}$ then it is done.
  \item If $\F{M : \<\G, x^L : U \v  T\>}{\l x^L. M : \<\G \v U \f T\>}$
    ($\deg(U \f  T) = \oslash$) then it is done.
  \item If $\F{\l x^L. M: \<\G \v  U_1\> \; \; \l x^L. M : \<\G \v
      U_2\>} {\l x^L. M : \<\G \v  U_1 \sqcap U_2\>}$ then $\deg(U_1 \sqcap
    U_2) = \deg(U_1) = \deg(U_2) = K$. By IH, we have four cases:
    \begin{itemize}
    \item If $U_1 = U_2 = \o^K$, then $U_1 \sqcap U_2 = \o^K$.
    \item If $U_1 = \o^K$, $U_2 = \sqcap_{i=1}^p \vec{e}_K (V_i \f T_i)$
      where $p \geq 1$ and $\forall 1 \leq i \leq p$, $M  : \<\G,x^L :
      \vec{e}_K V_i \v \vec{e}_K T_i\>$, then $U_1 \sqcap U_2 = U_2$
      ($\o^K$ is a neutral element).
    \item If $U_2 = \o^K$, $U_1 = \sqcap_{i=1}^p \vec{e}_K (V_i \f T_i)$
      where $p \geq 1$ and $\forall 1 \leq i \leq p$, $M  : \<\G ,x^L :
      \vec{e}_K V_i \v \vec{e}_K T_i\>$, then $U_1 \sqcap U_2 = U_1$
      ($\o^K$ is a neutral element).
    \item If $U_1 = \sqcap_{i= 1}^p \vec{e}_K (V_i \f T_i)$, $U_2 =
      \sqcap_{i= p+1}^{p+q} \vec{e}_K (V_{i} \f T_{i})$ (hence
      $U_1\sqcap U_2 = \sqcap_{i= 1}^{p+q} \vec{e}_K (V_i \f T_i)$)
      where $p,q \geq 1$, $\forall 1 \leq i \leq p+q$, $M : \<\G , x^L :
      \vec{e}_K V_i \v \vec{e}_K T_i\>$, we are done.
    \end{itemize}
  \item If $\F{\l x^L. M : \<\G \v U\>}{\l x^{i::L}. M^{+i} : \< \overline{e}_i\G
      \v \overline{e}_iU\>}$. $\deg(\overline{e}_iU) = i::\deg(U) = i::K' = K$. By IH, we
    have two cases:
    \begin{itemize}
    \item If $U = \o^{K'}$ then $\overline{e}_iU = \o^K$.
    \item If $U = \sqcap_{j=1}^p \vec{e}_{K'} (V_j \f T_j)$, where
      $p\geq 1$ and for all $1 \leq j \leq p$, $M : \<\G , x^L :
      \vec{e}_{K'} V_j \v  \vec{e}_{K'} T_j\>$. So $\overline{e}_iU =
      \sqcap_{j=1}^p \vec{e}_K (V_j \f T_j)$ and by $e$, for all $1
      \leq j \leq p$, $M^{+i} : \<\overline{e}_i\G , x^{i::L} : \vec{e}_K V_j \v
      \vec{e}_K T_j\>$.
    \end{itemize}
  \item Let $\F{\l x^L. M : \<\G\v U\> \;\;\;  \<\G\v U\> \sqsubseteq \<\G'\v
      U'\>}{\l x^L. M : \<\G'\v U'\>}$.  By
      lemma~\ref{env-Phisub}.\ref{Phisubthree}, $\G' \sqsubseteq \G$ and
      $U \sqsubseteq U'$ and by
      lemma~\ref{degnew-goodeg}.\ref{goodegoneone} $\deg(U) =
      \deg(U') = K$. By IH, we have two cases:
    \begin{itemize}
    \item If $U = \o^K$, then, by
      lemma~\ref{degnew-goodeg}.\ref{goodegonetwo}, $U' = \o^K$.
    \item If $U = \sqcap_{i=1}^p \vec{e}_K (V_i \f T_i)$, where $p\geq
      1$ and for all $1 \leq i \leq p$ $M : \<\G , x^L : \vec{e}_K V_i
      \v  \vec{e}_K T_i\>$. By lemma~\ref{degnew-goodeg}.\ref{goodegonefour}:
      \begin{itemize}
      \item Either $U' = \o^K$.
      \item Or $U' = \sqcap_{i=1}^q
        \vec{e}_K (V'_i \f T'_i)$, where $q\geq 1$ and $\forall 1 \leq
        i \leq q$, $\exists 1 \leq j_i \leq p$ such that $V'_i
        \sqsubseteq V_{j_i}$
        and $T_{j_i} \sqsubseteq T'_i$. Let $1 \leq i \leq q$. Since, by
        lemma~\ref{env-Phisub}.\ref{Phisubthree}, $\<\G ,x^L : \vec{e}_K V_{j_i} \v
        \vec{e}_K T_{j_i} \>\sqsubseteq \<\G' , x^L : \vec{e}_K V'_i \v
        \vec{e}_K T'_i \>$, then $M : \<\G', x^L : \vec{e}_K V'_i \v
        \vec{e}_K T'_i\>$.
      \end{itemize}
    \end{itemize}
  \end{itemize}
  \ref{newgenthree}.\ Similar as the proof of \ref{newgentwo}.\\
  \ref{newgenfour'}.\ By induction on the derivation $M \, x^L : \<\G
  ,x^L : U \v T\>$.
  We have two cases:
  \begin{itemize}
  \item
    Let $\F{M: \<\G \v  V \f T\> \;\;\; x^L : \<(x^L : U)  \v
      V\>\;\;\; \G \diamond (x^L : U)}{M
    \, x^L : \<\G , (x^L : U) \v T\>}$ (where, by 1.\, $U \sqsubseteq
    V$). Since $V \f T \sqsubseteq U \f T$, we have $M: \<\G \v U \f
    T\>$.
  \item
    Let $\F{M \; x^L : \<\G', (x^L:U') \v V'\> \;\;\; \<\G', (x^L:U')
    \v V'\> \sqsubseteq \<\G, (x^L:U) \v V\>}{M \; x^L : \<\G, (x^L
    :U) \v V\>}$ (by lemma~\ref{env-Phisub}).
    By lemma~\ref{env-Phisub}, $\G \sqsubseteq \G'$, $U\sqsubseteq
    U'$ and $V' \sqsubseteq V$. By IH, $M : \<\G' \v U' \f V'\>$ and
    by $\sqsubseteq$, $M : \<\G \v U \f V\>$.
  \end{itemize}
\end{proof}

\begin{proof}[Of lemma~\ref{substlem}]
  By lemma~\ref{typing'}.\ref{typing'one}, $M, N \in {\cal M}$,
  $\deg(N) = \deg(U)$, $\legalenv{\Delta}$ and $\legalenv{\G,x^L:U}$,
  so $\deg(N) = \deg(U) = L$. By
  lemma~\ref{env-Phisub}.\ref{lem:interlegalenv}, $\legalenv{\G \sqcap
    \Delta}$.  By lemma~\ref{degsub}.\ref{degsubthree}, $M[x^{L}:=N]
  \in {\cal M}$. By lemma~\ref{structyping}.\ref{structypingone},
  $x^{L} \in \fv{M}$.

  We prove the lemma by induction on the derivation $M:\<\G,x^L:U \v
  V\>$.
  \begin{itemize}
  \item If $\F{}{x^{\oslash} : \<(x^{\oslash}:T) \v T\>}$ and $N
    :\<\Delta\v T\>$, then $x^{\oslash}[x^{\oslash}:=N] = N : \<\Delta
    \v T\>$.
  \item If $\F{}{M : \<env^\o_{\fv{M}\setminus\{x^L\}},(x^L : \o^L) \v
      \o^{\deg(M)}\>}$ and $N :\<\Delta\v \o^L\>$ then by $\o$,
    $M[x^L:=N]:\<env^\o_{M[x^L:=N]} \v \o^{\deg(M[x^L:=N])}\>$. By
    lemma~\ref{degsub}.\ref{degsubthree} $\deg(M[x^L:=N]) =
    \deg(M)$. Since $x^L \in \fv{M}$ (and so $\fv{M[x^L:=N]} = (\fv{M}
    \setminus \{x^{L}\}) \cup \fv{N}$), by $\sqsubseteq$,
    $M[x^L:=N]:\<env^\o_{\fv{M}\setminus\{x^L\}} \sqcap \Delta \v
    \o^{\deg(M)}\>$.
  \item Let $\F{M : \<\G, x^L:U, y^K:U' \v T\>}{\l y^K. M : \<\G,
      x^L:U \v U' \f T\>}$ where $y^K \not \in \fv{N} \cup
    \{x^{L}\}$. So $(\l y^{K}. M)[x^{L}:=N] = \l y^{K}. M[x^{L}:=N]$.
    By lemma~\ref{degsub}.\ref{degsubone}, $M \diamond N$.  By IH,
    $M[x^L:=N] : \<\G \sqcap \Delta, y^K:U' \v T\>$. By $\f_I$, $(\l
    y^K.M)[x^L:=N] : \<\G\sqcap \Delta \v U'\f T\>$.
  \item Let $\F{M : \<\G, x^L:U \v T\>\;\;\; y^K \not \in \dom{\G,
        x^L:U}}{\l y^K. M : \<\G, x^L:U \v \o^K \f T\>}$ where $y^K
    \not \in \fv{N} \cup \{x^{L}\}$.  So $(\l y^{K}. M)[x^{L}:=N] = \l
    y^{K}. M[x^{L}:=N]$.  By lemma~\ref{degsub}.\ref{degsubone}, $M
    \diamond N$.  By lemma~\ref{structyping}.\ref{structypingone},
    $\fv{N} = \dom{\Delta}$, so $y^{K} \not \in \dom{\Delta}$.  By IH,
    $M[x^L:=N] : \<\G \sqcap \Delta \v T\>$. By $\f'_I$, $(\l
    y^K.M)[x^L:=N] : \<\G\sqcap \Delta \v \o^K \f T\>$.
  \item Let $\F{M_1 : \<\G_1, x^L:U_1 \v V \f T\> \;\;\; M_2 : \<\G_2,
      x^L:U_2 \v V\> \;\;\; \G_1 \diamond \G_2}{M_1 M_2 : \<\G_1
      \sqcap \G_2, x^L:U_1\sqcap U_2 \v T\>}$ (by
    lemma~\ref{structyping}.\ref{structypingone}) where $x^L \in
    \fv{M_1}\cap \fv{M_2}$, $N :\<\Delta \v U_1\sqcap U_2\>$. By
    lemma~\ref{degsub}.\ref{degsubone}, $M_1 \diamond N$ and $M_2
    \diamond N$. By $\sqcap_E$ and $\sqsubseteq$, $N :\<\Delta\v
    U_1\>$ and $N :\<\Delta\v U_2\>$. Now use IH and $\f_E$ (using the
    fact that $\G_1 \sqcap \Delta \diamond \G_2 \sqcap \Delta$, by
    lemma~\ref{structyping}.\ref{structypingone} and
    lemma~\ref{degsub}.\ref{degsubfour}).

    The cases $x^L \in \fv{M_1}\setminus \fv{M_2}$ or $x^L \in
    \fv{M_2}\setminus \fv{M_1}$ are similar.
  \item If $\F{M: \<\G, x^L:U \v U_1\> \;\; M : \<\G, x^L:U \v U_2\>}
    {M : \<\G, x^L:U \v U_1 \sqcap U_2\>}$ use IH and $\sqcap_I$.
  \item Let $\F{M : \<\G, x^L:U \v V\>}{M^{+i} : \< \overline{e}_i\G ,
      x^{i::L}:\overline{e}_iU \v \overline{e}_iV\>}$ and $N :\<\Delta\v \overline{e}_iU\>$.  By
    lemma~\ref{typing'}.\ref{typing'one}, $\deg(M) = \deg(\overline{e}_iU) =
    i::\deg(U)$.  By lemma~\ref{typing'}.\ref{typing'two'},
    $N^{-i}:\<\Delta^{-i} \v U\>$.  By lemma~\ref{deg+-f+}.\ref{four}
    and lemma~\ref{deg+-f+}.\ref{one''}, $(N^{-i})^{+i} = N$ and $M
    \diamond N^{-i}$.  By IH, $M[x^L:=N^{-i}] : \<\G\sqcap \Delta^{-i}
    \v V\>$. By $e$ and lemma~\ref{deg+-f+}.\ref{two},
    $M^{+i}[x^{i::L}:=N] : \<\overline{e}_i\G\sqcap \Delta \v \overline{e}_iV\>$.
  \item Let $\F{M : \<\G', x^L:U' \v V'\> \;\;\; \<\G', x^L:U' \v V'\>
      \sqsubseteq \<\G, x^L:U \v V\>}{M : \<\G, x^L:U \v V\>}$
    (lemma~\ref{env-Phisub}). By lemma~\ref{env-Phisub}, $\dom{\G} =
    \dom{\G'}$, $\G \sqsubseteq \G'$, $U \sqsubseteq U'$ and $V'
    \sqsubseteq V$. Hence $N :\<\Delta \v U'\>$ and, by IH, $M[x^L:=N]
    : \<\G'\sqcap \Delta\v V'\>$. It is easy to show that $\G\sqcap
    \Delta \sqsubseteq \G'\sqcap \Delta$. Hence, $\<\G'\sqcap \Delta\v
    V'\>\sqsubseteq \<\G\sqcap \Delta\v V\>$ and $M[x^L:=N] :
    \<\G\sqcap \Delta\v V\>$.
\hfill $\Box$
\end{itemize}
\end{proof}

The next lemma is needed in the proofs.
\begin{lemma}\label{restlem}
\begin{enumerate}
\item\label{restlemone}
If $\fv{N} \subseteq \fv{M}$, then $ env^M_\o \r_N = env^N_\o$.
\item\label{restlemtwo} If $\legalenv{\G_1}$, $\legalenv{\G_2}$,
  $\fv{M} \subseteq \dom{\G_1}$ and $\fv{N} \subseteq \dom{\G_2}$,
  then $(\G_1 \sqcap\G_2 )\r_{MN} \sqsubseteq (\G_1 \r_M)\sqcap (\G_2
  \r_N)$.
\item\label{restlemthree}
$\overline{e}_i(\G\r_M) = (\overline{e}_i\G)\r_{M^{+i}}$
\end{enumerate}
\end{lemma}

\begin{proof}
  \ref{restlemone}.\ Easy.  \ref{restlemtwo}.\ First, note that
  $\legalenv{\G_1 \sqcap \G_2}$ by
  lemma~\ref{env-Phisub}.\ref{lem:interlegalenv}, $\legalenv{\G_1
    \r_M}$, $\legalenv{\G_2 \r_N}$ and $\dom{(\G_1 \sqcap\G_2
  }\r_{MN}) = \fv{MN} = \fv{M} \cup \fv{N} = \dom{\G_1 \r_M} \cup
  \dom{\G_2\r_N} = \dom{(\G_1 \r_M}\sqcap(\G_2 \r_N))$. Now, we show
  by cases that if $(x^L:U_1) \in (\G_1 \sqcap\G_2 )\r_{MN}$ and
  $(x^L:U_2) \in (\G_1 \r_M)\sqcap(\G_2 \r_N)$ then $U_1 \sqsubseteq
  U_2$:
\begin{itemize}
\item If $x^L \in \fv{M}\cap \fv{N}$ then $(x^L:U'_1) \in \G_1$,
  $(x^L:U''_1) \in \G_2$ and $U_1 = U'_1\sqcap U''_1 = U_2$.
\item If $x^L \in \fv{M}\setminus \fv{N}$ then
\begin{itemize}
\item If $x^L \in \dom{\G_2}$ then  $(x^L:U_2) \in \G_1$, $(x^L:U'_1) \in
  \G_2$ and  $U_1 = U'_1 \sqcap U_2\sqsubseteq U_2$.
\item If $x^L \not \in \dom{\G_2}$ then $(x^L:U_2) \in \G_1$ and $U_1 = U_2$.
\end{itemize}
\item If $x^L \in \fv{N}\setminus \fv{M}$ then
\begin{itemize}
\item If $x^L \in \dom{\G_1}$ then  $(x^L:U'_1) \in \G_1$, $(x^L:U_2) \in
  \G_2$ and  $U_1 = U'_1 \sqcap U_2\sqsubseteq U_2$.
\item If $x^L \not \in \dom{\G_1}$ then $x^L:U_2 \in \G_2$ and $U_1 = U_2$.
\end{itemize}
\ref{restlemthree}.\ Let $\G = (x^{L_j}_j : U_j)_n$ and let $\fv{M} =
\{y^{K_1}_1,\dots,y^{K_m}_m\}$ where $m \leq n$ and $\forall 1 \leq k \leq
m$ $\exists 1 \leq j \leq n$ such that $y^{K_k}_k = x^{L_j}_j$. So
$\G\r_M = (y^{K_k}_k : U_k)_m$ and $\overline{e}_i(\G\r_M) = (y^{i::K_k}_k :
\overline{e}_iU_k)_m$. Since $\overline{e}_i\G = (x^{i::L_j}_j : \overline{e}_iU_j)_n$, $\fv{M^{+i}} =
\{y^{i::K_1}_1,\dots,y^{i::K_m}_m\}$ and $\forall 1 \leq k \leq
m$ $\exists 1 \leq j \leq n$ such that $y^{i::K_k}_k = x^{i::L_j}_j$ then
$(\overline{e}_i\G)\r_{M^{+i}} = (y^{i::K_k}_k : U_k)_m$.
\hfill $\Box$
\end{itemize}
\end{proof}

The next two theorems are needed in the proof of subject reduction.
\begin{theorem}\label{betatheo}
If $M : \<\G \v  U\>$ and $M \rhd_{\beta} N$, then $N : \<\G\r_N
\v  U\>$.
\end{theorem}

\begin{proof}
By induction on the derivation $M : \<\G \v  U\>$.
\begin{itemize}
\item
  Rule $\o$ follows by theorem~\ref{deg=}.\ref{deg=one} and
  lemma~\ref{restlem}.\ref{restlemone}.
\item If $\F{M : \<\G,(x^L:U) \v T\>}{\l x^L. M : \<\G \v U \f T\>}$
  then $N = \l x^L. N'$ and $M \rhd_{\beta} N'$. By IH, $N' :
  \<(\G,(x^L:U))\r_{N'} \v T\>$. If $x^L \in \fv{N'}$ then $N' :
  \<\G\r_{\fv{N'}\setminus \{x^L\}},(x^L:U) \v T\>$ and by $\f_I$, $\l
  x^L. N' : \<\G\r_{\l x^L. N'} \v U \f T\>$. Else $N' :
  \<\G\r_{\fv{N'}\setminus \{x^L\}} \v T\>$ so by $\f'_I$, $\l x^L. N'
  : \<\G\r_{\l x^L. N'} \v \o^L \f T\>$ and since by
  lemma~\ref{degnew-goodeg}.\ref{goodegtwo}, $U \sqsubseteq \o^{L}$,
  by $\sqsubseteq$, $\l x^L. N' : \<\G\r_{\l x^L. N'} \v U \f T\>$.
\item If $\F{M : \<\G \v T\>\;\;\; x^L \not \in \dom{\G}}{\l x^L. M :
    \<\G \v \o^L \f T\>}$ then $N = \l x^L N'$ and $M \rhd_{\beta}
  N'$. Since $x^L \not \in \fv{M}$, by
  theorem~\ref{deg=}.\ref{deg=one}, $x^L \not \in \fv{N'}$. By IH, $N'
  : \<\G\r_{\fv{N'}\setminus \{x^L\}} \v T\>$ so by $\f'_I$, $\l
  x^L. N' : \<\G\r_{\l x^L. N'} \v \o^L \f T\>$.
\item If $\F{M_1 : \<\G_1 \v U \f T\> \;\;\; M_2 : \<\G_2 \v U\>
    \;\;\; \G_1 \diamond \G_2} {M_1 \; M_2 : \<\G_1 \sqcap \G_2 \v
    T\>}$. Using lemma~\ref{restlem}.\ref{restlemtwo}, case $M_1
  \rhd_{\beta} N_1$ and $N = N_1M_2$ and case $M_2 \rhd_{\beta} N_2$
  and $N = M_1N_2$ are easy. Let $M_1 = \l x^L. M'_1$ and $N =
  M'_1[x^L := M_2]$.  By lemma~\ref{structyping}.\ref{structypingC}
  and lemma~\ref{degsub}.\ref{degsubone}, $M'_1 \diamond M_2$.  If
  $x^L \in FV (M'_1)$ then by lemma~\ref{newgen}.\ref{newgentwo},
  $M'_1 : \<\G_1, x^L : U \v T\>$. By lemma~\ref{substlem}, $M'_1[x^L
  := M_2] : \<\G_1 \sqcap \G_2 \v T\>$. If $x^L \not \in FV (M'_1)$
  then by lemma~\ref{newgen}.\ref{newgenthree}, $M'_1[x^L := M_2] =
  M'_1 : \<\G_1 \v T\>$ and by $\sqsubseteq$, $N : \<(\G_1 \sqcap
  \G_2)\r_N \v T\>$.
\item Case $\sqcap_I$ is by IH.
\item If $\F{M : \<\G \v U\>}{M^{+i} : \< \overline{e}_i\G \v \overline{e}_iU\>}$ and
  $M^{+i} \rhd_{\beta} N$, then by lemma~\ref{deg+-f+}.\ref{six},
  there is $P \in {\cal M}$ such that $P^{+i} = N$ and $M \rhd_{\beta}
  P$. By IH, $P : \<\G\r_P \v U\>$ and by $e$ and
  lemma~\ref{restlem}.\ref{restlemthree}, $N: \<(\overline{e}_i\G)\r_N \v
  \overline{e}_iU\>$.
\item If $\F{M : \<\G\v U\> \;\;\;\hspace{0.2in} \<\G\v U\>
    \sqsubseteq \<\G'\v U'\>}{M : \<\G'\v U'\>}$ then by IH,
  lemma~\ref{env-Phisub}.\ref{Phisubthree} and $\sqsubseteq$, $N :
  \<\G'\r_N \v U'\>$.
\hfill $\Box$
\end{itemize}
\end{proof}

\begin{theorem}\label{etatheo}
If $M : \<\G \v  U\>$ and $M \rhd_{\eta} N$, then $N : \<\G \v  U\>$.
\end{theorem}

\begin{proof}
By induction on the derivation $M : \<\G \v  U\>$.
\begin{itemize}
\item If $\F{}{M : \<env^\o_M \v \o^{\deg(M)}\>}$ then by
  lemma~\ref{deg=}.\ref{deg=two}, $\deg(M) = \deg(N)$ and $\fv{M} =
  \fv{N}$ and by $\o$, $N : \<env^\o_M \v \o^{\deg(M)}\>$.
\item If $\F{M : \<\G,(x^L:U) \v T\>}{\l x^L. M : \<\G \v U \f T\>}$
  then we have two cases:
  \begin{itemize}
    \item $M = Nx^L$ and so by lemma~\ref{newgen}.\ref{newgenfour'},
      $N : \<\G \v U \f T\>$.
    \item $N = \l x^L. N'$ and $M \rhd_{\eta} N'$. By IH, $N' :
      \<\G,(x^L:U) \v T\>$ and by $\f_I$, $N : \<\G \v U \f T\>$.
  \end{itemize}
\item if $\F{M : \<\G \v T\>\;\;\; x^L \not \in \dom{\G}}{\l x^L. M :
    \<\G \v \o^L \f T\>}$ then $N = \l x^L. N'$ and $M \rhd_{\eta}
  N'$. By IH, $N' : \<\G \v T\>$ and by $\f'_I$, $N : \<\G \v
  \o^L \f T\>$.
\item If $\F{M_1 : \<\G_1 \v U \f T\> \;\;\; M_2 : \<\G_2 \v U\>
    \;\;\; \G_1 \diamond \G_2}{M_1 M_2 : \<\G_1 \sqcap \G_2 \v
    T\>}$, then we have two cases:
  \begin{itemize}
  \item $M_1 \rhd_{\eta} N_1$ and $N = N_1M_2$. By IH $N_1 : \<\G_1 \v
    U \f T\>$ and by $\f_E$, $N : \<\G_1 \sqcap \G_2 \v T\>$.
  \item $M_2 \rhd_{\eta} N_2$ and $N = M_1N_2$. By IH $N_2 : \<\G_2 \v
    U\>$ and by $\f_E$, $N : \<\G_1 \sqcap \G_2 \v T\>$.
  \end{itemize}
\item Case $\sqcap_I$ is by IH and $\sqcap_I$.
\item If $\F{M : \<\G \v U\>}{M^{+i} : \< \overline{e}_i\G \v \overline{e}_iU\>}$ then by
  lemma~\ref{deg+-f+}.\ref{six}, there is $P
  \in {\cal M}$ such that $P^{+i} = N$ and $M \rhd_{\eta} P$. By IH, $P
  : \<\G \v U\>$ and by $e$, $N: \<\overline{e}_i\G \v \overline{e}_iU\>$.
\item If $\F{M : \<\G\v U\> \;\;\;\hspace{0.2in} \<\G\v U\>
    \sqsubseteq \<\G'\v U'\>}{M : \<\G'\v U'\>}$ then by IH,
  lemma~\ref{env-Phisub}.\ref{Phisubthree} and $\sqsubseteq$, $N : \<\G'
  \v U'\>$.
\hfill $\Box$
\end{itemize}
\end{proof}

The next auxiliary lemma is needed in proofs.
\begin{lemma}\label{envComp}
  Let $i \in \{1,2\}$ and $M: \<\G \v U\>$. We have:
  \begin{enumerate}
  \item\label{envComp1}
    If $(x^L:U_1) \in \G$ and $(y^K:U_2) \in \G$, then:
    \begin{enumerate}
    \item
      If $(x^L:U_1) \not = (y^K:U_2)$, then $x^L \not = y^K$.
    \item
      If $x = y$, then $L = K$ and $U_1 = U_2$.
    \end{enumerate}
  \item\label{envComp2}
    If $(x^L:U_1) \in \G$ and  $(y^K:U_2) \in \G$ and $(x^L:U_1) \not =
    (y^K:U_2)$, then $x \not = y$ and $x^L \not = y^K$.
  \end{enumerate}
\end{lemma}

\begin{proof}
\ref{envComp1}.\ If $x^{L} = Y^{K}$ then by definition $U_1 = U_2$, so
$(x^{L} : U_1) = (y^{K} : U_2)$. By
lemma~\ref{structyping}.\ref{structypingone}, $x^{L}, y^{K} \in
\fv{M}$. By lemma~\ref{degsub}.\ref{lem:refl+symm}, $M \diamond
M$. So, if $x = y$ then $L = K$ and by definition $U_1 = U_2$.
\ref{envComp2}.\ Corollary of \ref{envComp1}.
\hfill $\Box$
\end{proof}

\begin{proof}[Of theorem~\ref{subred}]
  Proofs are by induction on derivations using theorem~\ref{betatheo}
  and theorem~\ref{etatheo}.
\hfill $\Box$
\end{proof}

\section{Proofs for section~\ref{sexpsec}}

\begin{proof}[Of lemma~\ref{exp1}]
By lemma~\ref{typing'}.\ref{typing'one}, $M[x^{L}:=N] \in {\cal M}$,
so by definition, $M, N \in {\cal M}$ and $M \diamond N$ and $\deg(N)
= L$.
By induction on the derivation $M[x^L:= N] : \<\G \v U\>$.
\begin{itemize}
\item If $\F{}{y^{\oslash} : \<(y^{\oslash}:T) \v T\>}$ then $M =
  x^{\oslash}$ and $N = y^{\oslash}$. By $ax$, $x^{\oslash} :
  \<(x^{\oslash}:T) \v T\>$.
\item If $\F{}{M[x^L:= N] : \<env^\o_{M[x^L:= N]} \v \o^{\deg(M[x^L:=
      N])}\>}$ then by lemma~\ref{degsub}.\ref{degsubthree}, $\deg(M)
  = \deg(M[x^L:= N])$. By $\o$, $M :
  \<env^\o_{\fv{M}\setminus\{x^L\}},(x^L : \o^L) \v \o^{\deg(M)}\>$
  and $N : \<env^\o_N \v \o^L\>$ and because $\fv{M[x^{L}:=N]} =
  (\fv{M} \setminus \{x^{L}\}) \cup \fv{N}$,
  $env^\o_{\fv{M}\setminus\{x^L\}} \sqcap env^\o_N = env^\o_{M[x^L:=
    N]}$.
\item If $\F{M[x^L:=N] : \<\G, (y^K: W) \v T\>}{\l y^K. M[x^L:=N] :
    \<\G \v W \f T\>}$ where $y^K \not \in \fv{N} \cup \{x^{L}\}$. By
  IH, $\exists~V$ and $\exists~\G_1,\G_2$ type environments such that
  $M : \<\G_1, x^L : V \v T\>$, $N : \<\G_2 \v V\>$ and $\G, y^K : W =
  \G_1 \sqcap \G_2$. By lemma~\ref{structyping}.\ref{structypingone},
  $\fv{N} = \dom{\G_2}$ and $\fv{M} = \dom{\G_1} \cup \{y^{K}\}$.
  Since $y^K \in \fv{M}$ and $y^K \not \in \fv{N}$, $\G_1 = \D_1, y^K
  : W$. Hence $M : \<\D_1, y^K : W, x^L : V \v T\>$. By rule $\f_I$,
  $\l y^K. M : \<\D_1 , x^L : V \v W \f T\>$. Finally $\G = \D_1
  \sqcap \G_2$.
\item If $\F{M[x^L:=N] : \<\G \v T\> \;\;\;y^K \not \in \dom{\G}}{\l
    y^K.  M[x^L:=N] : \<\G \v \o^K \f T\>}$ where $y^K \not \in \fv{N}
  \cup \{x^{L}\}$. By IH, $\exists~V$ type and $\exists~\G_1,\G_2$
  type environments such that $M : \<\G_1, x^L : V \v T\>$, $N :
  \<\G_2 \v V\>$ and $\G = \G_1 \sqcap \G_2$. Since $y^K \neq x^L$,
  $\l y^K.M : \<\G_1, x^L : V \v \o^K \f T\>$.
\item
  If $\F{M_1[x^L:=N] : \<\G_1 \v  W \f T\> \;\;\; M_2[x^L:=N] :
    \<\G_2 \v W\> \;\;\; \G_1 \diamond \G_2} {M_1[x^L:=N] \;
    M_2[x^L:=N] : \<\G_1 \sqcap \G_2 \v T\>}$ where $M = M_1 M_2$,
  then we have three cases:


  \begin{itemize}
  \item If $x^L \in \fv{M_1}\cap \fv{M_2}$ then by IH,
    $\exists~V_1,V_2$ types and $\exists~\D_1,\D_2,\B_1,\B_2$ type
    environments such that $M_1 : \<\D_1, (x^L : V_1) \v W \f T\>$,
    $M_2 : \<\B_1, (x^L : V_2) \v W\>$, $N : \<\D_2 \v V_1\>$, $N :
    \<\B_2 \v V_2\>$, $\G_1 = \D_1 \sqcap \D_2$ and $\G_2 = \B_1
    \sqcap \B_2$. Because $\G_1 \diamond \G_2$, then $\D_1 \diamond
    \B_1$ and $\D_2 \diamond \B_2$ and because $\D_1, (x^L : V_1)$ and
    $\B_1, (x^L : V_2)$ are type environments, by lemma~\ref{envComp},
    $(\D_1, (x^L : V_1)) \diamond (\B_1, (x^L : V_2))$. Then, by rules
    $\sqcap_I$ and $\f_E$, $M_1 M_2: \<\D_1 \sqcap \B_1, (x^L : V_1
    \sqcap V_2) \v T\>$ and by $\sqcap'_I$, $N : \<\Delta_2 \sqcap
    \B_2 \v V_1 \sqcap V_2\>$. Finally, $\G_1 \sqcap \G_2 = (\D_1
    \sqcap \D_2) \sqcap (\B_1 \sqcap \B_2)$.
  \item If $x^L \in \fv{M_1}\setminus \fv{M_2}$ then by IH, $\exists~V$
    types and $\exists~\D_1,\D_1$ type environments such that $M_1 :
    \<\D_1, (x^L : V) \v W \f T\>$, $N : \<\D_2 \v V\>$ and $\G_1 =
    \D_1 \sqcap \D_2$. Since $\G_1 \diamond \G_2$, $\D_1 \diamond
    \G_2$ and since $\G_1 \sqcap \G_2$ is a type environment, by
    lemma~\ref{envComp}, $(\D_1, (x^L : V)) \diamond \G_2$. By $\f_E$,
    $M_1 M_2: \<\D_1 \sqcap \G_2, (x^L : V) \v T\>$ and $\G_1 \sqcap
    \G_2 = (\D_1 \sqcap \D_2) \sqcap \G_2$.
  \item If $x^L \in \fv{M_2}\setminus \fv{M_1}$ then by IH, $\exists~V$
    types and $\exists~\D_1,\D_2$ type environments such that $M_2 :
    \<\D_1, (x^L : V) \v W\>$, $N : \<\D_2 \v V\>$ and $\G_2 = \D_1
    \sqcap \D_2$. Since $\G_1 \diamond \G_2$, $\G_1 \diamond \D_1$ and
    since $\G_1 \sqcap \G_2$ is a type environment, by
    lemma~\ref{envComp}, $(\D_1, (x^L : V)) \diamond \G_1$. By $\f_E$,
    $M_1 M_2: \<\G_1 \sqcap \D_1, (x^L : V) \v T\>$ and $\G_1 \sqcap
    \G_2 = \G_1 \sqcap (\D_1 \sqcap \D_2)$.
  \end{itemize}
\item
  Let $\F{M[x^L:=N] : \<\G \v  U_1\> \;\; M[x^L:=N]  : \<\G \v U_2\>}
  {M[x^L:=N] : \<\G  \v  U_1 \sqcap U_2\>}$. By IH, $\exists~V_1,V_2$
  types and $\exists~\D_1,\D_2,\B_1,\B_2$ type environments such that
  $M : \<\D_1, x^L : V_1 \v U_1\>$, $M : \<\B_1, x^L : V_2 \v U_2\>$,
  $N : \<\D_2 \v V_1\>$, $N : \<\B_2 \v V_2\>$, $\G = \D_1
  \sqcap \D_2$ and $\G = \B_1 \sqcap \B_2$. Then, by rule
  $\sqcap'_I$, $M : \<\D_1 \sqcap \B_1, x^L : V_1
  \sqcap V_2 \v U_1 \sqcap U_2\>$ and $N : \<\D_2 \sqcap \B_2 \v V_1
  \sqcap V_2\>$. Finally, $\G = (\D_1 \sqcap \D_2) \sqcap (\B_1 \sqcap
  \B_2)$.
\item
  If $\F{M[x^L:=N] : \<\G \v U\>}{M^{+j}[x^{j::L}:=N^{+j}] : \< \overline{e}_j\G
  \v \overline{e}_jU\>}$ then by IH, $\exists~V$ type and $\exists~\G_1,\G_2$
  type environments such that $M : \<\G_1, x^L : V \v U\>$, $N :
  \<\G_2 \v V\>$ and $\G = \G_1 \sqcap \G_2$. So by $e$, $M^{+j} :
  \<\overline{e}_j\G_1, x^{j::L} : \overline{e}_jV \v \overline{e}_jU\>$, $N : \<\overline{e}_j\G_2 \v \overline{e}_jV\>$
  and $\overline{e}_j\G = \overline{e}_j\G_1 \sqcap \overline{e}_j\G_2$.
\item
  If $\F{M[x^L:=N] : \<\G' \v U'\> \;\;\; \<\G'\v  U'\> \sqsubseteq
    \<\G \v U\>}{M[x^L:=N] : \<\G \v U\>}$ then by
  lemma~\ref{env-Phisub}.\ref{Phisubtwo}, $\G \sqsubseteq \G'$ and $U'
  \sqsubseteq U$. By IH, $\exists~V$ type and $\exists~\G_1,\G_2$
  type environments such that $M : \<\G'_1, x^L : V \v U'\>$, $N :
  \<\G'_2 \v V\>$ and $\G' = \G'_1 \sqcap \G'_2$. Then by
  lemma~\ref{degnew-goodeg}.\ref{goodegfour}, $\G = \G_1 \sqcap \G_2$ and
  $\G_1 \sqsubseteq \G'_1$ and $\G_2 \sqsubseteq \G'_2$. So by
  $\sqsubseteq$, $M : \<\G_1, x^L : V \v U\>$ and $N : \<\G_2 \v
  V\>$.
\hfill $\Box\$
end{itemize}
\end{proof}

The next lemma is basic for the proof of subject expansion for $\beta$.
\begin{lemma}\label{exp2}
  If $M[x^L:=N] : \<\G \v U\>$, $\deg(U) = K$ and $L \succeq \deg(M)$,
  ${\cal U} = \fv{(\l x^L.M}N)$, then $(\l x^L.M)N :\<\G {\u^{\cal U}}
  \v U\>$.
\end{lemma}

\begin{proof}
  By lemma~\ref{typing'}.\ref{typing'one}, $M[x^{L}:=N] \in {\cal M}$,
  so $M, N \in {\cal M}$ and $M \diamond N$ and $\deg(N) = L$.  By
  definition $(\l x^{L}. M)N \in {\cal M}$.  By lemma~\ref{degsub}.\ref{degsubthree} and
  theorem~\ref{typing'}.\ref{typing'one}, $\deg(\G) \succeq \deg(U) =
  K = \deg(M[x^L:=N]) = \deg(M) = \deg((\l x^L.M)N)$.  So $L \succeq
  K$ and there exists $K'$ such that $L = K::K'$.  We have two cases:
  \begin{itemize}
  \item If $x^L \in \fv{M}$, then, by lemma~\ref{exp1}, $\exists~V$
    type and $\exists~\G_1,\G_2$ type environments such that $M :
    \<\G_1, x^L : V \v U\>$, $N : \<\G_2 \v V\>$ and $\G = \G_1 \sqcap
    \G_2$.  By lemma~\ref{typing'}.\ref{typing'one}, $\legalenv{\G_1}$
    and $\legalenv{\G_2}$.  By
    lemma~\ref{env-Phisub}.\ref{lem:interlegalenv}, $\legalenv{\G_1
      \sqcap \G_2}$.  So, it is easy to prove, using
    lemma~\ref{env-Phisub}.\ref{lem:legalomegaenv}, that
    $\legalenv{\G{\u^{\cal U}}}$.  By
    lemma~\ref{structyping}.\ref{structypingC}, $\G_1, x^L : V
    \diamond \G_2$, so $\G_1 \diamond \G_2$.  By
    lemma~\ref{typing'}.\ref{typing'one}, $\deg(\G_1) \succeq \deg(M)
    = \deg(U) = K$ and $L = \deg(N) = \deg(V) \preceq \deg(\G_2)$. By
    lemma~\ref{degnew-goodeg}, we have two cases :
    \begin{itemize}
    \item If $U = \o^K$, then by
      lemma~\ref{structyping}.\ref{structypingtwo}, $(\l x^L.M)N :\<\G
      {\u^{\cal U}} \v U\>$.
    \item If $U = \vec{e}_K \sqcap_{i=1}^p T_i$ where $p \geq 1$ and
      $\forall~1 \leq i \leq p$, $T_i \in {\mathbb T}$, then by
      theorem~\ref{typing'}.\ref{typing'two'}, $M^{-K} : \<\G_1^{-K},
      x^{K'} : V^{-K} \v \sqcap_{i=1}^p T_i\>$. By $\sqsubseteq$,
      $\forall~1 \leq i \leq p$, $M^{-K} : \<\G_1^{-K}, x^{K'} :
      V^{-K} \v T_i\>$, so by $\f_I$, $\l x^{K'}. M^{-K} : \<\G_1^{-K}
      \v V^{-K} \f T_i\>$. Again by
      theorem~\ref{typing'}.\ref{typing'two'}, $N^{-K} : \<\G_2^{-K}
      \v V^{-K}\>$ and since $\G_1 \diamond \G_2$, by
      lemma~\ref{env-Phisub}.\ref{Phisubfive}, $\G_1^{-K} \diamond
      \G_2^{-K}$, so by $\f_E$, $\forall~1 \leq i \leq p$, $(\l
      x^{K'}. M^{-K})N^{-K} : \<\G_1^{-K} \sqcap \G_2^{-K} \v
      T_i\>$. Finally by $\sqcap_I$ and $e$, $(\l x^L. M)N : \<\G_1
      \sqcap \G_2 \v U\>$, so $(\l x^L.M)N :\<\G {\u^{\cal U}} \v
      U\>$.
    \end{itemize}
  \item If $x^L \not \in \fv{M}$, then $M : \<\G \v U\>$.  By
    lemma~\ref{typing'}.\ref{typing'one}, $\legalenv{\G}$.  So, it is
    easy to prove, using
    lemma~\ref{env-Phisub}.\ref{lem:legalomegaenv}, that
    $\legalenv{\G{\u^{\cal U}}}$.  By lemma~\ref{degnew-goodeg}, we
    have two cases :
    \begin{itemize}
    \item If $U = \o^K$, then by
      lemma~\ref{structyping}.\ref{structypingtwo}, $(\l x^L.M)N :\<\G
      {\u^{\cal U}} \v U\>$.
    \item If $U = \vec{e}_K \sqcap_{i=1}^p T_i$ where $p \geq 1$ and
      $\forall~1 \leq i \leq p$, $T_i \in {\mathbb T}$, then by
      theorem~\ref{typing'}.\ref{typing'two'}, $M^{-K} : \<\G^{-K} \v
      \sqcap_{i=1}^p T_i\>$. By $\sqsubseteq$, $\forall~1 \leq i \leq
      p$, $M^{-K} : \<\G^{-K} \v T_i\>$. Using lemma~\ref{deg+-f+} and
      by induction on $K$, we can prove that $x^{K'} \not \in
      \fv{M^{-K}}$. So by
      lemma~\ref{structyping}.\ref{structypingone}, $x^{K'} \not \in
      \dom{\G^{-K}}$.  So by $\f'_I$, $\l x^{K'}. M^{-K} : \<\G^{-K}
      \v \o^{K'} \f T_i\>$.  By $(\o)$, $N^{-K} : \<env^{\o}_{N^{-K}}
      \v \o^{K'}\>$ and $N : \<env^{\o}_{N} \v \o^{L}\>$.  By
      theorem~\ref{typing'}.\ref{typing'one}, $\deg(env^{\o}_{N})
      \succeq \deg(N) = L$.  By
      lemma~\ref{structyping}.\ref{structypingC}, $\G \diamond
      env^{\o}_{N}$. By lemma~\ref{env-Phisub}.\ref{Phisubfive},
      $\G^{-K} \diamond env^{\o}_{N^{-K}}$.  By $\f_E$, $\forall~1
      \leq i \leq p$, $(\l x^{K'}. M^{-K})N^{-K} : \<\G^{-K} \sqcap
      env^{\o}_{N^{-K}} \v T_i\>$.  Finally by $\sqcap_I$ and $e$,
      $(\l x^L. M)N : \<\G \sqcap env^{\o}_{N} \v U\>$, so $(\l
      x^L.M)N :\<\G {\u^{\cal U}} \v U\>$.
\hfill $\Box$
    \end{itemize}
  \end{itemize}
\end{proof}

Next, we give the main block for the proof of subject expansion for $\beta$.
\begin{theorem}\label{betaexp}
If $N : \<\G \v  U\>$ and $M \rhd_{\beta} N$, then $M :
\<\G {\u^{M}} \v U\>$.
\end{theorem}

\begin{proof}By induction on the derivation $N : \<\G \v  U\>$.
\begin{itemize}
\item If $\F{}{x^{\oslash} : \<(x^{\oslash}:T) \v T\>}$ and $M
  \rhd_{\beta} x^{\oslash}$, then $M = (\l y^K. M_1)M_2$ and
  $x^{\oslash} = M_1[y^K:=M_2]$.  Because $M \in {\cal M}$ then $K
  \succeq \deg(M_1)$.  By lemma~\ref{exp2}, $M : \<(x^{\oslash}:T)
  {\u^{M}} \v T\>$.
\item If $\F{}{N : \<env_N^\o \v \o^{\deg(N)}\>}$ and $M \rhd_{\beta}
  N$, then since by theorem~\ref{deg=}.\ref{deg=one}, $\fv{N}
  \subseteq \fv{M}$ and $\deg(M) = \deg(N)$, $(env_N^\o) {\u^{M}} =
  env_M^\o$.  By $\o$, $M : \<env_M^\o \v \o^{\deg(M)}\>$.  Hence, $M
  : \<(env^N_\o) {\u^{M}} \v \o^{\deg(N)}\>$.
\item If $\F{N : \<\G , x^L :U \v T\>}{\l x^L. N : \<\G \v U \f T\>}$
  and $M \rhd_{\beta} \l x^L. N$, then we have two cases:
  \begin{itemize}
  \item If $M = \l x.M'$ where $M' \rhd_{\beta} N$, then by IH, $M' :
    \<(\G, (x^L:U)) {\u^{M'}} \v T\>$. Since by
    theorem~\ref{deg=}.\ref{deg=one} and
    lemma~\ref{structyping}.\ref{structypingone}, $x^L \in \fv{N}
    \subseteq \fv{M'}$, then we have $(\G, (x^L:U)) {\u^{\fv{M'}}} =
    \G {\u^{\fv{M'}\setminus\{x^L\}}}, (x^L:U)$ and $\G
    {\u^{\fv{M'}\setminus\{x^L\}}} = \G {\u^{\l x^L.M'}}$. Hence, $M'
    : \<\G {\u^{\l x^L.M'}}, (x^L:U) \v T\>$ and finally, by $\f_I$,
    $\l x^L.M' :\<\G {\u^{\l x^L.M'}} \v U \f T\>$.
  \item If $M = (\l y^K. M_1)M_2$ where $y^K \not \in \fv{M_2}$ and
    $\l x^L. N = M_1[y^K:=M_2]$, then, because $M \in {\cal M}$ then
    $K \succeq \deg(M_1)$ and by lemma~\ref{exp2}, Because
    $M_1[y^K:=M_2] : \<\G \v U \f T\>$, we have $(\l y^K. M_1)M_2 :
    \<\G {\u^{(\l y^K. M_1)M_2}} \v U \f T\>$.
  \end{itemize}
\item
  If $\F{N : \<\G \v T\>\;\;\;x^L \not \in \dom{\G}}{\l x^L. N : \<\G
  \v \o^L \f T\>}$ and $M \rhd_{\beta} N$  then similar to the above
  case.
\item
  If $\F{N_1 : \<\G_1 \v U \f T\> \;\;\; N_2 : \<\G_2 \v U\>
    \;\;\; \G_1 \diamond \G_2}{N_1 \; N_2 : \<\G_1 \sqcap \G_2 \v
    T\>}$ and $M \rhd_{\beta} N_1N_2$, we have three cases:
  \begin{itemize}
  \item $M = M_1N_2$ where $M_1 \rhd_{\beta} N_1$ and $M_1 \diamond
    N_2$. By IH, $M_1 : \<\G_1 {\u^{M_1}} \v U \f T\>$.  It is easy to
    show that $(\G_1 \sqcap \G_2) {\u^{M_1N_2}} = \G_1
    {\u^{M_1}}\sqcap \G_2$. Since $M_1 \diamond N_2$, by
    lemma~\ref{structyping}.\ref{structypingC}, $\G_1 {\u^{M_1}}
    \diamond \G_2$, hence use $\f_E$.
  \item
    $M = N_1M_2$ where $M_2 \rhd_{\beta} N_2$. Similar to the above
    case.
  \item If $M = (\l x^L.M_1)M_2$ and $N_1N_2 = M_1[x^L:=M_2]$ then,
    because $M \in {\cal M}$ then $L \succeq \deg(M_1)$ and by
    lemma~\ref{exp2}, $(\l x^L.M_1)M_2 :\<(\G_1 \sqcap \G_2) {\u^{(\l
        x^L.M_1)M_2}} \v T\>$.
  \end{itemize}
\item
  If $\F{N: \<\G \v U_1\> \;\;\;\hspace{0.2in} N : \<\G \v U_2\>}
  {N : \<\G \v U_1 \sqcap U_2\>}$ and $M \rhd_{\beta} N$ then use IH.
\item
  If $\F{N : \<\G \v U\>}{N^{+j} : \< \overline{e}_j\G \v \overline{e}_jU\>}$ then by
  lemma~\ref{deg+-f+}.\ref{five} then there is $P \in {\cal M}$ such
  that $M = P^{+j}$ and $P \rhd_{\beta} N$. By IH, $P : \<\G{\u^{P}}
  \v U\>$ and by $e$, $M : \< (\overline{e}_j\G) {\u^{M}} \v \overline{e}_jU\>$.
\item
  If $\F{N : \<\G\v U\> \;\;\;\hspace{0.2in} \<\G\v U\> \sqsubseteq
  \<\G'\v U'\>}{N : \<\G'\v U'\>}$ and $M \rhd_{\beta} N$.  By
  lemma~\ref{env-Phisub}.\ref{Phisubthree}, $\G' \sqsubseteq \G$ and $U
  \sqsubseteq U'$. It is easy to show that $\G' {\u^{M}} \sqsubseteq
  \G {\u^{M}}$ and hence by  lemma~\ref{env-Phisub}.\ref{Phisubthree},
  $\<\G {\u^{M}}\v U\> \sqsubseteq \<\G' {\u^{M}}\v U'\>$. By IH, $M
  {\u^{M}} : \<\G\v U\>$.  Hence, by $\sqsubseteq_{\<\>}$, we have $M
  : \<\G' {\u^{M}}\v U'\>$.
\hfill $\Box$
\end{itemize}
\end{proof}

\begin{proof}[Of theorem~\ref{finalbetaexp}]
By induction on the length of the derivation $M \rhd^*_{\beta} N$
using theorem~\ref{betaexp} and the fact that if $\fv{P} \subseteq
\fv{Q}$, then $(\G {\u^{P}}) {\u^{Q}} = \G {\u^{Q}}$.
\hfill $\Box$
\end{proof}

\section{Proofs of section~\ref{redsec}}

\begin{proof}[Of lemma~\ref{fx+}]
  \ref{fx+one}.\ and \ref{fx+two}.\ are easy.\\
  \ref{fx+three}.\ If $M\rhd_r^* N^{+i}$ where $N \in {\cal X}$, then,
  by lemma~\ref{deg+-f+}.\ref{five}, $M = P^{+i}$ such that $P \in
  {\cal M}$ and $P \rhd_r N$. As ${\cal X}$ is $r$-saturated, $P\in
  {\cal X}$ and so $P^{+i} = M \in {\cal X}^{+i}$.\\
  \ref{fx+four}.\ Let $M \in {\cal X} \fx {\cal Y}$ and $N \rhd_r^*
  M$. If $P \in {\cal X}$ such that $P \diamond N$, then by
  lemma~\ref{deg+-f+}.\ref{five'}, $P \diamond M$.  So, by definition,
  $MP \in {\cal Y}$.  Because ${\cal Y} \subseteq {\cal M}$, then $MP
  \in {\cal M}$.  Hence, $\deg(M) \preceq \deg(P)$.  By
  lemma~\ref{deg=}, $\deg(M) = \deg(N)$.  So $NP \in {\cal M}$ and $NP
  \rhd^{*}_r MP$.  Because $MP \in {\cal Y}$ and ${\cal Y}$ is
  $r$-saturated, then $NP \in {\cal Y}$. Hence, $N \in {\cal X} \fx
  {\cal Y}$.\\
  \ref{fx+five}.\ Let $M \in ({\cal X} \fx {\cal Y})^{+i}$, then $M =
  N^{+i}$ and $N \in {\cal X} \fx {\cal Y}$.  Let $P \in {\cal
    X}^{+i}$ such that $M \diamond P$.  Then $P = Q^{+i}$ such that $Q
  \in {\cal X}$.  Because $M \diamond P$ then by
  lemma~\ref{deg+-f+}.\ref{one''}, $N \diamond Q$. So $NQ \in {\cal
    Y}$. Because ${\cal Y} \subseteq {\cal M}$ then $NQ \in {\cal
    M}$. Because $(NQ)^{+i} = N^{+i}Q^{+i} = MP$ then $MP \in {\cal
    Y}^{+i}$. Hence, $M \in {\cal X}^{+i} \fx {\cal
    Y}^{+i}$.\\
  \ref{fx+six}.\ Let $M \in {\cal X}^{+i} \fx {\cal Y}^{+i}$ such that
  ${\cal X}^{+i} \wr {\cal Y}^{+i}$.  By hypothesis, there exists $P
  \in {\cal X}^{+i}$ such that $M \diamond P$.  Then $MP \in {\cal
    Y}^{+i}$.  Hence $MP = Q^{+i}$ such that $Q \in {\cal Y}$.
  Because ${\cal Y} \subseteq {\cal M}$ then $Q \in {\cal M}$ and by
  lemma~\ref{deg+-f+}.\ref{one'}, $MP \in {\cal M}$. Hence by
  definition $M \in {\cal M}$ and by lemma~\ref{deg+-f+}.\ref{one'},
  $\deg(M) = \deg(Q^{+i}) = i::\deg(Q)$. So by
  lemma~\ref{deg+-f+}.\ref{four}, there exists $M_1 \in {\cal M}$ such
  that $M = M_1^{+i}$.  Let $N_1 \in {\cal X}$ such that $M_1 \diamond
  N_1$.  By definition $N_1^{+i} \in {\cal X}^{+i}$ and by
  lemma~\ref{deg+-f+}.\ref{one''}, $M \diamond N_1^+$.  So, $MN_1^{+i}
  \in {\cal Y}^{+i}$. So $MN_1^{+i} = M'^{+i}$ such that $M' \in {\cal
    Y}$.  Because ${\cal Y} \subseteq {\cal M}$ then $M' \in {\cal
    M}$.  By lemma~\ref{deg+-f+}.\ref{one'}, $MN_1^{+i} \in {\cal M}$.
  So $M_1^{+i} \diamond N_1^{+i}$ and $\deg(M_1^{+i}) \preceq
  \deg(N_1^{+i})$.  By lemma~\ref{deg+-f+}.\ref{one'} and
  lemma~\ref{deg+-f+}.\ref{one''}, $M_1 \diamond N_1$ and $\deg(M_1)
  \preceq \deg(N_1)$. So $M_1N_1 \in {\cal M}$ and $(M_1 N_1)^{+i} =
  M_1^{+i} N_1^{+i} \in {\cal Y}^{+i}$.  Hence $M_1N_1 \in {\cal Y}$.
  Thus, $M_1 \in {\cal X} \fx {\cal Y}$ and $M = M_1^{+i} \in ({\cal
    X} \fx {\cal Y})^{+i}$.
\hfill $\Box$
\end{proof}

\begin{proof}[Of lemma~\ref{interpret-intsub}]
\ref{interpret}.\ref{interpretone} .\ By induction on $U$ using
lemma~\ref{fx+} and lemma~\ref{deg=}.\\
\ref{interpret}.\ref{interprettwo}.\ We prove $\forall x \in {\cal
  V}_1, {\cal N}_x^L \subseteq {\cal I}(U) \subseteq {\cal
M}^L$ by induction on $U$.
Case $U = a$: by definition.
Case $U = \o^L$: We have $\forall x \in {\cal V}_1, {\cal N}_x^L
\subseteq {\cal M}^L \subseteq {\cal M}^L$.
Case $U = U_1 \sqcap U_2$ (resp.\  $U = \overline{e}_i V$): use
IH since $\deg(U_1) = \deg(U_2)$ (resp.\ $\deg(U) = i::\deg(V)$,
$\forall x \in {\cal V}_1, ({\cal N}_x^K)^{+i} = {\cal N}_x^{i::K}$
and $({\cal M}^K)^{+i}
= {\cal M}^{i::K}$).
Case $U = V \f T$: by definition, $K=\deg(V)
\succeq \deg(T)=\oslash$.
\begin{itemize}
\item Let $x \in {\cal V}_1$, $N_1,...,N_k$ such that $\forall~1 \leq
  i \leq k$, $\deg(N_i) \succeq \oslash$ and $\diamond \{x^{\oslash},
  N_1, \dots, N_k\}$ and let $N \in {\cal I}(V)$ such that
  $(x^{\oslash} N_1... N_k) \diamond N$. By IH, $\deg(N) = K \succeq
  \oslash$. Again, by IH, $x^{\oslash} N_1... N_k N \in {\cal
    I}(T)$. Thus $x^{\oslash} N_1...N_k \in {\cal I}(V \f T)$.
\item Let $M \in {\cal I}(V \f T)$. Let $x \in {\cal V}_1$ such that
  $\forall L, x^L \not \in \fv{M}$. By IH, $x^K \in {\cal I}(V)$, then
  $Mx^K \in {\cal I}(T)$ and, by IH, $\deg(M x^K) = \oslash$. Thus
  $\deg(M) = \oslash$.
\end{itemize}
\ref{intsub}.\ By induction of the derivation $U \sqsubseteq V$.
\hfill $\Box$
\end{proof}

\begin{proof}[Of lemma~\ref{adeq}]
By induction on the derivation  $M : \<(x^{L_j}_j:U_j)_n \v U\>$.
\begin{itemize}
\item If $\F{}{x^{\oslash} : \<(x^{\oslash}:T) \v T\>}$ and $N \in
  {\cal I}(T)$, then $x^{\oslash}[x^{\oslash}:=N] = N \in {\cal
    I}(T)$.

\item If $\F{}{M : \<env^\o_M \v \o^{\deg(M)}\>}$. Let $env^\o_M =
  (x^{L_j}_j : U_j)_n$ so $\fv{M} =
  \{x^{L_1}_1,...,x^{L_n}_n\}$. Because, by
  lemma~\ref{typing'}.\ref{typing'one}, for all $j \in \{1, \dots,
  n\}$, $\deg(U_j) = L_j$ by
  lemma~\ref{interpret-intsub}.\ref{interpret}, ${\cal I}(U_j)
  \subseteq {\cal M}^{L_j}$, hence, $\deg(N_j) = L_j$. Because
  $M[(x^{L_j}_j:=N_j)_n] \in {\cal M}$, then $\diamond \{M\} \cup
  \{N_i$ / $i \in \{1, \dots, n\}\}$. Then, by
  lemma~\ref{degsub}.\ref{degsubthree}, $\deg(M[(x^{L_j}_j:=N_j)_n]) =
  \deg(M)$ and $M[(x^{L_j}_j:=N_j)_n]\in {\cal M}^{\deg(M)} = {\cal
    I}(\o^{\deg(M)})$.

\item If $\F{M : \< (x^{L_j}_j:U_j)_n ,(x^K:V) \v T\>} {\l x^K. M : \<
    (x^{L_j}_j:U_j)_n\v V \f T\>}$, $\forall 1 \leq j \leq n$, $N_j
  \in {\cal I}(U_j)$ and $N \in {\cal I}(V)$ such that $(\l
  x^K. M)[(x^{L_j}_j:=N_j)_n] \diamond N$. By
  lemma~\ref{typing'}.\ref{typing'one}, $\deg(V) = K$.  We have, $(\l
  x^K. M)[(x^{L_j}_j:=N_j)_n] = \l x^K. M[(x^{L_j}_j:=N_j)_n]$, where
  $\forall 1 \leq j \leq n, y^K \not \in \fv{N_j} \cup \{x^{L_j}_j\}$.
  Since $N \in {\cal I}(V)$ and by
  lemma~\ref{interpret-intsub}.\ref{interpret}, ${\cal I}(V) \subseteq
  {\cal M}^K$, $\deg(N) = K$.  By lemma~\ref{degsub}.\ref{degsubone}
  and lemma~\ref{degsub}.\ref{degsubthree}, $M[(x^{L_j}_j:=N_j)_n]
  \diamond N$ and $M[(x^{L_j}_j:=N_j)_n][x^{K}:=N] =
  M[(x^{L_j}_j:=N_j)_n,x^{K}:=N] \in {\cal M}$.  Hence, $(\l
  x^K. M[(x^{L_j}_j:=N_j)_n])N \in {\cal M}$ and $(\l
  x^K. M[(x^{L_j}_j:=N_j)_n])N \rhd_r
  M[(x^{L_j}_j:=N_j)_n,(x^K:=N)]$. By IH,
  $M[(x^{L_j}_j:=N_j)_n,(x^K:=N)] \in {\cal I}(T)$. Since, by
  lemma~\ref{interpret-intsub}.\ref{interpret} ${\cal I}(T)$ is
  $r$-saturated, then $(\l x^K. M[(x^{L_j}_j:=N_j)_n])N \in {\cal
    I}(T)$ and so $\l x^K. M[(x^{L_j}_j:=N_j)_n] \in {\cal I}(V) \fx
  {\cal I}(T) = {\cal I}(V \f T)$.

\item If $\F{M : \< (x^{L_j}_j:U_j)_n \v T\>\;\;\; x^K \not \in
    \dom{(x^{L_j}_j:U_j}_n)} {\l x^K. M : \<(x^{L_j}_j:U_j)_n\v \o^K
    \f T\>}$, $\forall 1 \leq j \leq n$, $N_j \in {\cal I}(U_j)$ and
  $N \in {\cal I}(\o^K)$ such that $(\l x^K. M)[(x^{L_j}_j:=N_j)_n]
  \diamond N$.  By lemma~\ref{structyping}.\ref{structypingone},
  $x^{K} \not \in \fv{M}$.  We have, $(\l x^K. M)[(x^{L_j}_j:=N_j)_n]
  = \l x^K. M[(x^{L_j}_j:=N_j)_n]$, where $\forall 1 \leq j \leq n,
  x^K \not \in \fv{N_j} \cup \{x^{L_j}_j\}$.  Since $N \in {\cal
    I}(\o^{K})$ and by lemma~\ref{interpret-intsub}.\ref{interpret},
  ${\cal I}(\o^{K}) = {\cal M}^K$ then $\deg(N) = K$.  By
  lemma~\ref{degsub}.\ref{degsubone} and
  lemma~\ref{degsub}.\ref{degsubthree}, $M[(x^{L_j}_j:=N_j)_n]
  \diamond N$ and $M[(x^{L_j}_j:=N_j)_n][x^{K}:=N] =
  M[(x^{L_j}_j:=N_j)_n,x^{K}:=N] = M[(x^{L_j}_j:=N_j)_n] \in {\cal
    M}$.  Hence, $(\l x^K. M[(x^{L_j}_j:=N_j)_n])N \in {\cal M}$ and
  $(\l x^K. M[(x^{L_j}_j:=N_j)_n])N \rhd_r
  M[(x^{L_j}_j:=N_j)_n,(x^K:=N)]$.  By IH, $M[(x^{L_j}_j:=N_j)_n] \in
  {\cal I}(T)$. Since, by lemma~\ref{interpret-intsub}.\ref{interpret}
  ${\cal I}(T)$ is $r$-saturated, then $(\l
  x^K. M[(x^{L_j}_j:=N_j)_n])N \in {\cal I}(T)$ and so $\l
  x^K. M[(x^{L_j}_j:=N_j)_n] \in {\cal I}(\o^{K}) \fx {\cal I}(T) =
  {\cal I}(\o^{K} \f T)$.

\item Let $\F{M_1 : \<\G_1 \v V \f T\> \;\;\; M_2 : \<\G_2 \v V\>
    \;\;\; \G_1 \diamond \G_2} {M_1 \; M_2 : \<\G_1 \sqcap \G_2 \v
    T\>}$ where $\G_1 = (x^{L_j}_j:U_j)_n, (y^{K_j}_j:V_j)_m$, $\G_2 =
  (x^{L_j}_j:U'_j)_n,(z^{S_j}_j:W_j)_p$ such that $\{y^{K_1}_1, \dots,
  y^{K_m}_m\} \cap \{z^{S_1}_1, \dots, z^{S_p}_p\} = \emptyset$ and
  $\G_1 \sqcap \G_2 = (x^{L_j}_j:U_j \sqcap
  U'_j)_n,(y^{K_j}_j:V_j)_m,(z^{S_j}_j:W_j)_p$.

  Let $\forall~1 \leq j \leq n, P_j \in {\cal I}(U_j \sqcap U'_j)$,
  $\forall~1 \leq j \leq m, Q_j \in {\cal I}(V_j)$ and $\forall~1 \leq
  j \leq p, R_j \in {\cal I}(W_j)$.  So, for all $j \in \{1, \dots,
  n\}$, $P_j \in {\cal I}(U_j)$ and $P_j \in {\cal I}(U'_j)$.  By
  hypothesis,
  $(M_1M_2)[(x^{L_j}_j:=P_j)_n,(y^{K_j}_j:=Q_j)_m,(z^{S_j}_j:=R_j)_p]
  = AB \in {\cal M}$ where using
  lemma~\ref{structyping}.\ref{structypingone}, $A =
  M_1[(x^{L_j}_j:=P_j)_n,(y^{K_j}_j:=Q_j)_m] \in {\cal M}$ and $B =
  M_2[(x^{L_j}_j:=P_j)_n,(z^{S_j}_j:=R_j)_p] \in {\cal M}$ and $A
  \diamond B$.

  By IH, $A \in {\cal I}(V) \fx {\cal I}(T)$ and $B \in {\cal I}(V)$.
  Hence, $AB \in {\cal I}(T)$.

\item Let $\F{M: \<(x^{L_j}_j:U_j)_n \v V_1\> \;\;\; M :
    \<(x^{L_j}_j:U_j)_n \v V_2\>} {M : \<(x^{L_j}_j:U_j)_n \v V_1
    \sqcap V_2\>}$. By IH, $M[(x^{L_j}_j:=N_j)_n] \in {\cal I}(V_1)$
  and $M[(x^{L_j}_j:=N_j)_n] \in {\cal I}(V_2)$. Hence,
  $M[(x^{L_j}_j:=N_j)_n] \in {\cal I}(V_1 \sqcap V_2)$.

\item Let $\F{M : \<(x^{L_k}_k:U_k)_n\v U\>}{M^{+j} : \<
    (x^{j::L_k}_k:\overline{e}_jU_k)_n \v \overline{e}_jU\>}$ and $\forall~1 \leq k \leq n$,
  $N_k \in {\cal I}(\overline{e}_jU_k) = {\cal I}(U_k)^{+j}$. Then $\forall~1
  \leq k \leq n$, $N_k = P_k^{+j}$ where $P_k \in {\cal I}(U_k)$.  By
  lemma~\ref{interpret-intsub}.\ref{interprettwo}, for all $k \in \{1,
  \dots, n\}$, $P_k \in {\cal M}^{L_k}$.  By the definition of the
  substitution, $\diamond \{M^{+j}\} \cup \{N_k$ / $k \in \{1, \dots,
  n\}\}$. By lemma~\ref{deg+-f+}.\ref{one'''}, $\diamond \{M\} \cup
  \{P_k$ / $k \in \{1, \dots, n\}\}$. By
  lemma~\ref{degsub}.\ref{degsubthree}, $M[(x^{L_k}_k:=P_k)_n] \in
  {\cal M}$.  By IH, $M[(x^{L_k}_k:=P_k)_n] \in {\cal I}(T)$. Hence,
  by lemma \ref{deg+-f+}, $M^{+j}[(x^{j::L_k}_k:=N_k)_n] =
  (M[(x^{L_k}_k:=P_k)_n])^{+j} \in {\cal I}(U)^{+j} = {\cal I}(\overline{e}_jU)$.

\item Let $\F{M: \Phi \;\;\; \Phi \sqsubseteq \Phi'} {M : \Phi'}$
  where $\Phi' = \<(x^{L_j}_j:U_j)_n \v U\>$. By lemma
  \ref{env-Phisub}, we have $\Phi = \<(x^{L_j}_j : U'_j)_n \v U'\>$,
  where for every $1 \leq j \leq n$, $U_j \sqsubseteq U'_j$ and $U'
  \sqsubseteq U$. By lemma~\ref{interpret-intsub}.\ref{intsub}, $N_j
  \in {\cal I}(U'_j)$, then, by IH, $M[(x^{L_j}_j:=N_j)_n] \in {\cal
    I}(U')$ and, by lemma \ref{interpret-intsub}.\ref{intsub},
  $M[(x^{L_j}_j:=N_j)_n] \in {\cal I}(U)$.
\hfill $\Box$
\end{itemize}
\end{proof}

\begin{proof}[Of lemma~\ref{exlem}]
\begin{enumerate}
\item Let $y \in {\cal V}_2$ and ${\cal X} = \{M \in {\cal M}^\oslash$
  / $M \rhd_\be^* x^\oslash N_1...N_k$ where $k\geq 0$ and $x \in
  {\cal V}_1$ or $M \rhd_\be^* y^\oslash\}$. ${\cal X}$ is
  $\be$-saturated and $\forall x \in {\cal V}_1, {\cal N}_x^\oslash
  \subseteq {\cal X} \subseteq {\cal M}^\oslash$. Take a
  $\be$-interpretation ${\cal I}$ such that ${\cal I}(a) = {\cal
    X}$. If $M \in [Id_0]_\be$, then $M$ is closed and $M \in {\cal X}
  \fx {\cal X}$. Since $y^\oslash \in {\cal X}$ and $M \diamond
  y^\oslash$ then $M y^\oslash \in {\cal X}$ and $M y^\oslash
  \rhd_\be^* x^\oslash N_1...N_k$ where $k\geq 0$ and $x \in {\cal
    V}_1$ or $M y^\oslash \rhd_\be^* y^\oslash$. Since $M$ is closed
  and $x^\oslash \not = y^\oslash$, by lemma~\ref{deg=}.\ref{deg=one},
  $My^\oslash \rhd_\be^* y^\oslash$. Hence, by
  lemma~\ref{trivial1}.\ref{trivial1three}, $M \rhd_\be^* \l
  y^\oslash.y^\oslash$ and, by lemma~\ref{deg=}, $M \in {\cal
    M}^\oslash$.

  Conversely, let $M \in {\cal M}^\oslash$ such that $M$ is closed and
  $M \rhd_\be^* \l y^\oslash.y^\oslash$. Let ${\cal I}$ be an
  $\be$-interpretation and $N \in {\cal I}(a)$ such that $M \diamond
  N$. By lemma~\ref{interpret-intsub}.\ref{interprettwo}, $N \in {\cal
    M}^{\oslash}$, so $MN \in {\cal M}^{\oslash}$.  Since ${\cal
    I}(a)$ is $\be$-saturated and $MN \rhd_\be^* N$, $MN \in {\cal
    I}(a)$ and hence $M \in {\cal I}(a) \fx {\cal I}(a)$. Hence, $M
  \in [Id_0]_\be$.

\item By lemma~\ref{[]} and lemma~\ref{fx+}, $[Id'_1]_\be = [\overline{e}_1a \f
  \overline{e}_1a]_\be = [\overline{e}_1(a \f a)]_\be = [Id_1] = [a \f a]_\be^{+1} =
  [Id_0]_\be^{+1}$.  By 1., $[Id_{0}]_\be^{+1} = \{M \in {\cal
    M}^{(1)}$ / $M \rhd_\be^* \l y^{(1)}.y^{(1)}\}$.

\item Let $y \in {\cal V}_2$, ${\cal X} = \{M \in {\cal M}^\oslash$ /
  $M \rhd_\be^* y^\oslash$ or $M \rhd_\be^* x^\oslash N_1...N_k$ where
  $k \geq 0$ and $x \in {\cal V}_1\}$ and ${\cal Y} = \{M \in {\cal
    M}^\oslash$ / $M \rhd_\be^* y^\oslash y^\oslash$ or $M \rhd_\be^*
  x^\oslash N_1...N_k$ or $M \rhd_\be^* y^\oslash (x^\oslash
  N_1...N_k)$ where $k \geq 0$ and $x \in {\cal V}_1\}$. ${\cal X}$,
  ${\cal Y}$ are $\be$-saturated and $\forall x \in {\cal V}_1, {\cal
    N}_x^\oslash \subseteq {\cal X},{\cal Y} \subseteq {\cal
    M}^\oslash$. Let ${\cal I}$ be a $\be$-interpretation such that
  ${\cal I}(a) = {\cal X}$ and ${\cal I}(b) = {\cal Y}$. If $M \in
  [D]_\be$, then $M$ is closed (hence $M \diamond y^\oslash$) and $M
  \in ({\cal X} \cap ({\cal X} \fx {\cal Y})) \fx {\cal Y}$. Since
  $y^\oslash \in {\cal X}$ and $y^\oslash \in {\cal X} \fx {\cal Y}$,
  $y^\oslash \in {\cal X} \cap ({\cal X} \fx {\cal Y})$ and $M
  y^\oslash \in {\cal Y}$. Since $x^\oslash \not = y^\oslash$, by
  lemma~\ref{deg=}.\ref{deg=one}, $My^\oslash \rhd_\be^* y^\oslash
  y^\oslash$. Hence, by lemma~\ref{trivial1}.\ref{trivial1three}, $M
  \rhd_\be^* \l y^\oslash.y^\oslash y^\oslash$ and, by
  lemma~\ref{deg=}, $\deg(M) = \oslash$ and $M \in {\cal M}^\oslash$.

  Conversely, let $M \in {\cal M}^\oslash$ such that $M$ is closed and
  $M \rhd_\be^* \l y^\oslash.y^\oslash y^\oslash$. Let ${\cal I}$ be
  an $\be$-interpretation and $N \in {\cal I}(a \sqcap (a \f b)) =
  {\cal I}(a) \cap ({\cal I}(a) \fx {\cal I}(b))$ such that $M
  \diamond N$.  By lemma~\ref{interpret-intsub}.\ref{interprettwo} and
  lemma~\ref{degsub}.\ref{lem:refl+symm}, $N \in {\cal M}^{\oslash}$
  and $N \diamond N$.  So $NN, MN \in {\cal M}^{\oslash}$.  Since
  ${\cal I}(b)$ is $\be$-saturated, $NN \in {\cal I}(b)$ and $MN
  \rhd_\be^* NN$, we have $MN \in {\cal I}(b)$ and hence $M \in {\cal
    I}(a \sqcap (a \f b)) \fx {\cal I}(b)$.  Therefore, $M \in
  [D]_\be$.

\item Let $f, y \in {\cal V}_2$ and take ${\cal X} = \{M \in {\cal
    M}^\oslash$ / $M \rhd_\be^* (f^\oslash)^n(x^\oslash N_1...N_k)$ or
  $M \rhd_\be^* (f^\oslash)^n y^\oslash$ where $k, n\geq 0$ and $x \in
  {\cal V}_1\}$. ${\cal X}$ is $\be$-saturated and $\forall x \in
  {\cal V}_1, {\cal N}_x^\oslash \subseteq {\cal X} \subseteq {\cal
    M}^\oslash$. Let ${\cal I}$ be a $\be$-interpretation such that
  ${\cal I}(a) = {\cal X}$. If $M \in [Nat_{0}]_\be$, then $M$ is
  closed and $M \in ({\cal X} \fx {\cal X}) \fx ({\cal X} \fx {\cal
    X})$. We have $f^\oslash \in {\cal X} \fx {\cal X}$, $y^\oslash
  \in {\cal X}$ and $\diamond \{M, f^\oslash, y^\oslash\}$ then $M
  f^\oslash y^\oslash\in {\cal X}$ and $M f^\oslash \, y^\oslash
  \rhd_\be^* (f^\oslash)^n (x^\oslash N_1...N_k)$ or $M f^\oslash
  y^\oslash \rhd_\be^* (f^\oslash)^n y^\oslash$ where $n \geq 0$ and
  $x \in {\cal V}_1$. Since $M$ is closed and $\{x^\oslash\} \cap
  \{y^\oslash, f^\oslash\} = \emptyset$, by
  lemma~\ref{deg=}.\ref{deg=one}, $Mf^\oslash y^\oslash \rhd_\be^*
  (f^\oslash)^n y^\oslash$ where $n \geq 1$. Hence, by
  lemma~\ref{trivial1}.\ref{trivial1three}, $M \rhd_\be^* \l
  f^\oslash.f^\oslash$ or $M \rhd_\be^* \l f^\oslash.\l
  y^\oslash.(f^\oslash)^n y^\oslash$ where $n \geq 1$. Moreover, by
  lemma~\ref{deg=}, $\deg(M) = \oslash$ and $M \in {\cal M}^\oslash$.

  Conversely, let $M \in {\cal M}^\oslash$ such that $M$ is closed and
  $M \rhd_\be^* \l f^\oslash.f^\oslash$ or $M \rhd_\be^* \l f^\oslash
  .\l y^\oslash . (f^\oslash)^n y^\oslash$ where $n \geq 1$. Let
  ${\cal I}$ be an $\be$-interpretation, $N \in {\cal I}(a \f a) =
  {\cal I}(a) \fx {\cal I}(a)$ and $N' \in {\cal I}(a)$ such that
  $\diamond \{M, N, N'\}$.  By
  lemma~\ref{interpret-intsub}.\ref{interprettwo}, $N, N' \in {\cal
    M}^{\oslash}$, so $MNN', (N)^m N' \in {\cal M}^{\oslash}$, where
  $m \geq 0$.  We show, by induction on $m \geq 0$, that $(N)^m N'\in
  {\cal I}(a)$. Since $MNN' \rhd_\be^* (N)^m N'$ where $m \geq 0$ and
  $(N)^m N'\in {\cal I}(a)$ which is $\be$-saturated, then $MNN' \in
  {\cal I}(a)$.  Hence, $M \in ({\cal I}(a) \fx {\cal I}(a)) \f ({\cal
    I}(a) \fx {\cal I}(a))$ and $M \in [Nat_{0}]_\be$.

\item By lemma \ref{[]}, $[Nat_{1}] = [\overline{e}_1Nat_{0}] =
  [Nat_{0}]^{+1}$.  By \ref{exfour}., $[Nat_1] = [Nat_0]^{+1} = \{M
  \in {\cal M}^{(1)}$ / $M \rhd_\be^*\l f^{(1)}.f^{(1)}$ or $M
  \rhd_\be^* \l f^{(1)}.\l y^{(1)}.(f^{(1)})^n y^{(1)}$ where $n \geq
  1\}$.

\item Let $f, y \in {\cal V}_2$ and take ${\cal X} = \{M \in {\cal
    M}^\oslash$ / $M \rhd_\be^* x^\oslash P_1...P_l$ or $M \rhd_\be^*
  f^\oslash(x^\oslash Q_1...Q_n)$ or $M \rhd_\be^* y^\oslash$ or $M
  \rhd_\be^* f^\oslash y^{(1)}$ where $l, n \geq 0$ and $\deg(Q_i)
  \succeq {(1)}\}$. ${\cal X}$ is $\be$-saturated and $\forall x \in
  {\cal V}_1, {\cal N}_x^\oslash \subseteq {\cal X} \subseteq {\cal
    M}^\oslash$. Let ${\cal I}$ be a $\be$-interpretation such that
  ${\cal I}(a) = {\cal X}$. If $M \in [Nat'_0]_\be$, then $M$ is
  closed and $M \in ({\cal X}^{+1} \fx {\cal X}) \fx ({\cal X}^{+1}
  \fx {\cal X})$. Let $N \in {\cal X}^{+1}$ such that $N \diamond
  f^\oslash$. We have $N \rhd_\be^* x^{(1)} P^{+1}_1...P^{+1}_k$ or $N
  \rhd_\be^* y^{(1)}$, then $f^\oslash N \rhd_\be^* f^\oslash
  (x^{(1)} P^{+1}_1...P^{+1}_k) \in {\cal X}$ or $N \rhd_\be^*
  f^\oslash y^{(1)} \in {\cal X}$, thus $f^\oslash \in {\cal X}^{+1}
  \fx {\cal X}$. We have $f^\oslash \in {\cal X}^{+1} \fx {\cal X}$,
  $y^{(1)} \in {\cal X}^{+1}$ and $\diamond\{M,f^\oslash,y^{(1)}\}$,
  then $M f^\oslash y^{(1)} \in {\cal X}$. Since $M$ is closed and
  $\{x^\oslash, x^{(1)}\} \cap \{y^{(1)}, f^\oslash\} = \emptyset$, by
  lemma~\ref{deg=}.\ref{deg=one}, $Mf^\oslash y^{(1)} \rhd_\be^*
  f^\oslash y^{(1)}$. Hence, by
  lemma~\ref{trivial1}.\ref{trivial1three}, $M \rhd_\be^* \l
  f^\oslash.f^\oslash$ or $M \rhd_\be^* \l f^\oslash.\l
  y^{(1)}.f^\oslash y^{(1)}$. Moreover, by lemma~\ref{deg=}, $\deg(M)
  = \oslash$ and $M \in {\cal M}^\oslash$.

  Conversely, let $M \in {\cal M}^\oslash$ such $M$ is closed and $M
  \rhd_\be^* \l f^\oslash.f^\oslash$ or $M \rhd_\be^* \l f^\oslash .\l
  y^{(1)}. f^\oslash y^{(1)}$. Let ${\cal I}$ be an
  $\be$-interpretation, $N \in {\cal I}(\overline{e}_1a \f a) = {\cal I}(a)^{+1}
  \fx {\cal I}(a)$ and $N' \in {\cal I}(a)^{+1}$ where $\diamond \{M,
  N, N'\}$.  By lemma~\ref{interpret-intsub}.\ref{interprettwo}, $N
  \in {\cal M}^{\oslash}$ and $N' \in {\cal M}^{(1)}$, so $MNN', NN'
  \in {\cal M}^{\oslash}$.  Since $MNN' \rhd_\be^* NN'$, $NN'\in {\cal
    I}(a)$ and ${\cal I}(a)$ is $\be$-saturated, then $MNN' \in {\cal
    I}(a)$. Hence, $M \in ({\cal I}(a)^{+1} \fx {\cal I}(a)) \f ({\cal
    I}(a)^{+1} \fx {\cal I}(a))$ and $M \in [Nat'_0]$.
\hfill $\Box$
\end{enumerate}
\end{proof}

\fi

\end{document}

%% file: jbw-switch.tex
\makeatletter
  \input{jbw-switch.sty}
\makeatother

%% file: compsem-big.bbl
\newcommand{\bibnoop}[1]{} \newcommand{\bibvonmagic}[2]{#2}
  \newcommand{\bibsingleletter}[1]{#1}
\begin{thebibliography}{10}
\expandafter\ifx\csname href\endcsname\relax
  \def\href#1#2{#2}\fi

\bibitem{Barendregt:LCSS-1984}
H.~P. Barendregt.
\newblock {\em The Lambda Calculus: Its Syntax and Semantics}.
\newblock North-Holland, revised edition, 1984.

\bibitem{Car+Pol+Wel+Kfo:ESOP-2004}
S.~Carlier, J.~Polakow, J.~B. Wells, A.~J. Kfoury.
\newblock
  \href{http://www.macs.hw.ac.uk/~jbw/papers/Carlier+Polakow+Wells+Kfoury:Syst%
em-E:ESOP-2004.pdf}{{S}ystem {E}: Expansion variables for flexible typing with
  linear and non-linear types and intersection types}.
\newblock In {\em Programming Languages \& Systems, 13th European Symp.\
  Programming}, vol. 2986 of {\em LNCS}. Springer-Verlag, 2004.

\bibitem{Car+Wel:ITRS-2004}
S.~Carlier, J.~B. Wells.
\newblock
  \href{http://www.macs.hw.ac.uk/~jbw/papers/Carlier+Wells:Expansion:ITRS-2004%
.pdf}{Expansion: the crucial mechanism for type inference with intersection
  types: {A} survey and explanation}.
\newblock In {\em Proc.\ 3rd Int'l Workshop Intersection Types \& Related
  Systems (ITRS 2004)}, 2005.
\newblock The ITRS '04 proceedings appears as vol.\ 136 (2005-07-19) of
  \emph{Elec.\ Notes in Theoret.\ Comp.\ Sci.}

\bibitem{Cop+Dez-Cia+Ven:HBC-1980}
M.~Coppo, M.~Dezani-Ciancaglini, B.~Venneri.
\newblock Principal type schemes and $\lambda$-calculus semantics.
\newblock In J.~R. Hindley, J.~P. Seldin, eds., {\em To {H}.~{B}.\ {C}urry:
  Essays on Combinatory Logic, Lambda Calculus, and Formalism}. Academic Press,
  1980.

\bibitem{Coquand:TLCA-2005}
T.~Coquand.
\newblock Completeness theorems and lambda-calculus.
\newblock In P.~Urzyczyn, ed., {\em TLCA}, vol. 3461 of {\em Lecture Notes in
  Computer Science}. Springer, 2005.

\bibitem{LNCS/1875}
G.~Goos, J.~Hartmanis, eds.
\newblock {\em $\l$-Calculus and Computer Science Theory, Proceedings of the
  Symposium Held in Rome, March 15-27, 1975}, vol.~37 of {\em Lecture Notes in
  Computer Science}. Springer-Verlag, 1975.

\bibitem{Hindley:ISOP-1982}
J.~R. Hindley.
\newblock The simple semantics for {C}oppo-{D}ezani-{S}all{\'e} types.
\newblock In M.~Dezani-Ciancaglini, U.~Montanari, eds., {\em International
  Symposium on Programming, 5th Colloquium}, vol. 137 of {\em LNCS}, Turin,
  1982. Springer-Verlag.

\bibitem{Hin2}
J.~R. Hindley.
\newblock The completeness theorem for typing $\l$-terms.
\newblock {\em Theoretical Computer Science}, 22, 1983.

\bibitem{Hin3}
J.~R. Hindley.
\newblock Curry's types are complete with respect to {F}-semantics too.
\newblock {\em Theoretical Computer Science}, 22, 1983.

\bibitem{Hindley:BSTT-1997}
J.~R. Hindley.
\newblock {\em Basic Simple Type Theory}, vol.~42 of {\em Cambridge Tracts in
  Theoretical Computer Science}.
\newblock Cambridge University Press, 1997.

\bibitem{kamnour}
F.~Kamareddine, K.~Nour.
\newblock A completeness result for a realisability semantics for an
  intersection type system.
\newblock {\em Ann. Pure Appl. Logic}, 146(2-3), 2007.

\bibitem{kamnourrahliwells}
F.~Kamareddine, K.~Nour, V.~Rahli, J.~B. Wells.
\newblock A complete realisability semantics for intersection types and
  infinite expansion variables.
\newblock Located at
  \url{http://www.macs.hw.ac.uk/~fairouz/papers/drafts/compsem-big.pdf}, 2008.

\bibitem{report}
F.~Kamareddine, K.~Nour, V.~Rahli, J.~B. Wells.
\newblock Realisability semantics for intersection type systems and expansion
  variables.
\newblock In ITRS'08. The file is Located at
  \url{http://www.macs.hw.ac.uk/~fairouz/papers/
  conference-publications/semone.pdf}, 2008.

\bibitem{krivine-lctm-1990}
J.~Krivine.
\newblock {\em Lambda-Calcul~: Types et Mod\`eles}.
\newblock Etudes et Recherches en Informatique. Masson, 1990.

\end{thebibliography}
